\documentclass[12pt]{amsart}

\usepackage{mathrsfs}
\usepackage{amsmath,amssymb,enumerate,tikz,pdfpages}
\usepackage{amsfonts,amssymb}
\usepackage{booktabs}
\usepackage{comment}
\usepackage[colorlinks,linkcolor=red,anchorcolor=blue,
citecolor=green,CJKbookmarks=True]{hyperref}
\usepackage[T1]{fontenc}
\usepackage{csquotes}

\usepackage{soul, xcolor, todonotes}

\newif\ifistoreview
\istoreviewtrue

\newcommand{\setreviewson}{\istoreviewtrue}

\setreviewson 

\theoremstyle{plain}
\newtheorem{theorem}{Theorem}[section]
\newtheorem{lemma}[theorem]{Lemma}
\newtheorem{corollary}[theorem]{Corollary}
\newtheorem{definition}[theorem]{Definition}

\newtheorem{proposition}[theorem]{Proposition}

\theoremstyle{definition}

\newtheorem{remark}[theorem]{Remark}
\newtheorem{question}[theorem]{Question}

\newtheorem{construction}[theorem]{Construction}

\newcommand{\tec}{Teichm\"uller }

\newcommand{\dth}{d_{Th}}

\newcommand{\HSR}{\mathrm{HSL}^+}
\newcommand{\PHSL}{\mathrm{PHSL}}
\newcommand{\HSL}{\mathrm{HSL}}
\newcommand{\HR}{\mathrm{HR}}
\newcommand{\hs}{}

\newcommand{\env}{\mathbf{Env}}

\newcommand{\outenv}{\mathbf{Out}}
\newcommand{\inenv}{\mathbf{In}}
\newcommand{\crl}{\mathcal{CRL}}

\newcommand{\far}{\mathrm{far}}
\newcommand{\cut}{\mathrm{cut}}
\newcommand{\cone}{\mathrm{Cone}}

\newcommand{\supp}{\mathrm{supp}}

\newcommand{\rpi}{\pi_r}
\newcommand{\lpi}{\pi_l}
\newcommand{\lomega}{\mathrm{Lbd}}
\newcommand{\romega}{\mathrm{Rbd}}
\newcommand{\lbd}{\mathrm{Lbd}}
\newcommand{\rbd}{\mathrm{Rbd}}
\newcommand{\bnd}{\mathrm{Bd}}

\newcommand{\Hopf}{\mathrm{Hopf}}

\newcommand{\ML}{\mathcal{ML}}
\newcommand{\MF}{\mathcal{MF}}

\newcommand{\PMF}{\mathcal{PMF}}

\newcommand{\str}{\mathrm{St}}

\newcommand{\T}{\mathcal{T}}
\newcommand{\R}{\mathbb R}

\renewcommand{\subset}{\subseteq}
\renewcommand{\supset}{\supseteq}

\title[Envelopes of  the Thurston metric]{Envelopes of  the Thurston metric  on  Teichm\"uller space}
\author{Huiping Pan}
\address{Huiping Pan,
School of Mathematics, South China  University of Technology, 510641, Guangzhou, China}
\email{panhp@scut.edu.cn} 

\author{Michael Wolf}
\address{Michael Wolf, School of Mathematics,
Georgia Institute of Technology,
Atlanta, GA USA 30332}
\email{mwolf40@gatech.edu}

\date{August 12, 2025}

\begin{document}

\begin{abstract}
  For the Thurston (asymmetric) metric on  Teichm\"uller space, the defect from being uniquely geodesic is described by the envelope, defined as the union of geodesics from the initial point to the terminal point. 
  
  Using the harmonic stretch lines we defined recently, we describe the shape of envelopes as a cone over a cone over a space, defined from a topological invariant of the initial and terminal points. In addition, we show that the envelope is always contractible. We prove that envelopes vary continuously with their endpoints. We also provide a parametrization of out-envelopes and in-envelopes in terms of straightened measured laminations complementary to the prescribed maximally stretched laminations.  
  
  We extend most of these results to the metrically infinite envelopes which have a terminal point on the Thurston boundary, illustrating some of the nuances of these with examples, and describing the accumulation set. Finally, we develop a new characterization of harmonic stretch lines that avoids a limiting process.
  \end{abstract}

\maketitle

\tableofcontents

 \section{Introduction}

\subsection{Overview} 
Aiming for a geometric understanding of the Teichm\"uller space from the point of view of hyperbolic geometry, Thurston introduced an asymmetric Finslerian metric, the Thurston metric, on the Teichm\"uller space and  proved that this metric is geodesic but not uniquely geodesic \cite{Thurston1986}.
 
The defect of an ordered pair of points from having a unique length-minimizing geodesic from the first to the second is described by the \emph{envelope}, which is defined as the union of geodesics from the initial point to the terminal point. For the one-dimensional Teichm\"uller space of the once-punctured torus, Dumas-Lenzhen-Rafi-Tao proved that the envelope is either a geodesic segment or a geodesic quadrilateral with the given initial and terminal points as opposite vertices, and that it varies continuously with its endpoints \cite[Theorem 1.1]{DLRT2020}.

The goal of this paper is to study envelopes for surfaces of higher complexity.
 We show that, for the Teichm\"uller space of closed orientable surfaces of genus at least two, an envelope is a cone over a cone over an ancestor space, and we describe that ancestor space in terms of the initial and terminal points of the envelope.  Moreover, we also show that these envelopes vary \emph{continuously} with their endpoints. Envelopes were found by \cite{DLRT2020} to be intersections of the closures of \emph{out-envelopes} and \emph{in-envelopes} -- spaces of geodesics with prescribed initial or terminal points, respectively -- and we also provide a parametrization of  those spaces in terms of \enquote{straightened measured laminations}. we show that an envelope is the closure of the intersection of the out-envelope and the in-envelope (improving a result of \cite{DLRT2020}).

 We extend most of these results to envelopes with one endpoint on the Thurston boundary.
 
 Our methods are different from those in \cite{DLRT2020}.  In particular, while those authors relied heavily upon the one-dimensional nature of the \tec spaces they were studying along with some sophisticated estimates on hyperbolic geometric quantities, we build upon analytical techniques we developed in \cite{PanWolf2022}.  In that article, we showed that Thurston geodesics could be approximated by harmonic map rays.  Moreover, certain uniqueness results--  tracing their roots back to the Jenkins-Serrin \cite{JS66} results on graphs with infinite boundary values -- and consequent uniform convergence  theorems allowed for us to pick out canonical Thurston geodesics $[X,Y]$ from $X$ to $Y$ in \tec space.  These canonical geodesics give rise to canonical foliations of the envelopes (as well as the out- and in-envelopes), and it is from these foliations that we are able to deduce our structure theorems.

 \subsection{Main Results.}
In this subsection, we state our results more precisely.
 Let $S$ be a closed  orientable surface of genus at least two. The Teichm\"uller space $\T(S)$ of $S$ is defined to be the equivalence classes of hyperbolic metrics on $S$, where two hyperbolic metrics, say $X=(S,g_X)$ and $Y=(S,g_{Y})$, are \emph{equivalent} if there exists an isometry $\xi: (S,g_X)\to (S,g_Y)$ between them which is isotopic to the identity map $\mathrm{id}:S\to S$. In the sequel, we will only refer to points in $\T(S)$ as $X$ or $Y$, omitting the metric description $g_X$ or $g_Y$.  
 
 For any two hyperbolic surfaces $X$ and $Y$ in $\T(S)$, the Thurston metric $d_{Th}$ is defined as:
\begin{equation*}
	d_{Th}(X,Y):=\log \inf_f \{L(f)  \}
\end{equation*}
where $L(f)$ is the Lipschitz constant of $f$ and where $f$ ranges over all Lipschitz homeomorphisms from $X$ to $Y$ in the homotopy class of the identity map $\mathrm{id}:S\to S$.  As we mentioned earlier, the Thurston metric is geodesic but not uniquely geodesic.  The \emph{envelope} $\env(X,Y)$ from $X$ to $Y$ is defined to be the union of \emph{all} geodesics from $X$ to $Y$ in the Thurston metric.

\subsubsection{Shape and continuity of Envelopes.}
Our first main result concerns the continuity of envelopes. 
  		  	 	\begin{theorem}[Continuity] \label{thm:env:top}
 Suppose $S$ is a closed orientable surface of genus at least two. The envelope $\env(X,Y)$ varies continuously {with} $X$ and $Y$ in $\T(S)$.
 \end{theorem}
(Of course, in the classical case of genus one, the Thurston metric is the (uniquely geodesic) hyperbolic metric 
  \cite{GJ2021,Saglam2021}, which being negatively curved has envelopes (i.e. in this case, geodesic arcs) which vary continuously from general considerations. For the punctured torus case, see \cite{DLRT2020}).

 Next, we shall describe the shapes of envelopes in more detail. Recall that the \emph{maximally stretched lamination} from $X$ to $Y$ in $\T(S)$, denoted by $\Lambda(X,Y)$,  is the union of all chain-recurrent laminations to which the restriction of every optimal Lipschitz map from $X$ to $Y$ in the homotopy class of the identity map takes leaves of $\lambda$ on $X$ to corresponding leaves of $\lambda$ on $Y$, multiplying arclength by a factor of $\exp(d_{Th}(X,Y))$.  Based on maximally stretched laminations, we define the \emph{out-envelopes} and \emph{in-envelopes} as:
 \begin{equation*}
 	 	\outenv(X,\lambda)=\{Z: \Lambda(X,Z)=\lambda\},\quad
 	\inenv(Y,\lambda)=\{Z:\Lambda(Z,Y)=\lambda\}.
 \end{equation*}

\begin{remark} \label{rem: append endpoint to Out and In}
    We will often study not only $\outenv(X,Y)$ but also $\{X\} \cup \outenv(X,Y)$, where the topology on this union is the one induced from the inclusion of both into $\T(S)$. Similar considerations also apply to $\{Y\} \cup \inenv(Y, \lambda)$. 
\end{remark}

 		  	 	Recall that for any two distinct points $X$ and $Z$ in $ \T(S)$, there exists a unique harmonic stretch line $\HSL(X,Z)$ proceeding from $X$ through $Z$ \cite[Theorem 1.6]{PanWolf2022}. Let $[X,Z]$ be the geodesic segment of $\HSL(X,Z)$ with endpoints $X$ and $Z$.    
 		  	 	
 		  	 	For any $X$ and $Y$ in $\T(S)$, the point $Z\in\env(X,Y)$ is called a \emph{right boundary point} from $X$ to $Y$ if 
 		  	 \begin{equation*}
 		  	 	\HSL(X,Z)\cap \env(X,Y)=[X,Z].
 		  	 \end{equation*}
 		  	 Similarly one can define \emph{left boundary points}.  A point is called a \emph{boundary point} if it is both a left boundary point and a right boundary point. We refer to Lemma \ref{lem:extendability} for a characterization of left/right boundary points in terms of maximally stretched laminations. 
 		  	 Let $\lbd(X,Y)$ and $\rbd(X,Y)$ be respectively the set of left boundary and the set of right boundary points from $X$ to $Y$.  Let $\bnd(X,Y)=\lbd(X,Y) \cap \rbd(X,Y)$ be the set of boundary points from $X$ to $Y$. 
 		  	 
 		  	 With these notions at hand, we now state our first main result on the geometric topological structure of an envelope. 
 		  	 	
 \begin{theorem}[Shapes] \label{thm:env:shape}
 Suppose $S$ is a closed orientable surface of genus at least two. Let  $X,Y\in\T(S)$ be distinct points in Teichm\"uller space. Then
 	\begin{enumerate}
 		\item $\env(X,Y)$ is a  geodesic segment if and only if the maximally stretched lamination $\Lambda(X,Y)$ is maximal as a chain-recurrent lamination.  
       \item $ \env(X,Y)=\overline{\outenv(X,\Lambda(X,Y))\cap\inenv(Y,\Lambda(X,Y))}$.
       \item   $\env(X,Y)$ is a cone over $\rbd(X,Y)$ as well as a cone over $\lbd(X,Y)$.
       \item Both $\rbd(X,Y)$ and $\lbd(X,Y)$ are cones over $\bnd(X,Y)$.
       \item The envelope $\env(X,Y)$ is contractible.
 	\end{enumerate}
 \end{theorem}

We remark that, in contrast to (ii) above, it was shown in \cite{DLRT2020} that $ \env(X,Y)=\overline{\outenv(X,\Lambda(X,Y))}\cap\overline{\inenv(Y,\Lambda(X,Y))}$ (cf. Proposition~\ref{prop:DLRT2020}).

[Also, here we understand the fourth statement to include the degenerate case where the envelope is a line segment, $\bnd(X,Y)$ is empty, and both $\rbd(X,Y) = \{Y\}$ and $\lbd(X,Y) = \{X\}$ are singletons.]

It may be worth observing, in relation to the fifth statement, that in the Riemannian setting, we have the following: (i) when the sectional curvature of a metric is non-positive, the envelope is a segment, (ii) there are spaces (e.g. a smoothing of a plane with a disk replaced by a hemispherical cap) with some positive sectional curvature where the envelope is contractible but may not be a segment, and (iii) there are spaces with regions of large positive sectional curvature, for example an ellipsoid with large long-axial ratio, where some envelopes are not contractible.

At this point, it may be useful to reflect on how this result might extend the Dumas-Lenzhen-Rafi-Tao \cite{DLRT2020} result on envelopes for the \tec space of the once-punctured torus: in that case and when the maximal stretch lamination was a simple closed curve, the envelope was a (non-trivial) quadrilateral in the plane, with opposite vertices given by $X$ and $Y$, as well as two additional points. Now, the paper \cite{PanWolf2022} only treated the case of closed surfaces, so the results of that paper, and hence this one, do not strictly apply to the situation of the punctured torus case treated in \cite{DLRT2020}, but we can suggest how Theorem~\ref{thm:env:shape} might be reflected in the once-punctured torus case. In that direction, the sides of the quadrilateral incident to $X$ are $\lbd(X,Y)$ and the sides of the quadrilateral incident to $Y$ are $\rbd(X,Y)$. Geodesics connecting $X$ and the sides incident to $Y$ are given by harmonic map stretch segments, as are those connecting $Y$ to the sides incident to $X$; in either case, those collections of geodesics foliate the envelope. Theorem~\ref{thm:env:shape} and its proof provide this structure for \tec spaces of compact surfaces.

\subsubsection{A topological description of the ancestor of the envelope.} We next seek to clarify the nature of these basic extremal points $\bnd(X,Y)$, from which the \enquote{cone on cone} is built (cf. Theorem \ref{thm:env:shape}(iii) and (iv)). Those results say that each of $\rbd(X,Y)$ and $\lbd(X,Y)$ is a cone over $\bnd(X,Y)$, and also that $\env(X,Y)$ is a cone over either $\rbd(X,Y)$ or $\lbd(X,Y)$: this justifies the phrase \enquote{cone on cone} for $\env(X,Y)$. 

Now, the basic topological quantity associated to an envelope is the maximally stretched lamination $\Lambda(X,Y)$, and so we aim for an answer that involves only that object.

 Moving towards that topological characterization of $\bnd(X,Y)$, our next result parametrizes out-envelopes and in-envelopes using \emph{$\lambda$-straightened laminations}. Let $X\in\T(S)$.  For any chain recurrent geodesic lamination $\lambda$,  let $\cone(\lambda)$ be the set of $\lambda$-straightened measured laminations in $X\setminus\lambda$ (see Definition \ref{def:ad:lamination}. Define $$\Theta^{out}(\lambda)=\cone(\lambda)\times[0,\infty)/_\sim$$
  where $(\gamma,t)\sim(\gamma',t')$ if $t=t'=0$.  Since there is a bijection between geodesic laminations on any two homeomorphic hyperbolic surfaces, we see that  $\Theta^{out}(\lambda)$ is independent of the choice of $X$. 
  
   Recall that any harmonic stretch ray starting at $X$ is constructed by a surjective harmonic diffeomorphism from some punctured surface to $X\backslash\lambda$ for some chain recurrent geodesic lamination $\lambda$. The horizontal foliation of the Hopf differential of this map gives rise to a $\lambda$-straightened measured lamination on $X\backslash\lambda$ (see Section \ref{subsec:hsr:construction} for more details). This defines a map
  $$\Psi_{out}:\{X\}\cup\outenv(X,\lambda)\to\Theta^{out}(\lambda)$$ which sends $Z\in\{X\}\cup \outenv(X,\lambda)$ to the pair $(\gamma,\dth(X,Z))$, where $\gamma$ is the $\lambda$-straightened measured lamination defining the harmonic stretch ray  which starts at $X$ and which passes through $Z$. For $Y\in\T(S)$, define $\Theta^{in}(\lambda)=\cone(\lambda)\times (-\infty,0]/_\sim$ and $\Psi_{in}:\{Y\}\cup\inenv(Y,\lambda) \to \Theta^{in}(\lambda)$ similarly.
  Let $\crl(S)$ be the space of chain recurrent geodesic laminations on $S$, equipped with the Hausdorff topology. 

  Our next main results on the structure of the envelope are Theorem~\ref{thm:structure:out:in} and Corollary~\ref{cor:bnd:description}.
  
\begin{theorem}[out-envelopes and in-envelopes]
    \label{thm:structure:out:in}
 Suppose $S$ is a closed orientable surface of genus at least two. 
 \begin{enumerate}
     \item For any $X\in\T(S)$ and any chain recurrent geodesic lamination $\lambda$, {the map} $\Psi_{out}:\{X\}\cup\outenv(X,\lambda)\to\Theta^{out}(\lambda)$, which sends geodesic rays in $\{X\}\cup\outenv(X,\lambda)$ starting at $X$ to rays in $\Theta^{out}(\lambda)$ starting at the cone point, is a homeomorphism.
     \item $\overline{\outenv(X,\lambda)}$ varies continuously in $(X,\lambda)\in\T(S)\times\crl(S)$, i.e. $\overline{\outenv(X_n,\lambda_n)}$ converges to $\overline{\outenv(X,\lambda)}$ locally uniformly as $X_n\to X$ in $\T(S)$ and $\lambda_n\to\lambda$ in $\crl(S)$ with the Hausdorff topology. 
 \end{enumerate}
 Similar results also hold for in-envelopes: 
 \begin{enumerate}
     \item[(i')] For any $Y\in\T(S)$ and any chain recurrent geodesic lamination $\lambda$, {the map} $\Psi_{in}:\{Y\}\cup\inenv(Y,\lambda)\to\Theta^{in}(\lambda)$,  which sends geodesic rays in $\{Y\}\cup\inenv(Y,\lambda)$ ending at $Y$ to rays in $\Theta^{in}(\lambda)$ ending at the cone point, is a homeomorphism.
     \item[(ii')] $\overline{\inenv(Y,\lambda)}$ varies continuously in $(Y,\lambda)\in\T(S)\times\crl(S)$, i.e. $\overline{\inenv(Y_n,\lambda_n)}$ converges to $\overline{\inenv(Y,\lambda)}$ locally uniformly as $Y_n\to Y$ in $\T(S)$ and $\lambda_n\to\lambda$ in $\crl(S)$ with the Hausdorff topology. 
 \end{enumerate} 
\end{theorem}

From these considerations, we may give a topological description of the ancestor space from which the \enquote{cone over cone} construction of the envelope descends as follows, where the topology of $\cup_{\lambda'}\cone(\lambda')$ is described in Remark \ref{rmk:union:cone}.

\begin{corollary}
\label{cor:bnd:description}
    The space $\bnd(X,Y)$ is homeomorphic to $\cup_{\lambda'}\cone(\lambda')$, where the union is taken over all chain recurrent geodesic laminations $\lambda’$ strictly containing $\Lambda(X,Y)$.
    \end{corollary}

 \begin{remark}
 Theorem \ref{thm:env:top},  Theorem \ref{thm:env:shape}, and Theorem \ref{thm:structure:out:in} answer the following questions: 
 \begin{itemize}
     \item (F. Gu\'eritaud, \cite[Problem 3.2]{Su2016}) Given $X,Y\in\T(S)$, describe the envelope $\env(X,Y)$.
     \item (K. Rafi, \cite[Problem 3.4]{Su2016}) Does $\env(X,Y)$ depend continuously on $X$ and $Y$?
 \end{itemize} 
 \end{remark}

 \color{black}
\subsubsection{Envelopes with an endpoint on the Thurston boundary.}

 In the paper \cite{PanWolf2022}, we observed that the harmonic stretch rays from an interior point $X \in \T(S)$ displayed \tec space as a cone over the Thurston boundary.  This suggests that we may extend the notion of the envelope $\env(X,Y)$ from geodesics that connect interior points of \tec space to a set $\env(X, \eta)$ of geodesics which connect an interior point $X \in \T(S)$ to a Thurston boundary point $\eta$.  Roughly put, our results on these semi-infinite envelopes are:
 
 \vskip.2cm
 (i) We extend our results (Theorem~\ref{thm:env:top} and Theorem~\ref{thm:env:shape}) on compact envelopes $\env(X,Y)$ to these envelopes $\env(X, \eta)$, proving a description of the shape of the envelope (Theorem \ref{thm:env:shape:tree}) and its continuity (Theorem \ref{thm:continuity2}). 

 \vskip.2cm

(ii) The accumulation set of $\env(X, \eta)$ on the Thurston boundary must include $\eta$, but could be larger: in Proposition~\ref{prop:msl:bnd2}, we give necessary and sufficient conditions for a Thurston boundary point to be in the accumulation set of $\env(X, \eta)$ (exhibiting that accumulation set as a star over $\eta$).

 \vskip.2cm

Along the way, we prove some naturality of the Thurston boundary as a visual boundary from a pole, i.e. that the visual cone structures vary continuously with their endpoints (Corollary \ref{cor:con:compactification}).

We also connect the geometry of the ray $[X,\eta]$ to the geometry of maps from the universal cover of $X$ to the $\R$-tree representing $\eta$. I.e. If $[X,\eta]$ is the harmonic stretch ray from $X\in \T(S)$ to $\eta \in \partial\T(S)$, and $T_\eta$ is the real tree representing $\eta$, then we show (Corollary \ref{cor:optimal:lip:tree}) that the maximally stretched lamination $\Lambda(X,\eta)$ along the ray (has a lift which) is exactly the intersection of the maximally stretched loci of  all (equivariant) optimal Lipschitz maps to the tree $T_\eta$. 

Section~\ref{sec:boundary} will develop the definitions and constructions for precise statements of these results.

\subsubsection{A characterization of harmonic stretch lines using the stretch lamination.}
As an application of the parametrization of out-envelopes, we obtain a new characterization of harmonic stretch lines. Recall that for any $X\in\T(S)$ and any geodesic lamination $\lambda$, a  \emph{piecewise harmonic stretch line} (or a \emph{piecewise harmonic stretch line defined by $(X, \lambda, \gamma)$}),  which passes through $X$ and which maximally stretches exactly $\lambda$, is constructed by taking a surjective harmonic diffeomorphism from some (possibly disconnected) punctured surface to $X\backslash\lambda$ \cite[Theorem 1.5]{PanWolf2022} and then extending the image across $\lambda$. The aforementioned harmonic diffeomorphism is determined by
 a measured lamination $\gamma$ complementary to $\lambda$ \cite[Theorem 1.13]{PanWolf2022} (see also Section \ref{subsec:theta:hs}).  A  {piecewise harmonic stretch line} is called a \emph{harmonic stretch line} if it is a limit of harmonic map rays in $\T(S)$.  
Harmonic stretch lines provide unique Thurston geodesics between points in \tec space. In \cite{PanWolf2022}, we gave a characterization of harmonic stretch lines in terms of energy difference of harmonic maps. Here we give a new characterization using straightened measured laminations.  
\begin{theorem}[Characterization of harmonic stretch lines]
\label{thm:phs:hs}
	Let $X\in\T(S)$. Let $\lambda$ be a chain recurrent geodesic lamination and $\gamma$ a measured geodesic lamination on $X\backslash\lambda$. Then the piecewise harmonic stretch line defined by $(X,\lambda,\gamma)$ is a harmonic stretch line if and only if $\gamma\in\cone(\lambda)$.
\end{theorem}

This characterization then removes the apparent reliance of the definition of harmonic stretch lines upon a limiting process for maps. 

We provide, in an appendix, a characterization of straightened laminations.

 \subsection{Related results}  The Thurston metric was introduced and investigated by Thurston in 1986 \cite{Thurston1998}.  Thurston proved that metric is Finsler, asymmetric, geodesic but not uniquely geodesic. He constructed a special type of geodesic lines, the Thurston stretch lines, and constructed geodesics for any ordered pair of distinct points via a concatenation of stretch lines.  Recently, the authors of the current paper constructed another type of geodesic lines, the {\it harmonic stretch lines}, and proved the existence and uniqueness of such lines for any ordered pair of points in the Teichm\"uller space. For other constructions of Thurston geodesics and optimal Lipschitz maps, we refer to \cite{PapadopoulosTheret2012, PapadopoulosYamada2017, HuangPapadopoulos2019, GueritaudKassel2017, AlessandriniDisarlo2019, CalderonFarre2021,DaskalopoulosUhlenbeck2020, DaskalopoulosUhlenbeck2022}. For the horofunction boundary of the Thurston metric, we refer to \cite{Walsh2014}. For the rigidity of Thurston metric,   we refer to \cite{Walsh2014, DLRT2020,Pan2020,HOP2021}. For the coarse geometry of the Thurston metric, we refer to \cite{ChoiRafi2007, LRT2012,LRT2015, LMRT2024}.  For the behaviour of geodesics and length functions along geodesics, we refer to \cite{Papadopoulos1991,Theret2007,Theret2008,Theret2010,Theret2014,DLRT2020,Tel2019}. For expositions on Thurston metric, we refer to \cite{PapadopoulosTheret2007,PapadopoulosSu2015,Su2016,PanSu2023}. 
  For an intensive study of the case of the punctured torus, we refer to \cite{DLRT2020}. For envelopes,  Dumas, Lenzhen, Rafi, and Tao proved that the envelope in $\T(S_{1,1})$ is either a geodesic segment or a quadrilateral region bounded by four Thurston stretch segments, and varies continuously on endpoints \cite{DLRT2020}. Later,   Bar-Natan \cite{Bar-NatanThesis} proved that envelopes in $\T(S_{1,1})$ and $\T(S_{0,4})$ have uniformly bounded width:   Thurston geodesics with the common initial point and terminal point stay a uniformly-bounded distance from each other. For generalizations to outer spaces and higher \tec spaces, we refer to \cite{FrancavigliaMartino2011,CDPW2022}.

  \subsection{Organization of the paper}    This paper is organized as follows. In Section \ref{sec:preliminary}, we collect some properties about the Thurston metric. In Section \ref{sec:out-in:env}, we parametrize the out-envelopes and in-envelopes using \emph{straightened laminations}, and prove the continuity of the closure of out-envelopes and in-envelopes. In Section \ref{sec:extendability}, we introduce the notion of extendability and describe the shape of envelopes as \emph{cones over cones}. In Section \ref{sec:continuity}, we prove the continuity of envelopes.  In Section \ref{sec:boundary}, we extend some of our previous results to the Thurston boundary. The paper concludes with the appendix on straightened laminations.

 \subsection*{Acknowledgements}
 The first author is supported by the National Natural Science Foundation of China NSFC 12371073 as well as the Guangzhou Basic and Applied Basic Research Foundation grant number 2024A04J3636. The second author gratefully acknowledges support from the U.S. NSF grant DMS-2005551 and DMS-2429005. We are grateful to Weixu Su and Alexander Nolte for conversations about Theorem \ref{thm:GeodesicLimit} and Remark~\ref{rem:GeodesicLimit}. We would like to thank Xian Dai for her comments on an earlier version of the manuscript, and especially on Remark \ref{rem:concatenation}.

 The authors are very grateful to the referee for the very careful reading and the many thoughtful and useful suggestions.

 \section{Preliminaries}\label{sec:preliminary}

\subsection{\tec space and the Thurston metric}
Let $S$ be an orientable closed surface of genus at least two. The \tec space $T(S)$ is the space of equivalence classes of complex structures on $S$, where two complex structures $X$ and $Y$ are said to be \emph{equivalent} if there exists a conformal map $X\to Y$ which is homotopic to the identity map.  By the Uniformization Theorem,   the \tec space $T(S)$ is also the space of equivalence classes of hyperbolic structures on $S$, where two hyperbolic structures $X$ and $Y$ are said to be \emph{equivalent} if there exists an isometry $X\to Y$ which is homotopic to the identity map.  For simplicity, we  denote the equivalence class of the complex/hyperbolic structure   $X$  by $X$ itself.

For any two hyperbolic surfaces $X$ and $Y$ in $\T(S)$, the Thurston metric $d_{Th}$ is defined as:
\begin{equation*}
	d_{Th}(X,Y):=\log \inf_f \{L(f)  \}
\end{equation*}
where $L(f)$ is the Lipschitz constant of $f$ and where $f$ ranges over all Lipschitz homeomorphisms from $X$ to $Y$ in the homotopy class of the identity map. Thurston \cite[Theorem 8.5]{Thurston1998} gave another characterization of this metric in terms of simple closed curves: 
\begin{equation*}
	d_{Th}(X,Y)=\log \sup_{\alpha} \frac{\ell_\alpha(Y)}{\ell_\alpha(X)}
\end{equation*}
where $\ell_\alpha(\cdot)$ represents the length of the (unique) geodesic representative of $\alpha$ and where $\alpha$ ranges over all homotopically nontrivial simple closed curves. 

Using stretch maps between ideal triangles, Thurston constructed special types of geodesics under this metric,  the \emph{Thurston stretch lines}. Although Thurston stretch lines are rare, in the sense that for any $X\in\T(S)$, the union of unit tangent vectors of Thurston stretch lines at $X$ have zero Hausdorff dimension in the unit tangent sphere of $\T(S)$ over $X$ \cite[Theorem 10.5]{Thurston1998},  any two ordered points in $\T(S)$ can be connected by a geodesic which is a concatenation of Thurston stretch segments \cite[Theorem 8.5]{Thurston1998}.  In particular, the optimal Lipschitz constant $\exp(d_{Th}(X,Y))$ is realized by some homeomorphism from $X$ to $Y$.

 \subsection{Neighbourhoods of laminations.}\label{subsec:nghd}
   Let $X$ be a hyperbolic surface and $\lambda$ a geodesic lamination on $X$. There are two related notions of neighbourhoods of $\lambda$. To avoid confusion, we make the following convention. For any $\epsilon>0$, the \emph{$\epsilon$ neighbourhood} $N_\epsilon(\lambda)$ on $X$ is defined to be:
   \begin{equation*}
       N_\epsilon(\lambda):=\{p\in X: d_X(p,\lambda)<\epsilon\}.
   \end{equation*}
   The \emph{Hausdorff $\epsilon$ neighbourhood} $N_\epsilon^H(\lambda)$  of $\lambda$ in the space of geodesic laminations is defined to be the set of geodesic laminations $\lambda'$ such that 
   \begin{equation*}
       d_H(\lambda,\lambda')<\epsilon,\text{ and }~~d_H(\lambda',\lambda)<\epsilon
   \end{equation*}
   where $d_H(\cdot,\cdot)$ represents the Hausdorff metric on the space of closed subsets of $X$.

\subsection{Balls and the Hausdorff topology}\label{subsec:hausdorff}
Since the Thurston metric is asymmetric, there are \emph{out-going balls} 
\begin{equation*}
	B_{out}(X,r):=\{Y\in\T(S): \dth(X,Y)<r\},
\end{equation*}
 \emph{in-coming balls} 
 \begin{equation*}
 	B_{in}(X,r):=\{Y\in\T(S):\dth(Y,X)<r\},
 \end{equation*}
 and \emph{symmetric balls}
 \begin{equation*}
 	B_{sym}(X,r)= B_{out}(X,r)\cap B_{in}(X,r).
 \end{equation*}
They all define the same topology as the usual one on $\T(S)$ \cite{Liu2001,PapadopoulosTheret2007a}. 
Moreover, for $r< \infty$, the closed balls $\overline{B_{out}(X,r)}$, $\overline{B_{in}(X,r)}$, and $\overline{B_{sym}(X,r)}$ are compact subsects of $\T(S)$ \cite[Proposition 5]{PapadopoulosTheret2007a}.

Recall that the Thurston metric is Finsler \cite[Theorem 5.1]{Thurston1998} (see also \cite[Theorem 2.3]{PapadopoulosSu2015}), geodesic \cite[Theorem 8.5]{Thurston1998}, and complete \cite[Theorem 2]{PapadopoulosTheret2007a}. Thus, for any compact subset $K\subset\T(S)$, there exists a constant $C$ depending on $K$ such that for any $X,Y\in K$, we have
\begin{equation*}
    \frac{1}{C}\cdot \dth(Y,X)\leq \dth (X,Y)\leq C\cdot  \dth(Y,X).
\end{equation*}
   
  	  Consider the symmetrization$$\dth ^S(X,Y):=\max\{\dth(X,Y),\dth(Y,X)\}.$$Throughout this paper, whenever we mention the \emph{Hausdorff topology} on the set of closed subsets of $\T(S)$, we mean the one induced by the symmetrization $\dth^S$. Given a subset $M\subset \T(S)$ and a constant $\epsilon>0$, let $\mathcal{N}_\epsilon(M)$ be the $\epsilon$ neighbourhood of $M$ under the symmetrization $\dth^S$.

\subsection{Maximally stretched laminations} 
Given distinct $X$ and $Y$ in $\T(S)$, the \emph{maximally stretched lamination} $\Lambda(X,Y)$ from $X$ to $Y$\footnote{This is referred to in Thurston's seminal work \cite{Thurston1986} as the \enquote{maximal ratio-maximizing chain recurrent lamination}.} is the largest chain recurrent geodesic lamination $\lambda$  with  following property: there exists an $\exp(d_{Th}(X,Y))$-Lipschitz map, homotopic to the identity map,  from a neighbourhood of $\lambda$ on $X$ to a neighbourhood of $\lambda$ on $Y$, which takes leaves of $\lambda$ on $X$ to corresponding leaves of $\lambda$ on $Y$ by multiplying arclength by a factor of $\exp(d_{Th}(X,Y))$. It is also the union of all chain recurrent geodesic laminations with the aforementioned property \cite[Theorem 8.2]{Thurston1998}.  Equivalently, as a second version, the maximally stretched lamination is also the union of chain-recurrent laminations to which the restriction of every  Lipschitz map from $X$ to $Y$ with the extremal Lipschitz constant, which is homotopic to the identity,  takes leaves of $\lambda$ on $X$ to corresponding leaves of $\lambda$ on $Y$ by multiplying arclength by a factor of $\exp(d_{Th}(X,Y))$.  (For the equivalence between these two descriptions, we refer to \cite[Section 9]{GueritaudKassel2017}, see also \cite[Section 5]{DaskalopoulosUhlenbeck2022} for a related discussion.)

We recall the following theorem of Thurston on limits of maximally stretched laminations.

\begin{theorem}[\cite{Thurston1998}, Theorem 8.4] \label{thm:maximally:stretched:lamination} Let $X$ and $Y$ be two distinct hyperbolic surfaces in $\T(S)$.
 Let $X_n,Y_n\in\T(S)$ with $X_n\to X$ and $Y_n\to Y$ as $n\to\infty$. Then  $\Lambda(X,Y)$ contains any lamination in the limit set of $\Lambda(X_i,Y_i)$ with respect to the Hausdorff topology.  \end{theorem}

\begin{remark}
    In Section \ref{subsec:msl:trees}, we extend the discussion of maximally stretched laminations to equivariant Lipschitz maps to $\R$-trees. 
\end{remark}

We end this subsection with the following observation. 
\begin{lemma}\label{lem:multi:curve:lip} Let $X,Y\in\T(S)$ be two distinct hyperbolic surfaces.
    Let $\alpha_n$ be a sequence of multicurves that converges to a sublamination of $\Lambda(X,Y)$ with respect to the Hausdorff topology. Then the ratio
    \begin{equation*}
        \frac{\ell_Y(\alpha_n)}{\ell_X(\alpha_n)} 
    \end{equation*}
    converges to $\exp(\dth(X,Y))$ as $n\to\infty$.
\end{lemma}
\begin{proof}
    Let $d=\exp(\dth(X,Y))$. For simplicity, we assume that $\alpha_n$ is a simple closed curve. For any $\delta>0$, there exists $N=N(\delta)>0$ such that for $n>N$, the geodesic $\alpha_n$ is contained in the $\delta/d$ neighbourhood of $\Lambda(X,Y)$.  We approximate $\alpha_n$ by a leaf segment $\alpha_n^1$ of $\Lambda(X,Y)$ followed by a geodesic arc $\alpha_n^2$ of length less than $2\delta/d$. In particular,
\begin{equation}\label{eq:length:12:alpha}
    \ell_{X}(\alpha_n)\leq \ell_{X}(\alpha_n^1)+\ell_{X}(\alpha_n^2).
\end{equation}
 On the other hand, by the definition of maximally stretched lamination, we see that there exists an $d-$Lipschitz homeomorphism $f$ from a neighbourhood of  $\Lambda(X,Y)$ on $X$ to a neighbourhood of $\Lambda(X,Y)$ on $Y$ that sends the geodesic lamination $\Lambda(X,Y)$ on $X$ to the geodesic lamination $\Lambda(X,Y)$ on $Y$ and expands the arclength along $\Lambda(X,Y)$ by the factor $d$. In particular, the lengths of $f(\alpha_n^1)$ and $f(\alpha_n^2)$ on $Y$ satisfy:
\begin{eqnarray}
    \ell_Y(f(\alpha_n^1))&=& d \cdot \ell_{X}(\alpha_n^1),
    \label{eq:length:1:alpha}\\ 
     \ell_Y(f(\alpha_n^2))&\leq&  d\cdot \ell_{X}(\alpha_n^2)<2\delta.\label{eq:length:2:alpha}
\end{eqnarray}
Choosing $\delta$ to be sufficiently small and applying the Anosov Closing Lemma (cf. \cite[Lemma 13.1]{Anosov1969}), we see that the length $\ell_Y(\alpha_n)$ of the geodesic homotopic to $\alpha_n$ on $Y$ satisfies
\begin{equation}\label{eq:length:12Y:alpha}
    |\ell_{Y}(\alpha_n)-\ell_{Y}(f(\alpha_n^1))-\ell_{Y}(f(\alpha_n^2))|<\epsilon
\end{equation}
for some positive constant $\epsilon$ which tends to zero as $\delta\to0$.
 Combining the discussion above, we infer that
\begin{eqnarray*}
    d&\geq& \frac{\ell_Y(\alpha_n)}{\ell_{X}(\alpha_n)} \qquad (\text{by }\eqref{eq:length:12Y:alpha} )
    \\&\geq & 
    \frac{\ell_{Y}(f(\alpha_n^1))+\ell_{Y}(f(\alpha_n^2))-\epsilon}{\ell_{X}(\alpha_n)} \\
  & \geq&  \frac{\ell_{Y}(f(\alpha_n^1))-\epsilon}{\ell_{X}(\alpha_n)}
  \\&
    =& \frac{d\cdot \ell_{X}(\alpha_n^1)-\epsilon}{\ell_{X}(\alpha_n)} \qquad(\text{by }\eqref{eq:length:1:alpha})\\
   &\geq &\frac{d\cdot (\ell_{X}(\alpha_n)-\ell_{X}(\alpha_n^2))-\epsilon}{\ell_{X}(\alpha_n)}\qquad(\text{by }\eqref{eq:length:12:alpha})\\
   &\geq &\frac{d\cdot \ell_{X}(\alpha_n)-2\delta-\epsilon}{\ell_{X}(\alpha_n)}\qquad (\text{by }\eqref{eq:length:2:alpha})\\
   &= & d-\frac{2\delta+\epsilon}{\ell_{X}(\alpha_n)} \\
   &\geq & d-\frac{2\delta+\epsilon}{\mathrm{syst(X)}} 
\end{eqnarray*}
  where $\mathrm{syst(X)}$ is the length of the shortest closed geodesics on $X$. The lemma then follows from the arbitrariness of $\delta$ and the fact that $\epsilon\to0$ as $\delta\to0$. 
\end{proof}
 
  \subsection{Envelopes}
 
 For any distinct points $X,Y$ in $\T(S)$,  the \emph{envelope} $\env(X,Y)$ is defined to be the union of Thurston geodesics from $X$ to $Y$, that is,  \begin{equation*}
 	\env(X,Y)=\{Z: Z\text{ is contained in some geodesic from $X$ to $Y$}\}.
 \end{equation*}
  From the definition, we see that $\env(X,Y)$, as  a closed subset of the (compact \cite{PapadopoulosTheret2007a}) right $\exp(d_{Th}(X,Y))$-ball centered at $X$,  is  compact.

    From the definition of envelope, we see that $Z\in\env(X,Y)$ if and only if $\dth(X,Y)=\dth(X,Z)+\dth(Z,Y)$. Hence, we have
  \begin{equation}\label{eq:char:envelope}
      \env(X,Y)=\{Z: [X,Z]\cup[Z,Y] \text{ is a geodesic}\}.
  \end{equation}
   \color{black}

 Instead of fixing the endpoint of geodesics from $X$, one can also prescribe the maximally stretched lamination. For any $X\in\T(S)$ and any chain recurrent geodesic lamination $\lambda$, the \emph{out-envelopes} $\outenv(X,\lambda)$  and \emph{in-envelopes} $\inenv(X,\lambda)$ are defined as:
 \begin{eqnarray*}
 	\outenv(X,\lambda)&=&\{Z: \Lambda(X,Z)=\lambda\},\\
 	\inenv(X,\lambda)&=&\{Z:\Lambda(Z,X)=\lambda\}.
 \end{eqnarray*}

 We note that this definition allows that the frontier of $\outenv(X,\lambda)$  (respectively, $\inenv(X,\lambda)$) has points $Z$ where $\Lambda(X,Z)$ (respectively, $\Lambda(Z,X)$) strictly contains $\lambda$.

 Note that, by convention, we have that $X\notin \outenv(X,\lambda)$  and $Y\notin \inenv(Y,\lambda)$.  (This will have some consequences in describing intersections of $\outenv(X,\lambda)$ and $\inenv(Y,\lambda)$ with various other sets.)

\begin{remark}\label{rmk:envelope:inclusion}
The following observations follow directly from the definitions of envelope, out-envelope, and in-envelope, and we will use them frequently. Let $X,Y\in\T(S)$ be two distinct hyperbolic surfaces and $\lambda$ be a chain recurrent geodesic lamination. 
 \begin{enumerate}
     \item For any $Z\in\env(X,Y)\setminus\{X,Y\}$, we have $$\env(X,Z)~{\subsetneq}~ \env(X,Y) \text{ and } \env(Z,Y)~{\subsetneq}~ \env(X,Y).$$
     \item  For any $Z\in \outenv(X,\lambda)$, we have $$ \outenv(Z,\lambda)~ {\subsetneq} ~\outenv(X,\lambda).$$
     \item  For any $Z\in\inenv(Y,\lambda)$, we have $$ \inenv(Z,\lambda)~{\subsetneq} ~\inenv(Y,\lambda). $$
 \end{enumerate}
\end{remark}

 In preparation for the next proposition, we say that a chain recurrent geodesic lamination is \emph{maximally CR} if it is not contained in any other chain recurrent geodesic lamination. (Of course, a maximally CR lamination is not necessarily maximal (in the sense of having only ideal triangles as its complementary domains).) 
 
   \begin{proposition}[\cite{DLRT2020}, Propositions 5.2 and 5.3] \label{prop:DLRT2020} For any $X,Y\in\T(S)$ and any chain-recurrent lamination $\lambda$,  the out-envelopes, in-envelopes, and envelopes have the following properties. 
\begin{enumerate}
	\item If $\lambda$ is maximally CR, then for any completion $\hat{\lambda}$ (as a geodesic lamination) of $\lambda$, the set $\outenv(X,\lambda)$ is the Thurston stretch ray starting at $X$ associated with $\hat{\lambda}$, and the set $\inenv(Y,\lambda)$ is the Thurston stretch ray associated with $\hat\lambda$ ending at $Y$.  
	\item The closure of $\outenv(X,\lambda)$ consists of points $Y$ with $\lambda\subset\Lambda(X,Y)$. Similarly, the closure of $\inenv(Y,\lambda)$   is the set of points $Z$ with $\lambda\subset \Lambda(Z,Y)$. 
	\item  If $\lambda$ is a simple closed curve, then $\outenv(X,\lambda)$ and $\inenv(Y,\lambda)$  are open sets.
	\item $\env(X,Y)=\overline{\outenv(X,\Lambda(X,Y))}\cap \overline{\inenv(Y,\Lambda(X,Y))}$.
\end{enumerate}
\end{proposition} 

Here the proposition identifies, in the first and third cases, the two extremal situations for $\outenv(X,\lambda)$. Note that, in the first case, the set $\outenv(X,\lambda)$ is independent of the completion chosen: the important property is that $\lambda$ is maximally CR. (cf. Theorem~\ref{thm:env:shape}.)

  \begin{remark}\label{rem:EnvelopeByStretchLamination}
      For any $Z\in\env(X,Y)$, the union $[X,Z]\cup[Z,Y]$ is a geodesic from $X$ to $Y$, which gives an optimal Lipschtiz map from $X$ to $Y$  that maximally stretches exactly along $\Lambda(X,Z)\cap\Lambda(Z,Y)$, so by (the second version of) the definition of maximally stretched lamination, we have $\Lambda(X,Y)\subset (\Lambda(X,Z)\cap \Lambda(Z,Y))$. To see the reverse inclusion, note that (the first version of) the definition of maximally stretched lamination implies that $\Lambda(X,Y)$ contains $\Lambda(X,Z)\cap\Lambda(Z,Y)$. Therefore, we have the identity $\Lambda(X,Y)=\Lambda(X,Z)\cap \Lambda(Z,Y)$.  On the other hand, if $Z\in\T(S)$ satisfies the identity $\Lambda(X,Y)=\Lambda(X,Z)\cap \Lambda(Z,Y)$, then $[X,Z]\cup[Z,Y]$ is a geodesic from $X$ to $Y$, hence $Z\in\env(X,Y)$. In summary:\end{remark}
      \begin{corollary}\label{cor:env in terms of stretch lams}
    For any distinct $X$ and $Y$ in $\T(S)$, the envelope $\env(X,Y)$ can be characterized as: 
      \begin{equation*}
          \env(X,Y)=\{Z\in\T(S):\Lambda(X,Y)=\Lambda(X,Z)\cap \Lambda(Z,Y)\}.
      \end{equation*}
      \end{corollary}
 Note that on the right hand side of the displayed equation, we adopt the convention that $\Lambda(X,X)=X$.
      
  \subsection{(Piecewise) harmonic stretch lines}\label{subsec:PHSL}

  A principal tool in this paper will be \enquote{piecewise harmonic stretch lines}.  These will be built from \enquote{piecewise harmonic stretch maps}, the basic result (\cite[Theorem 1.5]{PanWolf2022}) for which we present next.

\begin{theorem}[piecewise harmonic stretch map]\label{thm:generalized:stretchmap}
    Let  $Y\in\T(S)$ be any closed hyperbolic surface, and let  $\lambda$ be any geodesic lamination. Then for any harmonic diffeomorphism $f:X\to Y\setminus \lambda$ from some (possibly disconnected) punctured surface $X$,  there is a new hyperbolic surface
    \begin{equation*}
     Y_t:= \PHSL_{Y,\lambda,f}(t)\in\T(S)
    \end{equation*}
    depending analytically on $\{t>0\}$ such that
    \begin{enumerate}[(i)]
    \item the induced map $f_t:X\to Y_t\setminus\lambda$ is a   harmonic diffeomorphism with  Hopf differential  $t\Hopf(f)$;
      \item for any $0< s<t$, the map ($f_t\circ f_s^{-1}$)  extends to a homeomorphism from $Y_s$ to $Y_t$ that is $\sqrt{t/s}$-Lipschitz with (pointwise) Lipschitz constant strictly less than  $\sqrt{t/s}$ in $Y_s\backslash\lambda$, but
       exactly expands arc length on $\lambda$ by the constant factor $\sqrt{t/s}$.
    \end{enumerate}
  \end{theorem}
  The family of hyperbolic structures $\PHSL_{Y,\lambda,f}(t)$  constructed above is called a {\it piecewise harmonic stretch line}. It admits a canonical orientation coming from the orientation of the  positive real ray $\{t>0\}$. In that orientation, a piecewise harmonic stretch line is a (reparametrized) geodesic in the Thurston metric.  Whenever we say a piecewise harmonic stretch line, we mean a directed line.

A piecewise harmonic stretch line is defined to be a \emph{harmonic stretch line} if it appears as the limit of a sequence of \emph{harmonic map rays} in $\T(S)$ (see  \cite[Definition 12.1]{PanWolf2022}). 

In terms of harmonic stretch lines, we have the basic existence and uniqueness result.

 \begin{theorem} [\cite{PanWolf2022}]\label{thm:PanWolf2022:uniqueness} For any two distinct points $X,Y\in\T(S)$, we have the following.
 \begin{enumerate}
     \item 
 	There exists a unique harmonic stretch line proceeding from $X$ through $Y$.
  \item The harmonic stretch line proceeding from $X$ through $Y$ maximally stretches exactly  $\Lambda(X,Y)$. 
 \end{enumerate}
 \end{theorem}
 \begin{proof}
     The first statement is exactly \cite[Theorem 1.6]{PanWolf2022}. The second statement follows from the construction of harmonic stretch lines. Notice that a harmonic stretch line is also a piecewise harmonic stretch line. 
     By \cite[Theorem 1.5]{PanWolf2022}, there exist a geodesic lamination $\lambda$ and a piecewise harmonic stretch map from $X$ through $Y$ which is an $\exp(\dth(X,Y))$-Lipschitz homeomorphism that maximally stretches exactly $\lambda$. Moreover, by 
     \cite[Definition 12.1]{PanWolf2022}, that harmonic stretch line proceeding from $X$ through $Y$ is a limit of harmonic map rays.
     It then follows from \cite[Lemma 4.6]{PanWolf2022} that the geodesic lamination $\lambda$ previously mentioned is chain-recurrent. 
     The statement now follows from the definition of $\Lambda(X,Y)$, since the piecewise harmonic stretch map from $X$ to $Y$ is an $\exp(\dth(X,Y))$-Lipschitz homeomorphism that maximally stretches the chain-recurrent geodesic lamination $\lambda$: we conclude that $\Lambda(X,Y) = \lambda$.
 \end{proof}
 
  For any distinct $X$ and $Y$ in $\T(X)$, let $\HSL(X,Y)$ be the harmonic stretch line proceeding from $X$ through $Y$. The closed subinterval of this line from $X$ to $Y$ is called the \emph{harmonic stretch segment} from $X$ to $Y$, denoted by $[X,Y]$.

 Note that the uniqueness of harmonic stretch lines implies that for any $V,W$ in $[X,Y]$ with the ordering $(X,V,W,Y)$, the harmonic stretch segment $[V,W]$ is exactly the subset of points of $[X,Y]$ which lie in between $V$ and $W$.

 \begin{proposition}[\cite{PanWolf2022}, Proposition 12.14]\label{prop:PanWolf2022}  Let $X$ and $Y$ be two distinct points in $\T(S)$. Let $X_n,Y_n\in\T(S)$ with $X_n\to X$ and $Y_n\to Y$ as $n\to\infty$, then the harmonic stretch lines  $\HSL(X_n,Y_n)$  converge to $\HSL(X,Y)$ locally uniformly as $n\to\infty$. In particular, $[X_n,Y_n]$ converges to $[X,Y]$ uniformly. 
 \end{proposition}

By {\it locally uniform convergence} mentioned in the statement of Proposition \ref{prop:PanWolf2022}  we mean the following. We parametrize $\HSL(X_n,Y_n)$ and $\HSL(X,Y)$ by arclength, and setting $\HSL(X_n,Y_n)(0)=X_n$ and $\HSL(X,Y)(0)=X$. Then $\HSL(X_n,Y_n)$ converges to $\HSL(X,Y)$ locally uniformly as $n\to\infty$ if for any $r>0$, the restriction $\HSL(X_,Y_n)|_{[-r,r]}$ converges to $\HSL(X,Y)|_{[-r,r]}$ uniformly. The parameterizations for $[X_n,Y_n]$ and $[X,Y]$ are slightly differently due to the fact that $\dth(X_n,Y_n)$ and $\dth(X,Y)$ may not be equal. We renormalize speed of parametrization so that they have the same domain $[0,\dth(X,Y)]$.     

Alternatively, we may describe the convergence using the Hausdorff topology on $\T(S)$ induced from the symmetrization of the Thurston metric (see Section \ref{subsec:hausdorff}). Then $[X_n,Y_n]$ converges uniformly to $[X,Y]$ means that, as closed subsets of $\T(S)$, the segment $[X_n,Y_n]$ converges to $[X,Y]$ with respect to the Hausdorff topology. 

To apply this viewpoint to the setting of 
$\HSL(X_n,Y_n)$ and $\HSL(X,Y)$ we need to add some more explanation.
More precisely, for any $Z\in\T(S)$ and for any $\rho>0$, let
$$B(Z,\rho):=\{Z'\in\T(S): \dth(Z,Z')<\rho,~\dth(Z',Z)<\rho\}.$$
Then $\HSL(X_n,Y_n)$ converges to $\HSL(X,Y)$ locally uniformly means that, for any $\rho>\max\{\dth(X,Z),\dth(Z,X)\}$, the component of $\HSL(X_n,Y_n)\cap \overline{B(Z,\rho)}$ that contains $X_n$  converges with respect to the Hausdorff topology to the component of $\HSL(X,Y)\cap \overline{B(Z,\rho)}$ that contains $X$, as closed subsets of $\T(S)$
\footnote{Here we consider components containing $X_n$ and $X$ respectively, instead of $\HSL(X_n,Y_n)\cap \overline{B(Z,\rho)}$ and $\HSL(X,Y)\cap \overline{B(Z,\rho)}$, to avoid the possibility that $\HSL(X,Y)\cap \overline{B(Z,\rho)}$ may contain  singleton components (disjoint from $X$) while $\HSL(X_n,Y_n)\cap \overline{B(Z,\rho)}$ does not.}.

\subsection{Concatenation of stretch lines} 
As mentioned earlier, Thurston \cite{Thurston1998} constructed many geodesics using concatenations of Thurston stretch lines. The key point of this construction is the behaviour of the maximally stretched lamination (Theorem \ref{thm:maximally:stretched:lamination}).  Accordingly, this method applies to piecewise harmonic stretch lines as well.  

In particular, after two preparatory lemmas on the structure of the envelope $\env(X,Y)$ near $X$ and $Y$ as well as on how piecewise harmonic stretch lines whose stretch laminations contain $\Lambda(X,Y)$ meet the envelope $\env(X,Y)$, we will describe (Remark~\ref{rem:concatenation}) how we may traverse a path from $X$ to $Y$ along piecewise harmonic stretch segments that are contained within the envelope $\env(X,Y)$.

We begin with the following lemma describing the local neighborhoods in the envelope of the initial and terminal points.

\begin{lemma}\label{lem:near:XY}
 	Let $X,Y\in\T(S)$. Then there exists neighbourhoods $U$ and $V$  of  $X$ and $Y$, respectively, such that 
 	\begin{equation*}
 		\env(X,Y) \cap U= \overline{\outenv(X,\Lambda(X,Y))} \cap U 
 	\end{equation*}
  and
  \begin{equation*}
     \quad\env(X,Y) \cap V=  \overline{\inenv(Y,\Lambda(X,Y))}\cap V.
  \end{equation*}
 \end{lemma}

  \begin{proof}
To simplify notation, let $\Lambda:=\Lambda(X,Y)$.
   Recall that any hyperbolic surface $Z\in\T(S)$ admits a measured horocycle (partial) foliation in a neighbourhood of $\Lambda$ on $Z$ and transverse to $\Lambda$. 
      Let us denote this foliation by $F_\Lambda(Z)$.  
   By \cite[Theorem 8.5]{Thurston1998},  for any $Z\in \overline{\outenv(X,\Lambda)}$, there exists a Thurston geodesic connecting $X$ to $Z$ which is a concatenation of Thurston stretch lines and whose maximally stretched lamination contains $\Lambda$. Accordingly, these partial foliations $F_\Lambda(Z)$ and $F_\Lambda(X)$ satisfy: $F_\Lambda(Z)=\exp(d_{Th}(X,Z))F_\Lambda(X)$. Similarly, we have $F_\Lambda(Y)=\exp(\dth(X,Y))F_\Lambda(X)$.  Therefore,  
   \begin{equation*}
   	F_\Lambda(Y)=\frac{\exp(\dth(X,Y))}{\exp(d_{Th}(X,Z))}F_\Lambda(Z). 
   \end{equation*}
  Consider the map from a neighbourhood of $\Lambda$ on $Z$ to a neighbourhood of $\Lambda$ on $Y$, which maps leaves of $F_\Lambda(Z)$ to leaves of $F_\Lambda(Y)$ and takes leaves of $\Lambda$ on $Z$ linearly to corresponding leaves of $\Lambda$ on $Y$. This gives a Lipschitz homeomorphism from a neighbourhood of $\Lambda$ on $Z$ to a neighbourhood of $\Lambda$ on $Y$ which maximally stretches exactly along $\Lambda$.  
  
On the other hand,  by Theorem \ref{thm:maximally:stretched:lamination}, for any $\epsilon>0$, there exists a neighbourhood $U=U(X,Y,\epsilon)$ of $X$ in $\T(S)$ such that for any $Z\in U$, the maximally stretched lamination $\Lambda(Z,Y)$ is contained in the $\epsilon$ neighbourhood of $\Lambda$, see Section \ref{subsec:nghd} for the clarification of \enquote{neighbourhood}. (To see this, suppose to the contrary that there exists $\epsilon>0$ such that for any neighbourhood $U$ of $X$ in $\T(S)$, there exists $Z\in U$ such that the maximally stretched lamination $\Lambda(Z,Y)$ is not contained in the $\epsilon$ neighbourhood of $\Lambda$. In particular, this implies that, by taking $U$ arbitrary small, there exists $ Z_n\in\T(S)$ that converges to $Z$ as $n\to\infty$ such that $\Lambda(Z_n,Y)$ is not contained in the $\epsilon$ neighbourhood of $\Lambda$. This constradicts Theorem \ref{thm:maximally:stretched:lamination}.) 
 Hence, for any $Z\in \overline{\outenv(X,\Lambda)}\cap U$, the Lipschitz homeomorphism constructed in the last paragraph implies that $\Lambda(Z,Y)=\Lambda$. In particular, this means that $\Lambda(X,Y)=\Lambda=\Lambda(X,Z)\cap \Lambda(Z,Y)$.    By Corollary \ref{cor:env in terms of stretch lams}, this proves that $Z\in \env(X,Y)$. The arbitrariness of $Z\in \overline{\outenv(X,\Lambda)}\cap U$ implies that $(\overline{\outenv(X,\Lambda)}\cap U) \subset ( \env(X,Y)\cap U)$.  On the other hand, by item (iv) in Proposition \ref{prop:DLRT2020}, we see that 
  \begin{equation*}\label{eq:local:env}
  	(\env(X,Y)\cap U)\subset (\overline{\outenv(X,\Lambda)}\cap U)
  \end{equation*}
   holds for any open subset $U\subset\T(S)$. Therefore, we have 
 \begin{equation*}
 	\env(X,Y)\cap U= \overline{\outenv(X,\Lambda)}\cap U.
 \end{equation*}
 Similarly, one can prove that $\env(X,Y)\cap V= \overline{\inenv(Y,\Lambda)}\cap V$ for some neighbourhood $V$ of $Y$.
   \end{proof}
    
We next begin to relate (piecewise) harmonic stretch lines to envelopes: these relationships will underlie our approach to the structure of envelopes.

 \begin{lemma}\label{lem:concatenation}
	Let $S$ be an orientable closed surface of genus at least two. Let $X,Y\in\T(S)$ be two distinct points with maximally stretched lamination $\Lambda(X,Y)$. Let  $\lambda$ be an arbitrary geodesic lamination which contains $\Lambda(X,Y)$. Then
	\begin{enumerate}
		\item   for any piecewise harmonic stretch line $\mathrm{PHSL}(X,\lambda)$  through $X$ which maximally stretches $\lambda$,  the intersection 
	\begin{equation*}
		\mathrm{PHSL}(X,\lambda)\cap \env(X,Y)
	\end{equation*}
	 is a non-trivial, oriented, closed geodesic segment with initial point $X$;
	 \item 
	 for any piecewise harmonic stretch line $\mathrm{PHSL}(Y,\lambda)$  through $Y$ which maximally stretches $\lambda$,  the intersection 
	\begin{equation*}
		\mathrm{PHSL}(Y,\lambda)\cap \env(X,Y)
	\end{equation*}
	 is a non-trivial, oriented, closed geodesic segment with terminal point $Y$.
	\end{enumerate}
	\end{lemma}
\begin{proof}
(i) The point  $X$ cuts $\PHSL(X,\lambda)$ into two oriented geodesic rays. Let $R^+$ and $R^-$ be respectively the one that starts at $X$ and the one that ends at $X$. For any point $Z'\in R^-$,  the maximally stretched lamination  $\Lambda(Z',X)$ from $Z'$ to $X$ is exactly $\lambda$ which also (by hypothesis) contains $\Lambda(X,Y)$.  
By Lemma \ref{lem:multi:curve:lip}, we see that for any multicurve $\alpha_n$ converging to $\Lambda(X,Y)$ with respect to the Hausdorff topology, the ratio
\begin{equation}\label{eq:ratio:ZX}
    \frac{\ell_X(\alpha_n)}{\ell_{Z'}(\alpha_n)}
\end{equation}
converges to $\exp(\dth(Z',X))>1$. 
This implies that the maximally stretched lamination $\Lambda(X,Z')$ from $X$ to $Z'$ does not contain $\Lambda(X,Y)$, otherwise  the ratio above would converge to  $\exp(-\dth(X,Z'))<1$. In particular, by Corollary~\ref{cor:env in terms of stretch lams}, the point $Z'$ is not contained in $\env(X,Y)$.  Therefore,  $\PHSL(X,\lambda)\cap\env(X,Y)$  is disjoint from $R^-$. 

Next we claim that $\PHSL(X,\lambda)\cap\env(X,Y)$ is a non-trivial segment of $R^+\cup\{X\}$. Since by assumption $\lambda$ contains $\Lambda(X,Y)$, it then follows from Proposition \ref{prop:DLRT2020}(ii) together with Theorem~\ref{thm:generalized:stretchmap}(ii)  that $\PHSL(X,\lambda)$ is contained in the closure $\overline{\outenv(X,\Lambda(X,Y))}$. Combining with Lemma~\ref{lem:near:XY}, we see that the intersection $\PHSL(X,\lambda)\cap\env(X,Y)$ is non-trivial. I.e. the intersection $\PHSL(X,\lambda)\cap\env(X,Y)$ is not the singleton $\{X\}$.
For any $Z\in \PHSL(X,\lambda)\cap\env(X,Y)$ different from $X$,  the envelope $\env(X,Z)$ is contained in $\env(X,Y)$.  In particular, the subinterval of $\PHSL(X,\lambda)$, whose  endpoints are $X$ and $Z$, is also contained in $\env(X,Y)$.  Therefore, the intersection $\PHSL(X,\lambda)\cap\env(X,Y)$ is a non-trivial oriented geodesic segment with $X$ being the initial point. That $\PHSL(X,\lambda)\cap\env(X,Y)$ is a closed set follows directly from the definitions of $\PHSL(X,\lambda)$ and $\env(X,Y)$.

(ii) Similarly, one can prove the intersection  $\PHSL(Y,\lambda)\cap\env(X,Y)$ is a non-trivial, oriented, closed geodesic segment with $Y$ being the terminal point.
\end{proof}

\begin{remark} \label{rem:concatenation}
	We can use Lemma \ref{lem:concatenation} to construct geodesics from $X$ to $Y$ using concatenations of piecewise harmonic stretch lines, generalizing Thurston's construction of concatenations of Thurston stretch lines.  Let $\lambda_1$ be a geodesic lamination which contains $\Lambda(X,Y)$ and $\PHSL(X,\lambda_1)$ a piecewise harmonic stretch line through $X$ which maximally stretches exactly $\lambda_1$. 
    (There are many such piecewise stretch lines. Consider the crowned hyperbolic surface $X\setminus\lambda_1$. By \cite[Proposition A.5]{PanWolf2022}, there exists (at least) a harmonic diffeomorphism $f_\infty:X_\infty\to X\setminus\lambda_1$ from some (possibly disconnected) punctured Riemann surface $X_\infty$. We then apply Theorem \ref{thm:generalized:stretchmap} to get a piecewise harmonic stretch line that passes through $X$ and maximally stretches exactly along $\lambda_1$.)
    By Lemma \ref{lem:concatenation}, the intersection 
	$\PHSL(X,\lambda_1)\cap \env(X,Y)$ is a non-trivial oriented segment with initial point $X$. Let $X_1$ be the terminal point of this segment. Then $\Lambda(X_1,Y)$ strictly contains $\Lambda(X,Y)$. Otherwise, since $\PHSL(X,\lambda_1)$ is a piecewise harmonic stretch line which passes through $X_1$ and which maximally stretches exactly $\lambda_1$, then by Lemma \ref{lem:concatenation}, $$\PHSL(X,\lambda_1)\cap \env(X_1,Y)$$ 
	is an oriented segment with initial point $X_1$, meaning that $X_1$ is not the terminal point of $\PHSL(X,\lambda_1)\cap\env(X,Y)$ (since $\env(X_1, Y) \subset \env(X,Y)$).    Repeating the preceding construction for $X_1$ and $Y$ several times, we get a sequence of points $X_i\in\env(X,Y)$ and  a sequence of piecewise harmonic stretch segments from $X_i$ to $X_{i+1}$, such that $\Lambda(X_{i+1},Y)$ strictly contains $\Lambda(X_i,Y)$.  Since the underlying surface $S$ is of finite type,  we see that there exists a constant $K$ depending only on $S$ such that the strictly increasing sequence 
	\begin{equation*}
		\Lambda(X,Y)\subsetneq \Lambda(X_1,Y)\subsetneq\cdots\subsetneq \Lambda(X_n,Y)\subsetneq\cdots
	\end{equation*}
	 must terminate for some finite index $k\leq K$, i.e. $X_k=Y$. 
	This means that there is a geodesic from $X$ to $Y$ which is a concatenation of a finite number of piecewise harmonic segments and which starts with a segment contained in $\PHSL(X,\lambda_1)$.
\end{remark}

\subsection{Train tracks}\label{sec:traintrack}
Comment: We added this subsection about train tracks, which will be used in Section \ref{sec:out-in:env} and the Appendix.
In this subsection, we briefly introduce \emph{train tracks}. For more details, we refer to \cite{PennerHarer1992}.

A \emph{train track} on $S$ is a finite collection $\tau$  of embedded CW-complexes, each comprising vertices (called \emph{switches}) and edges (called \emph{branches}), such that
\begin{enumerate}
\item (non-degeneracy) any switch is contained in the interior of some $C^1$ path embedded in $\tau$.
\item (smoothness) $\tau$ is smooth away from switches; any branch incident to a switch is tangent to the above-mentioned $C^1$ path at the switch.
\item (geometry) Suppose that $\Sigma$ is a component of $S\backslash\tau$. Let $D(\Sigma)$ be the double of $\Sigma$ along the $C^1$ frontier edges of $\Sigma$. Thus, non-smooth points in the frontier of $\Sigma$ give rise to punctures of $D(\Sigma)$. We require that the Euler characteristic $\chi(D(\Sigma))$ of $D(\Sigma)$ be negative.
\end{enumerate}

A (geodesic) lamination $\alpha$ is said to be \emph{carried} by $\tau$ if there exists a $C^1$ map $\phi:S\to S$ such that
\begin{itemize}
    \item $\phi(\alpha)\subset\tau$,
    \item $\phi$ is homotopic to the identity, and
    \item the restriction of the differential $d\phi$ to the tangent line to $\alpha$ at $p$ is nonzero for every $p\in \alpha$. 
\end{itemize}
For each switch $v$ of $\tau$, we fix a direction in the tangent line $T_v(\tau)$ to $\tau$ at $v$. The end $e$ of a branch $b$ of $\tau$ which is incident on $v$ may then be called \emph{incoming} if the direction of the one-sided tangent vector to $e$ at $v$ agrees with this direction, \emph{outgoing} or not.
A \emph{weight system} $w$ on $\tau$ is a function that assigns to each branch $b$ of $\tau$ a nonnegative real number $w(b)\in \R\cup\{\infty\}$, called the \emph{weight} on $b$, which satisfies for each switch $v$ of $\tau$ the \emph{switch condition}:
\begin{equation*}
    w(e_1)+\cdots w(e_r)=w(e_{r+1})+\cdots + w(e_{r+t})
\end{equation*}
where $e_1,\cdots,e_r$ are incoming ends of branches which are incident on $v$ and $e_{r+1},\cdots,e_{r+t}$ are the outgoing ones (and the weight $w(e_i)$ is the weight of the branch that contains $e_i$). 

 Every geodesic lamination is carried by some train track $\tau$ (\cite[Theorem 1.6.5]{PennerHarer1992}), and any transverse measure on that lamination induces a weight system on $\tau$ (\cite[Proposition 1.7.5]{PennerHarer1992}). Conversely, every weight system on a train track $\tau$ gives a measured geodesic lamination carried by $\tau$ (\cite[Construction 1.7.7]{PennerHarer1992}).
\color{black}

 \section{Out-envelopes and in-envelopes}\label{sec:out-in:env}

The goal of this section is to describe $\outenv(X,\lambda)$ and $\inenv(Y,\lambda)$ for any closed hyperbolic surfaces $X$ (and $Y$) and for any chain-recurrent geodesic lamination $\lambda$ on $X$ (or $Y$). As we shall see, this is closely related to the \enquote{tangent space} of the space of measured laminations,  developed by Thurston in \cite[Section 6]{Thurston1998}. 

\subsection{Measured laminations on $X\setminus\lambda$} \label{sec:ML:crown} In this subsection, we shall associate to any measured geodesic lamination on $X$,  whose support is \enquote{nearly parallel} to $\lambda$, a measured lamination in $X\setminus\lambda$, following \cite[Section 6]{Thurston1998}. 
For our purpose, some modifications are required, but the result will be a space $\cone(\lambda)$ comprising  measured laminations on $X\backslash\lambda$ which, near $\lambda$, approximate $\lambda$.

A \emph{geodesic lamination on the crowned surface $X\backslash\lambda$} is a closed subset foliated by simple geodesics (called \emph{leaves}), each of which is either a simple closed geodesic or a bi-infinite simple geodesic.   A \emph{measured geodesic lamination} (or \emph{measured lamination} for short) on $X\backslash\lambda$ is a geodesic
lamination equipped with a transverse invariant measure, which associates to every
arc transverse to the lamination a Radon measure (see \cite[Page 11-13]{Bonahon2001}). Typical examples of measured laminations on $X\backslash\lambda$ are weighted simple closed geodesics or weighted bi-infinite simple geodesics that limit on ideal vertices of $X\backslash\lambda$ or spiral around a closed boundary component of $X\backslash\lambda$ in either direction.  Let $\far_\lambda$ be the space of measured laminations on $X\setminus\lambda$.

Let $\alpha$ and $\beta$ be two intersecting geodesic laminations on $X$. At each intersection point $p$ of $\alpha$ and $\beta$, there are two angles. We denote by $\theta_p(\alpha,\beta)$ the one which is at most $\pi/2$. Define  the {\it  intersection angle} $\theta(\alpha,\beta)$ to be the supremum of $\theta_p(\alpha,\beta)$, where $p$ ranges over all intersection points. If $\alpha$ and $\beta$ are disjoint, then $\theta(\alpha,\beta)=0$. Let $0<\delta<\pi/2$ be a fixed positive number. Let $U_\delta(\lambda)\subset\ML(S)$ be the subset of measured laminations $\gamma$ with $\theta(\lambda,\gamma)\leq \delta$.  It is clear that $U_\delta(\lambda)\subset U_{\delta'}(\lambda)$ if $\delta<\delta'$.

 Let $U_\delta(\lambda)$ and $\far_\lambda$ be defined as above.
Define 
\begin{equation}\label{eq:far}
	\cut_\lambda:U_\delta(\lambda)\longrightarrow \far_\lambda
\end{equation}
 as follows. Let $\gamma\in U_\delta(\lambda)$. Note that the support of $\gamma$ is the (disjoint) union of minimal sublaminations (cf. \cite[Proposition 1.7.2 and Corollary 1.7.3]{PennerHarer1992}). In particular, every leaf $\ell$ of $\gamma$ is recurrent and dense in the minimal sublamination that contains $\ell$. Therefore, the intersection $\ell\cap (X\backslash\lambda)$ is either an empty set,  $\ell$ itself, or the union of proper arcs in $X\backslash\lambda$. Let $X_0$ be a component of $X\setminus\lambda$. The intersection $\gamma\cap X_0$ has two disjoint subsets: one is a compact sublamination in $X_0$ and the other is the union of proper geodesic arcs in $X_0$.  We shall \enquote{straighten} those proper geodesic arcs to ideal geodesic arcs. More precisely, let $c$ be an arbitrary proper geodesic arc of $\gamma\cap X_0$ with endpoints $p_1$ and $p_2$ (on $\lambda$), see Figure \ref{fig:straighten}.  Let $\lambda_i$ be the boundary component of $X_0$ which contains $p_1$ or $p_2$ respectively.  Here we view $\lambda_i$ as an bi-infinite geodesic, although it might be a simple closed geodesic. The endpoint $p_i$ cuts $\lambda_i$ into two half-infinite geodesic rays, say $\lambda_i'$ and $\lambda_i''$.  The assumption that $\gamma\in U_\delta(\lambda)$ implies that one of  $\lambda_i'$ and $\lambda_i''$, say $\lambda_i'$, makes an angle bigger than $\pi/2$ with $c$ at $p$.  We then straighten the arc $\lambda_1'\cup c\cup \lambda_2'$ to its geodesic representative $\gamma(c)$ (see Figure \ref{fig:straighten}), which is an ideal geodesic in $X_0$ or a boundary ideal geodesic component of $X_0$.  Recall that any pair of proper arcs $c,c'$ in $\gamma\cap X_0$ are disjoint. Hence, their straightenings $\gamma(c)$ and $\gamma(c')$ either coincide or are disjoint.
  Consider the subset $\Gamma_\gamma$ comprising those $\gamma(c)$ such that $c$ is a proper arc of $\gamma\cap(X\setminus\lambda)$ and $\gamma(c)$ is not a boundary ideal geodesic component of $X_0$.  Since elements in $\Gamma_\gamma$  are pairwise disjoint ideal geodesics, it follows that $\Gamma_\gamma$  contains finitely many ideal geodesics.  We rewrite $\Gamma_\gamma$ as $\Gamma_\gamma=\{\gamma_1,\cdots,\gamma_k\}$. 

   \begin{figure}
    \begin{tikzpicture}[scale=1.5]
        \draw (-1,1) ..controls (-2,0.7) and (-3,0.7).. (-4,1)
        (-1,-1) ..controls (-2,-0.7) and (-3,-0.8).. (-4,-1);
        \draw[line width=1, blue] (-3.2,0.82)..controls (-2.8,0.5) and (-2.8,-0.5)..(-3.2,-0.855);
        \draw[line width=1, red] 
        (-4,-1).. controls (-3,-0.5) and (-3,0.5)..(-4,1);
        \draw (-2.5,-1.3)node{(a)} (-1.5,-0.7)node{$\lambda_2$} (-1.5,0.7)node{$\lambda_1$} (-3.1,-1.05)node{$p_2$}   (-3.1,1)node{$p_1$} (-2.8,0) node{$c$}  (-3.6,0)node{$\gamma(c)$};

         \draw (1,1) ..controls (2,0.7) and (3,0.7).. (4,1)
        (1,-1) ..controls (2,-0.7) and (3,-0.8).. (4,-1) ;\draw[line width=1, blue](1.7,0.85)--(3,-0.84);
       \draw [line width=1, red]  (1,1)--(4,-1); \draw (2.5,-1.3)node{(b)} (1.5,-0.7)node{$\lambda_2$} (3.5,0.7)node{$\lambda_1$} (3.1,-1.05)node{$p_2$}   (1.9,1)node{$p_1$} (2.6,-0.5) node{$c$}  (1.6,0.3)node{$\gamma(c)$};
        
          \draw (-1,0-3) ..controls (-2,0.1-3) and (-3,-3).. (-4,1-3)
        (-1,-3) ..controls (-2,-0.1-3) and (-3,-3).. (-4,-1-3);
        \draw[line width=1, blue] (-3.2,0.42-3)..controls (-3,0.2-3) and (-3,-0.2-3)..(-3.2,-0.42-3);
       \draw[line width=1, red]  (-4,-1-3).. controls (-3.5,-0.5-3) and (-3.5,0.5-3)..(-4,1-3); 
           \draw (-2.5,-1.3-2.8)node{(c)} (-2,-0.5-2.8)node{$\lambda_2$} (-2,0.1-2.8)node{$\lambda_1$} (-3.1,-.8-2.8)node{$p_2$}   (-3.1,0.4-2.8)node{$p_1$} (-3.2,-0.2-2.8) node{$c$}  (-4,-0.2-2.8)node{$\gamma(c)$};

            \draw (1,0-3) ..controls (2,0.1-3) and (3,-3).. (4,1-3);
       \draw[line width=1, red] (1,-3) ..controls (2,-0.1-3) and (3,-3).. (4,-1-3) ;
        \draw[line width=1,blue](2.4,0.15-3)--(3,-0.31-3);
          \draw (2.5,-1.3-2.8)node{(d)} (4.2,-0.7-2.8)node{$\lambda_2=\gamma(c)$} (3.4,0.6-2.8)node{$\lambda_1$} (3,-0.7-2.8)node{$p_2$}   (2.3,0.1-2.8)node{$p_1$} (2.8,-0.2-2.8) node{$c$} ;
         \end{tikzpicture}
      \caption{The straightening $\gamma(c)$ of the geodesic arc $c$.}
      \label{fig:straighten}
  \end{figure}

  Next,  to each ideal geodesic $\gamma_i\in \Gamma_\gamma$ we associate a weight $w_i$. To this end, choose an arc $\tau_i$ transverse to the measured lamination $\gamma\in U_\delta(\lambda)$   such that 
  \begin{itemize}
  \item $\tau_i\subset X\setminus\lambda$
   \item the endpoints of $\tau_i$ are not contained in $\gamma\cap(X\setminus\lambda)$,
  	\item $\tau_i$ intersects  exactly those proper arcs $c$ such that $\gamma(c)=\gamma_i$.
  \end{itemize} 
  Such a $\tau_i$ could be found by choosing an arc transverse to a train track that carries $\gamma$ sufficiently closely.
  
  The weight $w_i$ is then defined to be the intersection number $i(\tau_i,\gamma)$.

  Finally, we define $\cut_\lambda(\gamma):= \sum_{i}w_i\gamma_i+\gamma_{cpt}$ where $\gamma_{cpt}$ is the compact sublamination of $\gamma\cap(X\setminus\lambda)$.  The following observations follow directly from the construction of the map $\cut_\lambda$.

\begin{lemma}\label{lem:straighten}
    For any small positive constant $\epsilon$, there exists $\delta_0$ depending only on $\epsilon$ such that every $\gamma\in U_{\delta}(\lambda)$ is contained in the $\epsilon$ neighbourhood (on $X$) of $\lambda\cup \cut_{\lambda}(\gamma)$. 
\end{lemma}
\begin{proof}
  To see this, recall that we may write $\cut_\lambda(\gamma)$ as $\cut_\lambda(\gamma)= \sum_{i}w_i\gamma_i+\gamma_{cpt}$ where $\gamma_{cpt}$ is the compact sublamination of $\gamma\cap(X\setminus\lambda)$, $w_i>0$,  and $\gamma_i$ is a proper bi-infinite simple geodesic on $X\backslash\lambda$. Clearly, we see that $\gamma_{cpt}$, as a sublamination of $\gamma$, is contained in the $\epsilon$ neighbourhood  of $\lambda\cup\cut_\lambda(\gamma)$. It remains to consider the leaves of $\gamma$ intersecting $\lambda$. Let $c$ be a proper arc of $\gamma\cap (X\backslash\lambda)$. By the definition of $\cut_\lambda$, the arc $c$ intersects two boundary geodesic arcs of $X\backslash\lambda$, say, $\lambda_1$ and $\lambda_2$ (possibly $\lambda_1 = \lambda_2$), at intersection points $p_1\in\lambda_1$ and $p_2\in\lambda_2$; moreover, the intersection angle at $p_i$ is less than $\delta_1$ or larger than $\pi-\delta_1$ (counted in the counterclockwise from $\lambda_i$ to $c$), see Figure \ref{fig:straighten}. Let $\gamma(c)$ be the straightening of $c$. Note that $\gamma(c)$ may be one of $\gamma_i$ or a boundary leaf of $X\backslash\lambda$ (see Figure \ref{fig:straighten}(d)). By elementary computation, we see that, for each of the four cases listed in Figure \ref{fig:straighten} depending on whether $\lambda_1$ and $\lambda_2$ are asymptotic in one direction, there exists $\delta_0$ such that for every $\gamma\in U_{\delta_0}(\lambda)$, the arc $c\subset \gamma$ is contained in the $\epsilon$ neighbourhood of $\lambda_1\cup\lambda_2\cup\gamma(c)$. Combining the discussion above, we see that $\gamma$ is contained in the $\epsilon$ neighourhood of $\lambda\cup\cut_\lambda(\gamma)$. This proves the lemma.
  \end{proof}

\begin{lemma}\label{lem:image:cut} 
(i) The image $\cut_\lambda(U_\delta(\lambda))$  stablizes as $\delta\to0$. In other words, there exists $\delta_0=\delta_0(X,\lambda
)>0$ such that for any $0<\delta<\delta'<\delta_0$,  one has  $\cut_\lambda(U_\delta(\lambda))=\cut_\lambda(U_{\delta'}(\lambda))$. 

 (ii) Furthermore, for any $\gamma\in U_\delta (\lambda)$ with $\delta<\delta_0$ there exists $\gamma_n\in U_\delta(\lambda)$ such that 
\begin{itemize}
    \item[(ii-a)] $\cut_\lambda(\gamma_n)=\cut_\lambda(\gamma)$;
    \item[(ii-b)] $\supp(\gamma_n)\to \lambda\cup \supp(\gamma)$ with respect to the Hausdorff topology.
\end{itemize}
  \end{lemma}
  \begin{proof}

      Notice that for any $0<\delta<\delta'<\pi/2$,  we have $U_{\delta}(\lambda)\subset U_{\delta'}(\lambda)$. Hence, $\cut_\lambda(U_\delta(\lambda))\subset \cut_\lambda(U_{\delta'}(\lambda))$.

We now turn to the reverse inclusion. Let $\epsilon=\epsilon(X,\lambda)$ be a small constant such that the complementary region of the $\epsilon$ neighbourhood of $\lambda$ on $X$ stabilizes. Let $\delta_0=\delta_0(X,\lambda)$ be the constant in Lemma \ref{lem:straighten}.  Let $\gamma\in U_{\delta'}(\lambda)$ with $\delta'<\delta_0$. Let $\tau$ be a train track that carries $\lambda\cup \cut_\lambda(\gamma)$ and is contained in the $\epsilon$ neighbourhood of $\lambda \cup \cut_\lambda(\gamma)$ (see Section \ref{sec:traintrack} definitions related to train track).  By Lemma \ref{lem:straighten}, the measured geodesic lamination $\gamma$ is also carried by $\tau$. Let $w(\gamma)$ be the weight system on $\tau$ induced from $\gamma$.
Since $\lambda$ is chain recurrent,  there exists a sequence of multicurves $\alpha_n$ that converges to $\lambda$ with respect to the Hausdorff topology. In particular, for $n$ sufficiently large, the multicurve $\alpha_n$ is also carried by $\tau$. Let $w(\alpha_n)$ be the weight system on $\tau$ induced from $\alpha_n$. Finally, let $\gamma_n$ be the measured geodesic lamination on $X$ constructed from the weight system $w(\gamma)+nw(\alpha_n)$ (\cite[Construction 1.1.7]{PennerHarer1992}). Since for large $n$, the straightening $\cut_\lambda(n\alpha_n)$ is trivial (i.e. the straightening operation takes $n\alpha_n$ into $\lambda$ which does not meet $X \setminus \lambda$), we have $\cut_\lambda(\gamma_n)=\cut_\lambda(\gamma)$, proving (ii-a). Furthermore, as $n\to\infty$, the support of $\gamma_n$ converges to $\lambda\cup \supp(\cut_\lambda(\gamma))$ with respect to the Hausdorff topology.  This proves (ii-b) and implies that $\gamma_n\in U_\delta(\lambda)$ for large $n$. Therefore, we have $\cut_\lambda(\gamma)=\cut_\lambda(\gamma_n)\in\cut_\lambda(U_\delta(\lambda))$ and accordingly the inclusion $\cut_\lambda(U_{\delta'}(\lambda)\subset \cut_\lambda(U_\delta(\lambda))$. This completes the proof of (i) and also the proof of the lemma.
  \end{proof}
 In view of Lemma \ref{lem:image:cut}, let $\cone(\lambda)=\cut_\lambda(U_\delta(\lambda))$ for any $\delta<\delta_0(X,\lambda)$. Notice that $\cone(\lambda)$ is independent of the choice of hyperbolic metrics on $S$.   Intuitively, $\cone(\lambda)$ consists of the images under $\cut_\lambda$ of those measured laminations whose support, restricted to a small neighbourhood of $\lambda$ on $X$, is nearly parallel to $\lambda$. 
\begin{definition}\label{def:ad:lamination} 
	A measured lamination $\gamma$ in $\far_\lambda$ is said to be $\lambda$-straightened if it is contained in $\cone(\lambda)$. 
\end{definition}

The set $\cone(\lambda)$ is called the set of \emph{$\lambda$-straightened measured laminations} complementary to $\lambda$.

  \begin{remark}
      \label{rmk: intersection} Note that a measured lamination $\gamma\in\far_\lambda$ is determined by its intersection number with all proper geodesic arcs on $X\backslash\lambda$ and the orientation along which it spirals around each closed boundary geodesic component of $X\backslash\lambda$.
  \end{remark} 
  \begin{remark}\label{rmk:union:cone}
     Later, we also need $\cup_{\lambda'}\cone(\lambda')$ where $\lambda'$ is a chain recurrent geodesic lamination on $X$ that strictly contains $\lambda$. To define the topology of this set, we extend each $\gamma\in\cone(\lambda')$ to a \emph{generalized measured lamination} on $X$ which is a formal sum $\infty\cdot \lambda'+\gamma$ on $X$. We require that for any geodesic arc $k$ on $X$ that is transverse to $\lambda'$, we have $i(\gamma+\infty\cdot \lambda',k)= \infty$. Then we say $\gamma_n\in\cone(\lambda'_n)$ converges to $\gamma\in \cone(\lambda') $ if $\gamma_n+\infty \cdot \lambda'_n$ converges to $\gamma+\infty\cdot\lambda'$, that is, 
     \begin{itemize}
         \item for any geodesic arc $k$ on $X$ that is transverse to $\lambda'$ we have $i(\gamma_n+\infty\cdot \lambda'_n,k)\to \infty$;
         \item for any compact geodesic arc $k$ on $X\backslash\lambda'$ we have $i(\gamma_n+\infty\cdot \lambda_n',k)\to i(\gamma,k)$;
         \item for any simple closed geodesic $\alpha$ in $\lambda'$ and for sufficiently large $n$, the curve $\gamma_n$ wraps along $\alpha$ in the same direction as $\gamma$.
     \end{itemize}
  \end{remark}

  The following observation is due to Thurston at the end of \cite[Page 28]{Thurston1998}, together with a hint of the proof using train track coordinates;    we provide a proof in the appendix.
Recall that the \emph{stump} of a geodesic lamination is the maximal compact sublamination of $\lambda$, with respect to inclusion,  which admits a transverse invariant measure of full support. The stump of $\lambda$ is the union of all compact sublaminations of $\lambda$,  each of which admits a transverse invariant measure of full support.
  \begin{lemma}[\cite{Thurston1998}, the end of  page 28]\label{lem:image:cut2}
  	Let $\lambda$ be a chain recurrent geodesic lamination and $\gamma\in\far_\lambda$. Then  $\gamma$ is $\lambda$-straightened  if and only if for any orientable component of the stump of $\lambda$, the total measure of leaves of $\gamma$ tending to that component in the two directions must be equal. 
  \end{lemma}

 Define $$\Theta^{out}(\lambda)=\cone(\lambda)\times[0,\infty)/_\sim$$
 and 
  $$\Theta^{in}(\lambda)=\cone(\lambda)\times(-\infty,0]/_\sim$$
  where $(\gamma,t)\sim(\gamma',t')$ if $t=t'=0$.
  The remainder of this section is devoted to showing that $\Theta^{out}(\lambda)$ (resp. $\Theta^{in}(\lambda)$) provides a parametrization  $\outenv(X,\lambda)$ (resp.  $\inenv(Y,\lambda)$), and to prove the continuity of $\overline{\outenv(X,\lambda)}$ and $\overline{\inenv(Y,\lambda)}$ based on these parametrizations.

 \subsection{From $\outenv(X,\lambda)$ to $\Theta^{out}(\lambda)$  and from $\inenv(Y,\lambda)$ to $\Theta^{in}(\lambda)$}\label{subsec:hsr:construction}
 Notice from Theorem~\ref{thm:PanWolf2022:uniqueness} that any point $Z\in \outenv(X,\lambda)$ is contained in a unique harmonic stretch line $\HSL(X,Z)$. We next define a map from $\outenv(X,\lambda)$ to $\Theta^{out}(\lambda)$  (and from $\inenv(Y,\lambda)$ to $\Theta^{in}(\lambda)$ by analogy). 

 \begin{proposition}\label{prop:psi}  
  For any $X,Y\in\T(S)$ and any chain recurrent geodesic lamination $\lambda$,
     there are well-defined maps \begin{equation*}
 	\Psi_{out}: \{X\}\cup\outenv(X,\lambda)\longrightarrow \Theta^{out}(\lambda),\quad \Psi_{in}: \{Y\}\cup\inenv(Y,\lambda)\longrightarrow \Theta^{in}(\lambda)
 \end{equation*}
 such that
 \begin{enumerate}
     \item $\Psi_{out}$ sends $X$ to the cone point of $\Theta^{out}(\lambda)$ and sends $Z\in\outenv(X,\lambda)$ to the pair $(\mu,\dth(X,Z))$, where $\mu$ is the defining measured lamination of the harmonic stretch line $\HSL(X,Z)$,
     \item $\Psi_{in}$ sends $Y$ to the cone point of $\Theta^{in}(\lambda)$ and sends $Z\in\inenv(Y,\lambda)$ to the pair $(\mu,-\dth(Z,Y))$, where $\mu$ is the defining measured lamination of the harmonic stretch line $\HSL(Z,Y)$.    \end{enumerate} 
 \end{proposition} 
 \begin{proof}
  By the definition of harmonic stretch lines, there exists a sequence of harmonic maps $f_n:X_n\to X$ which,  in the sense of \cite[Definition 4.5]{PanWolf2022}, converges to some  harmonic diffeomorphism $f_\infty:X_\infty\to X\setminus\lambda$ defining $\HSL(X,Z)$. Removing all half-planes and half-infinite cylinders from the \emph{horizontal foliation} of the Hopf differential of $f_\infty$ gives a measured foliation on $X\setminus\lambda$, which corresponds to a measured lamination $\mu\in\far_\lambda$. Consider the measured lamination $\mu_n$ corresponding to the horizontal foliation of the Hopf differential of $f_n$. It is clear that $\cut_\lambda(\mu_n)$ converges to $\mu$ as $n\to\infty$. In particular, $\mu\in\cone(\lambda)$.
 
  To conclude that this construction defines a map from $\outenv(X,\lambda)$ to $\Theta^{out}(\lambda)$, it remains to show that the resulting measured lamination $\mu$ is independent of the choice of the harmonic map $f_\infty:X_\infty\to X$. To see this, let $f_\infty':X_\infty'\to X\setminus\lambda$ be a  harmonic diffeomorphism from a punctured surface $X_\infty'$ which also defines $\HSL(X,Z)$.  Let  $\mu'\in\far_\lambda$ be the measured lamination obtained from the \emph{horizontal foliation} of the Hopf differential of $f_\infty'$ by removing all half-planes and half-infinite cylinders.   By \cite[Proposition 12.9]{PanWolf2022}, there exists a surjective conformal map $\psi:X_\infty'\to X_\infty$ such that $f_\infty'=f_\infty\circ\psi$.  In particular, $\psi$ sends the horizontal foliation of the Hopf differential of $f'_\infty$ to the horizontal foliation of the Hopf differential of $f_\infty$. This implies that $\mu=\mu'$.

 In summary, we arrive at a well-defined map  from $\outenv(X,\lambda)$ to $\Theta^{out}(\lambda)$:
 \begin{equation*}
 	\widehat\Psi_{out}: \outenv(X,\lambda)\longrightarrow \Theta^{out}(\lambda)
 \end{equation*}
 which sends $Z\in\outenv(X,\lambda)$ to the pair $(\mu,\dth(X,Z))$.  Sending $X$ to the cone point of $\Theta^{out}(\lambda)$, we may extend $\widehat\Psi_{out}$ to $$\Psi_{out}:\{X\}\cup\outenv(X,\lambda)\to\Theta^{out}(\lambda).$$ Similarly, one can define a map  $\Psi_{in}:\{Y\}\cup\inenv(Y,\lambda)\to\Theta^{in}(\lambda)$. 
 \end{proof}
 The remainder of this section is devoted to proving the next theorem.
 
   \begin{theorem}\label{thm:out:in:env}
	For any $X\in\T(S)$ and any chain recurrent geodesic lamination $\lambda$,  both $\Psi_{in}$ and $\Psi_{out}$ are homeomorphisms. 
\end{theorem}
 
As a consequence, we see that

\begin{corollary}\label{cor:env:dim}
	The dimension of $\env(X,Y)$ is given by $$\dim \env(X,Y)=1+\dim(\cone(\lambda)).$$\end{corollary}

As a corollary of Theorem~\ref{thm:out:in:env}, we will obtain Corollary~\ref{cor:phs:ps:characterization} a characterization of harmonic stretch lines that does not explicitly involve a limiting process of maps.

 \subsection{From $\Theta^{out}(\lambda)$ to $\outenv(X,\lambda)$ and from $\Theta^{in}$ to $\inenv(Y,\lambda)$} \label{subsec:theta:hs}

In this subsection, we shall assign to any measured lamination in $\cone(\lambda)$ a specific harmonic stretch line that passes through $X$.  \begin{definition}\label{def:modified:ml}
     Let $\gamma\in  \far_\lambda$ and let $\alpha$ be an ideal geodesic boundary arc of $X\backslash\lambda$ or a simple closed geodesic boundary of $X\backslash\lambda$ that is not spiralled by $\gamma$. We formally add the infinitely weighted component $\infty\cdot \alpha$ to $\gamma$. 
     The resulting measured lamination $ \mathfrak{m}(\gamma)$ is called \emph{the modified measured lamination} of $\gamma$. To clarify:  for any $\gamma\in\far_\lambda$,   the restriction to $X\setminus\lambda$ of any arc in $X$ which is transverse to $\lambda$ has infinite intersection number with the modified lamination $\mathfrak m(\gamma)$. 
 \end{definition}\label{rmk:modified:ml}

      \begin{remark}
      	(i) From the point of view of measured foliations,  the measured foliation corresponding to 
      	$\mathfrak{m}(\gamma)$ is obtained from the measured foliation corresponding to $\gamma$ by adding half-infinite cylinders for some simple geodesic boundary components of $X\backslash\lambda$ and by adding half-planes for all ideal geodesic boundary components of $X\backslash\lambda$. The modified measure lamination $\mathfrak{m}(\gamma)$ here corresponds to an \emph{admissible measured foliation}  defined in \cite[Definition 5.2]{PanWolf2022}.

        (ii) The reader will note a resemblance in the definition here of a \enquote{modified measured lamination} and the \enquote{generalized measured lamination} introduced in Remark~\ref{rmk:union:cone}.  The distinction between the two is that the generalized measured lamination is defined on the {\it closed} surface $X$, while the modified measured lamination is defined on the {\it crowned} surface $X \setminus \lambda$.
      \end{remark}

\begin{construction}
    \label{phsl:construction} For any $\gamma\in\far_\lambda$,
      by \cite[Theorem 1.13]{PanWolf2022}, there exists a unique punctured Riemann surface $X_\infty$ homeomorphic to $X\setminus\lambda$ and a  harmonic diffeomorphism $f_\infty: X_\infty\to X\setminus\lambda$ such that the horizontal foliation of the Hopf differential of $f_\infty$ is equivalent to the modified measured lamination  $\mathfrak m(\gamma)$. Combined with \cite[Theorem 1.5]{PanWolf2022}, this defines a piecewise harmonic stretch line $\PHSL(X,\lambda,\gamma)$ in $\T(S)$ with data $(X,\lambda,
      \gamma)$. 
      \end{construction}

The following lemma strengthens the conclusion of this construction when $\gamma\in\cone(\lambda)$.

  \begin{lemma}\label{lem:hs}   For any $\gamma\in\cone(\lambda)$, the piecewise harmonic stretch line $\PHSL(X,\lambda,\gamma)$ is a harmonic stretch line, denoted by $\HSL(X,\lambda,
  \gamma)$. 
         \end{lemma}

   \begin{proof}   Let $\{\gamma_n\}$ be as in Lemma \ref{lem:image:cut}(ii).   
       By \cite[Theorem 3.1]{Wolf1998}, there exists a unique Riemann surface $X_n$ such that the horizontal foliation of the Hopf differential of the (unique) harmonic map $X_n\to X$ is equivalent to $\gamma_n$. Let $\HR(X_n,X)$ be the harmonic map ray comprising those $X_{n,t}$ satisfying
       \begin{itemize}
       	\item $X=X_{n,1}$ for all $n$;
       	\item the horizontal foliation of the Hopf differential of the (unique) harmonic map $X_n\to X_{n,t}$ is $\sqrt{t}\gamma_n$.
       \end{itemize}
         By \cite[Theorem 1.3]{PanWolf2022}, there exists a subsequence of $\{\gamma_n\}$, still denoted by $\{\gamma_n\}$, such that the sequence of harmonic map rays $\HR(X_n,X)$ not only subconverges, but converges locally uniformly to some harmonic stretch line $R'_\lambda(X)$. Moreover, by \cite[Lemma 4.5 and Proposition 8.1]{PanWolf2022}, we see that  $R'_\lambda(X)$  is defined by some  harmonic diffeomorphism $f'_\infty:X_\infty'\to X\setminus\lambda$ which is the limit of the sequence of the harmonic maps $f_n:X_n\to X$ in the sense of \cite[Definition 4.5]{PanWolf2022}. This implies that the horizontal foliation of the Hopf differential of $f'_\infty$ is the limit of the horizontal foliation of the Hopf differential of $f_n:X_n\to X$, which is exactly $\mathfrak m(\gamma)$. Recalling the map $f_\infty: X_\infty \to Y\setminus \lambda$ that defined $\PHSL(X,\lambda,\gamma)$, it then follows from \cite[Theorem 1.13]{PanWolf2022} that $f'_\infty:X_\infty'\to X\setminus\lambda$ and $f_\infty:X_\infty\to X\setminus\lambda$ are equivalent, i.e. there exists a conformal map $\eta:X_\infty\to V'_\infty$ such that $f_\infty=f'_\infty\circ \eta$. Combined with \cite[Lemma 12.8]{PanWolf2022}, this implies  that $\PHSL(X,\lambda,\gamma)$ and $R'_\lambda(X)$ coincide. In particular, $\PHSL(X,\lambda,\gamma)$ is a harmonic stretch line. 
     \end{proof}

            As a direct consequence of Lemma \ref{lem:hs},  we have the following. 
          \begin{corollary}\label{cor:Psi:surjection}
          	For any $X,Y\in\T(S)$ and any chain recurrent geodesic lamination $\lambda$, both $\Psi_{out}:\{X\}\cup\outenv(X,\lambda)\to\Theta^{out}(\lambda)$ and $\Psi_{in}:\{Y\}\cup\inenv(Y,\lambda)\to \Theta^{in}(\lambda)$ are surjective. 
          \end{corollary}
          \begin{proof}
          For any $\gamma\in\cone(\lambda)$,   consider the harmonic stretch line $\HSL(X,\lambda, \gamma)$    obtained from Lemma \ref{lem:hs}. Notice that $X$ cuts $\HSL(X,\lambda, \gamma)$ into two half-infinite rays. We denote  by $\HSR(X,\lambda,\gamma)$ the one starting at $X$.  For any pair $(\gamma,t)\in\Theta^{out}(\lambda)$, consider the point  $Z_{\gamma,t}$ on the harmonic stretch ray $\HSR(X,\lambda,\gamma)$ whose distance from $X$ is $t$. It is clear that $\Psi_{out}(Z_{\gamma,t})=(\gamma,t)$.  This shows that $\Psi_{out}$ is surjective. That $\Psi_{in}$ is surjective follows similarly. 
          \end{proof}

\begin{proof}[Proof of Theorem \ref{thm:out:in:env}] 
That $\Psi_{in}$ and $\Psi_{out}$ are both injective follows from the uniqueness part of \cite[Theorem 1.13]{PanWolf2022}. To see this, note that each $\gamma\in\cone(\lambda)$ defines a real tree $T_\gamma$, which is the leaf space of the lift to the universal cover $\widetilde{X\backslash\lambda}$ of the measured foliation corresponding to the modified measured lamination $\mathfrak{m}(\gamma)$ of $\gamma$, and an admissible orientation-preserving boundary correspondence. The uniqueness part of \cite[Theorem 1.13]{PanWolf2022} implies that there exists at most one $\pi_1(X\backslash \lambda)$-equivariant minimal graph in $\widehat{Y\backslash\lambda}\times T_\gamma$. Equivalently, this means that there exists at most one pair $(X_\infty,f_\infty)$ of a (possibly disconnected) punctured Riemann surface $X_\infty$ and a harmonic diffeomorphism $f_\infty:X_\infty\to X\setminus\lambda$ whose Hopf differential has $\mathfrak{m}(\gamma)$ as the horizontal foliation. Since every harmonic stretch line is constructed by such a pair, we see that there exists at most one harmonic stretch line through $X$ whose defining modified measured lamination is $\mathfrak m(\gamma)$. This proves the uniqueness of $\Psi_{in}$ and $\Psi_{out}$.   Corollary \ref{cor:Psi:surjection} asserts
that $\Psi_{in}$ and $\Psi_{out}$ are both surjective. 
The continuity of $\Psi_{in},\Psi_{out}$ and their inverse maps follow from Proposition \ref{prop:PanWolf2022}.
\end{proof}

We are now in a position to prove Theorem~\ref{thm:phs:hs}, as promised after Corollary~\ref{cor:env:dim}.

	\begin{corollary}[Characterization of harmonic stretch lines]\label{cor:phs:ps:characterization}
		Let $X\in\T(S)$. Let $\lambda$ be a chain recurrent geodesic lamination and $\gamma\in\far_\lambda$. Then the piecewise harmonic stretch line with data $(X,\lambda,\gamma)$  is a harmonic stretch line if and only if $\gamma\in\cone(\lambda)$.
	\end{corollary}

\begin{proof}[Proof of Corollary \ref{cor:phs:ps:characterization}] Let $\PHSL(X,\lambda,\gamma)$ be the piecewise harmonic stretch line defined by the triple $(X,\lambda,\gamma)$. For one direction, suppose that  $\PHSL(X,\lambda,\gamma)$ is a harmonic stretch line. From the construction of the map $\widehat{\Psi}_{out}:\outenv(X,\lambda)\to\Theta^{out}(\lambda)$, we see that $\gamma\in\cone(\lambda)$. The converse direction is the content of Lemma \ref{lem:hs}.
\end{proof}
 The precise relationship between harmonic stretch lines and piecewise harmonic stretch lines remains unclear.  To begin, the following basic question is unknown.
	
\begin{question}
	Does there exist any piecewise harmonic stretch line whose support, as a subset of $\T(S)$, is different from the support of any harmonic stretch line?
\end{question}

\begin{remark} 
   Some explanation about the question may be needed. Notice that, by the item (i) of Theorem \ref{thm:env:shape}, for any $X\in\T(S)$ and for any maximally CR geodesic lamination $\lambda$, the out-envelope $\outenv(X,\lambda)$ is a single harmonic stretch ray starting at $X$. This implies that,  for  any geodesic lamination $\lambda'$ containing $\lambda$ and for any measured geodesic lamination $\gamma'\in \far_{\lambda'}$, the  piecewise harmonic stretch line $\PHSL(X,\lambda',\gamma')$ with  data $(X,\lambda',\gamma')$   -- as a set in the \tec space--  coincides with the harmonic stretch line which contains $\outenv(X,\lambda)$. The question above asks if this agreement between piecewise harmonic stretch lines and harmonic stretch lines always holds.
\end{remark}

\subsection{Harmonic stretch line and extended measured laminations}
Recall that the generalized measured lamination of $\gamma\in\cone(\lambda)$ is the formal sum $ \gamma+\infty\cdot \lambda$ (Remark \ref{rmk:union:cone}). Every harmonic stretch line through $X$ is determined by the triple $(\lambda,\gamma)$  with $\gamma\in \cone (\lambda)$, so is equivalently determined by the generalized measured lamination $\gamma+\infty\cdot\lambda$. Given $X\in\T(S)$ and  a measured geodesic lamination $\lambda$ on $X$, let $\HR_{X,\lambda}:[0,\infty)\to\T(S)$ be the harmonic map ray such that $X=\HR_{X,\lambda}(1)$ and that the Hopf differential of the harmonic diffeomorphism from $\HR_{X,\lambda}(0)$ to $\HR_{X,\lambda}(t)$ is $\sqrt{t}\lambda$ (see \cite[Theorem 3.1]{Wolf1998} for the existence and uniqueness of $\HR_{X,\lambda}(0)$). 
\begin{proposition}\label{hr:hsl}Let $X\in\T(S)$ and $\lambda_n$ be a sequence measured geodesic lamination on $X$. Then the sequence of harmonic map rays $\HR_{X,\lambda_n}$ converges to a harmonic stretch line $\HSL(X,\lambda,
\gamma)$ with $\lambda$ being a chain recurrent geodesic lamination and $\gamma\in\cone(\lambda)$ if and only if $\lambda_n$ converges to the generalized measured lamination $\gamma+\infty\cdot\lambda$ in the following sense:     \begin{itemize}
        \item  for any geodesic arc $k$ transverse to $\lambda$, the intersection number satisfies $i(\lambda_n,k)\to\infty$;
        \item  for any compact geodesic arc $k$ on $X\backslash\lambda$ that is transverse to $\gamma$, the intersection number satisfies $i(\lambda_n,k)\to i(\gamma,k)$;
        \item for any simple closed geodesic $\alpha$ of $\lambda$ and for sufficiently large $n$, we have $\lambda_n$ wraps around $\alpha$ in the same direction as $\gamma$. 
    \end{itemize}
\end{proposition}
\begin{proof}
    (i) Suppose that $\HR_{X,\lambda_n}$ converges to the harmonic stretch line $\HSL(X,\lambda,\gamma)$. Let $X_n=\HR_{X,\lambda_n}(0)$. Then by \cite[Proposition 8.1]{PanWolf2022}, the there exists a subsequence, still denoted by $X_n$ for simplicity, such that the harmonic map $f_n:X_n\to X$ (homotopic to the identity) converges to a harmonic diffeomorphism $f_\infty:X_\infty\to X\backslash\lambda$ for some (possibly disconnected) Riemann surface $X_\infty$ in the sense of \cite[Definition 4.5]{PanWolf2022}. Let $\Phi_n$ and $\Phi_\infty$ be respectively the Hopf differential of $f_n$ and $f_\infty$.
    Let $\mathfrak{p}_n$ and $\mathfrak{p}$ be respectively the set of zeros of $\Phi_n$ and $\Phi_\infty$.  Then, by \cite[Lemma 4.1 and the first paragraph in Section 4.2]{PanWolf2022}, the pointed singular flat surface $(|\Phi_n|,\mathfrak{p}_n)$ subconverges to the pointed  singular flat  surface $(|\Phi_\infty|,\mathfrak{p})$.  
    Consider the horizontal foliation of $\Phi_n$,  which is $\lambda_n$ by assumption. We have, 
    \begin{itemize}
        \item  for any geodesic arc $k$ transverse to $\lambda$, the intersection number satisfies $i(\lambda_n,k)\to\infty$;
        \item  for any compact geodesic arc $k$ on $X\backslash\lambda$ that is transverse to the horizontal foliation $\gamma$ of $\Phi_\infty$, the intersection number satisfies $i(\lambda_n,k)\to i(\gamma,k)$;
        \item for any simple closed geodesic $\alpha$ of $\lambda$ and for sufficiently large $n$, we have $\lambda_n$ wraps around $\alpha$ in the same direction as $\gamma$. 
    \end{itemize}
    This proves that $\lambda_n\to \gamma+\infty\cdot \lambda$.

    (ii) Now we turn to the inverse direction. Suppose that $\lambda_n\to \gamma+\infty\cdot \lambda$. By \cite[Theorem 1.3,  Proposition 8.1, and Definition 12.1]{PanWolf2022}, any sequence of $\HR_{X,\lambda_n}$ contains a subsequence that converges to some harmonic stretch line $\HSL(X,\lambda',\gamma')$. By the discussion in the preceding paragraph and the assumption that $\lambda_n\to\gamma+\infty\cdot\lambda$, we see that $\gamma+\infty\cdot \lambda=\gamma'+\infty\cdot\lambda'$. Hence $\HSL(X,\lambda',\gamma')=\HSL(X,\lambda,\gamma)$. The arbitrariness of the convergent subsequence of $\HR_{X,\lambda_n}$ then implies that $\HR_{X,\lambda_n}$ convergs to $\HSL(X,\lambda,\gamma)$. This completes the proof.
\end{proof}

\begin{proposition}\label{prop:hsl:extended}
   Let $\lambda_n$ be a chain recurrent geodesic lamination on $X\in\T(S)$ and $\gamma_n\in \cone(\lambda_n)$. Then the sequence of harmonic stretch lines $\HSL(X,\lambda_n,\gamma_n)$  converges to $\HSL(X,\lambda_0,\gamma_0)$ if and only if $\gamma_n+\infty\cdot \lambda_n\to \gamma_0+\infty\cdot \lambda_0$. 
\end{proposition}
\begin{proof}
   Note that there are only countably many homotopy classes of proper geodesics on $X\backslash\lambda_0$.  The proposition then follows from Proposition \ref{hr:hsl} by a diagonal argument. 
\end{proof}

\color{black}

\subsection{Continuity of out-envelopes and in-envelopes}\label{sec:continuity:out:in:cpt}
In this subsection, we consider the continuity of out-envelopes and in-envelopes. The convexity of balls under the Thurston metric is currently unknown. The intersection between a Thurston geodesic and a ball may have more than one component. To simplify the discussion, we need a definition. Let $X\in\T(S)$ and let $\Omega\subset\T(S) $ be a subset of $\T(S)$ that contains $X$, define the \emph{star} $\str(\Omega,X)$ of $\Omega$ centered at $X$ to be the subset of points $Z\in\Omega$ such that the harmonic stretch segment $[X,Z]$ is also contained in $\Omega$:
\begin{equation*}
    \str(\Omega,X):=\{Z\in\Omega:[X,Z]\subset \Omega\}.
\end{equation*} 
\color{black}
\begin{theorem}\label{thm:out:in:continuity}
	Let $X,X_n\in\T(S)$ with $X_n\to X$. Let $\lambda,\lambda_n$ be chain recurrent geodesic laminations such that $\lambda_n\to\lambda$ with respect to the Hausdorff topology. Then for any closed ball $\mathscr K\subset \T(S)$ centered at $X$ of positive radius,  we have
	\begin{equation*}		{\str\left(\overline{\outenv(X_n,\lambda_n)}\cap \mathscr K),X_n\right)\longrightarrow {\str}\left(\overline{\outenv(X,\lambda) }\cap \mathscr K,X\right) }
		\end{equation*}
		and
		\begin{equation*}
		{ \str\left(\overline{\inenv(X_n,\lambda_n)}\cap \mathscr K,X_n\right)\longrightarrow \str\left(\overline{\inenv(X,\lambda) }\cap \mathscr K,X\right)}
	\end{equation*}
	with respect to the Hausdorff topology, as $n\to\infty$.
\end{theorem}
\begin{proof}
	Since the proofs for out-envelopes and in-envelopes are similar, we shall prove the conclusion for out-envelopes.  To this end, it suffices to prove the following two claims.
	\begin{enumerate}
		\item  For any convergent sequence  $[X_n,Z_n]\subset \overline{\outenv(X_n,\lambda_n) }\cap \mathscr K$, the limit is contained in $ \overline{\outenv(X,\lambda) }\cap \mathscr K$.
		\item Any segment  $[X,Z]\subset \overline{\outenv(X,\lambda) }\cap \mathscr K$ can be approximated by some sequence $[X_n,Z_n]\subset \overline{\outenv(X_n,\lambda_n) }\cap \mathscr K$.
	\end{enumerate}
	
	For claim (i), let $Z$ be the limit of $Z_n$. By Proposition \ref{prop:PanWolf2022}, $[X_n,Z_n]$ converges to $[X,Z]$. It then follows from Theorem \ref{thm:maximally:stretched:lamination} that the maximally stretched lamination along $[X,Z]$ contains the Hausdorff limit of $\lambda_n$, which is $\lambda$ by assumption. Hence $[X,Z]\subset\overline{\outenv(X,\lambda)}$. That $[X,Z]\subset \mathscr K$ follows from the assumption that $[X_n,Z_n]\subset \mathscr K$. This proves claim (i).
	
	For claim (ii), it suffices to consider the case where $Z\in\outenv(X,\lambda)$. Let $Z$ be such a point. By Theorem \ref{thm:out:in:env}, there exists    	$\gamma\in \cone(\lambda)$ which defines the harmonic stretch line $\HSL(X,Z)$ proceeding from $X$ through $Z$.  
	
	 Since $\lambda_n\to \lambda$ with respect to the Hausdorff topology, it follows  that there exists $N>0$  such that $\cut_{\lambda_n}:U_\delta(\lambda)\to \far_{\lambda_n}$ is well defined for $n>N$ for $\delta$ sufficiently small and fixed. (To see this, note that $\lambda_n$ is $C^1$ close to $\lambda$, so there exists $N>0$ such that for $n>N$, any $\gamma\in U_\delta(\lambda)$ belongs to $U_{2\delta}(\lambda_n)$. In particular, $\cut_{\lambda_n}$ is well defined on $U_\delta(\lambda)$.)  By 
    the second part of Lemma \ref{lem:image:cut}, there exists $\gamma_n\in U_\delta(\lambda)$    such that $\cut_\lambda(\gamma_n)=\gamma$ and that $\supp(\gamma_n)\to \lambda\cup \supp(\gamma)$. 
        For each $\gamma_n$  with $n>N$,  by Lemma \ref{lem:hs}, the measured lamination $\cut_{\lambda_n}(\gamma_n)$, which belongs to $\cone(\lambda_n)$,  defines a harmonic stretch line $\HSL(X,\lambda_n,\cut_{\lambda_n}(\gamma_n))$ maximally stretching exactly along $\lambda_n$. Take $Y_n\in\HSL(X,\lambda_n,\cut_{\lambda_n}(\gamma_n))$ so that $\HSL(X,\lambda_n,\cut_{\lambda_n}(\gamma_n))$ proceeds from $X$ through $Y_n$ and that $\dth(X,Y_n)=1$. Then $\{Y_n\}$ is contained in the closed unit ball centered at $X$, which is compact. Therefore, up to a subsequence, we may assume that $Y_n$ converges to some point $Y\in\T(S)$.  By Proposition \ref{prop:PanWolf2022},  we see that  $\HSL(X,\lambda_n,\cut_{\lambda_n}(\gamma_n))$  converges to some harmonic stretch line $R$ which proceeds from $X$ through $Y$.  By Theorem \ref{thm:maximally:stretched:lamination}, the maximally stretched lamination $\lambda'$ of $R$ contains the Hausdorff limit of $\lambda_n$; by assumption, that limit is $\lambda$. Again by Theorem \ref{thm:out:in:env}, the harmonic stretch line $R$ is defined by a measured lamination $\gamma'\in \cone(\lambda')$.  Hence,  we see that $\gamma'\subset X\backslash\lambda'\subset X\backslash\lambda$. In particular, this means that $\gamma'\in\far_\lambda$.

         On the other hand, the construction of $\gamma_n$ and $\cut_{\lambda_n}$ implies that
         \begin{itemize}
         	\item[(a)] $\lambda_n \cup \supp(\cut_{\lambda_n}(\gamma_n))$ converges to $\lambda\cup \supp(\gamma)$ with respect to the Hausdorff topology;  
            \item[(b)] for any geodesic arc $k$ on $X$ that is transverse to $\lambda$ we have $i(\cut_{\lambda_n}(\gamma_n)+\infty\cdot \lambda_n,k)\to \infty$,
         	\item[(c)] for any compact geodesic arc $k\subset X\setminus\lambda$, we have $i(\cut_{\lambda_n}(\gamma_n)+\infty \cdot \lambda_n,k)\to i(\gamma,k)$;
            \item[(d)] for any simple closed geodesic $\alpha$ of $\lambda$ and for sufficiently large $n$, the lamination $\gamma_n$ wraps along $\alpha$ in the same direction as $\gamma$.
         \end{itemize}
From this, we infer that $\cut_{\lambda_n}(\gamma_n)+\infty\cdot \lambda_n\to \gamma+\infty\cdot\lambda$. On the other hand, combining  Proposition \ref{prop:hsl:extended} and the assumption that $\HSL(X,\lambda_n,\cut_{\lambda_n}(\gamma_n))$ converges to $R$, we see that $\cut_{\lambda_n}(\gamma_n)+\infty\cdot \lambda_n\to \gamma'+\infty\cdot\lambda'$. Hence $ \gamma+\infty\cdot\lambda=\gamma'+\infty\cdot\lambda'$.
   In particular, the harmonic stretch line $R$  coincides with $\HSL(X,Z)=\HSL(X,\lambda,\gamma)$. This completes the proof of claim (ii). 
\end{proof}

For the convenience of reference, we also need the following corollary. Recall that (Section \ref{subsec:hausdorff}) given a subset $M\subset \T(S)$ and a positive constant $\epsilon$, we denote by $\mathcal{N}_\epsilon(M)$  the $\epsilon$ neighbourhood of $M$ under the symmetric metric $\dth^S$ of the Thurston metric.
\begin{corollary}\label{cor:out:in:approximation}
	Let $X,X_n\in\T(S)$ with $X_n\to X$. Let $\lambda,\lambda_n$ be chain recurrent geodesic laminations such that $\lambda$ contains the Hausdorff limit of any convergent subsequence of $\{\lambda_n\}$ with respect to the Hausdorff topology. Then for any closed ball $\mathscr K\subset \T(S)$ centered at $X$ of positive radius and any $\epsilon>0$, there exists $ D=(X,\lambda,\mathscr K,\epsilon)$ such that for $n>D$, 
	\begin{equation*}
			{\str\left(\overline{\outenv(X,\lambda) }\cap \mathscr K,X\right)\quad\subset 	\quad \mathcal N_\epsilon\left(\str\left(\overline{\outenv(X_n,\lambda_n)}\cap \mathscr K,X_n\right) \right)}
		\end{equation*}
		and
		\begin{equation*}
{\str\left(\overline{\inenv(X,\lambda) }\cap \mathscr K,X\right)\quad\subset 	\quad \mathcal{N}_\epsilon\left(\str\left(\overline{\inenv(X_n,\lambda_n)}\cap \mathscr K,X_n\right)\right).}
	\end{equation*}
\end{corollary}
\begin{proof}

	Notice that $\lambda$  comprises finitely many minimal sublaminations and finitely many isolated leaves. In particular, $\lambda$ has finitely many sublaminations, say $\lambda^1,\cdots,\lambda^k$.  In view of this, we can divide  $\{n\}_{n\geq 1}$ into  $k$ infinite subsequences  $\{n^i_j\}_{j\geq 1}$ with $1\leq i\leq k$ such that for each $1\leq i\leq k$,  we have $\lambda_{n^i_j}\to \lambda^i$ as $j\to\infty$. 
    By Theorem \ref{thm:out:in:continuity},  we have 
	\begin{equation*}
			\str\left(\overline{\outenv(X,\lambda^i) }\cap \mathscr K,X\right)\quad = 	\quad \lim_{j\to\infty}\str\left(\overline{\outenv(X_{n^i_j},\lambda_{n^i_j})}\cap \mathscr K,X_{n^i_j}\right)
		\end{equation*}
		and
		\begin{equation*}
\str\left(\overline{\inenv(X,\lambda^i) }\cap \mathscr K,X\right)\quad= 	\quad \lim_{j\to\infty}\str\left(\overline{\inenv(X_{n^i_j},\lambda_{n^i_j})}\cap \mathscr K,X_{n^i_j}\right).
	\end{equation*}	
    Then, for any $1\leq i\leq k$ and any $\epsilon>0$, there exists $D^i=D^i(X,\lambda,\mathscr K,\epsilon,i)$ such that for $j>D^i$, 
	\begin{equation*}
			\str\left(\overline{\outenv(X,\lambda^i) }\cap \mathscr K,X\right)\quad\subset 	\quad \mathcal N_\epsilon\left(\str\left(\overline{\outenv(X_{n^i_j},\lambda_{n^i_j})}\cap \mathscr K,X_{n^i_j}\right) \right)
		\end{equation*}
		and
		\begin{equation*}
\str\left(\overline{\inenv(X,\lambda^i) }\cap \mathscr K,X\right)\quad\subset 	\quad \mathcal{N}_\epsilon\left(\str\left(\overline{\inenv(X_{n^i_j},\lambda_{n^i_j})}\cap \mathscr K,X_n\right)\right).
	\end{equation*}
	Since $\lambda^i$ is a sublamination of $\lambda$, we see that  
		\begin{equation*}
			\overline{\outenv(X,\lambda) }\subset \overline{\outenv(X,\lambda^i) },\qquad 
			\overline{\inenv(X,\lambda) }\subset \overline{\inenv(X,\lambda^i)}.
		\end{equation*} 
        Therefore, for any $j>D^i$, we have
        \begin{equation*}
			\str\left(\overline{\outenv(X,\lambda) }\cap \mathscr K,X\right)\quad\subset 	\quad \mathcal N_\epsilon\left(\str\left(\overline{\outenv(X_{n^i_j},\lambda_{n^i_j})}\cap \mathscr K,X_{n^i_j}\right) \right)
		\end{equation*}
		and
		\begin{equation*}
\str\left(\overline{\inenv(X,\lambda) }\cap \mathscr K,X\right)\quad\subset 	\quad \mathcal{N}_\epsilon\left(\str\left(\overline{\inenv(X_{n^i_j},\lambda_{n^i_j})}\cap \mathscr K,X_{n^i_j}\right)\right).
\end{equation*} The finiteness  of $1\leq i\leq k$ then implies that there exists $D=D(X,\lambda,\mathscr K, \epsilon)$ such that \begin{equation*}
			\str\left(\overline{\outenv(X,\lambda) }\cap \mathscr K,X\right)\quad\subset 	\quad \mathcal N_\epsilon\left(\str\left(\overline{\outenv(X_n,\lambda_n)}\cap \mathscr K,X_{n^i_j}\right) \right)
		\end{equation*}
		and
		\begin{equation*}
\str\left(\overline{\inenv(X,\lambda) }\cap \mathscr K,X\right)\quad\subset 	\quad \mathcal{N}_\epsilon\left(\str\left(\overline{\inenv(X_n,\lambda_n)}\cap \mathscr K,X_n\right)\right).
\end{equation*}
This proves the lemma.
\color{black}
	\end{proof}

\begin{remark}
 The continuity of envelopes does not follow directly from Proposition \ref{prop:DLRT2020} and Theorem \ref{thm:out:in:continuity}, because in general $\lim\limits_{n\to\infty} A_n=A$ and $\lim\limits_{n\to\infty} B_n=B$ does not imply $\lim\limits_{n\to\infty} A_n\cap B_n=A\cap B$.  
 \end{remark}

\subsection{Proof of Theorem \ref{thm:structure:out:in}}
\begin{proof}[Proof of Theorem \ref{thm:structure:out:in}]
Theorem  \ref{thm:structure:out:in}(i) follows from Theorem \ref{thm:out:in:env}, and Theorem \ref{thm:structure:out:in}(ii)  follows from Theorem \ref{thm:out:in:continuity}. 
\end{proof}

 \section{Extendability}\label{sec:extendability}
In this section, we shall introduce the notion of \emph{extendability} of points in envelopes, and describe the shapes of envelopes as \emph{cones over cones}, identifying the ancestor object in terms of the initial and terminal points of the envelope.

\subsection{Extendability} 
Recall that  any pair of distinct points $(X,Z)$ determines a harmonic stretch line $\HSL(X,Z)$ proceeding from $X$ through $Z$.  Notice that $X$ cuts $\HSL(X,Z)$ into a harmonic stretch half-infinite segment and a harmonic stretch ray: the segment ends at $X$, and the ray starts at $X$.
\begin{definition}
    We denote by $\HSR(X,Z)\subset\HSL(X,Z)$ the harmonic stretch ray that starts at $X$, and by $\HSL^-(X,Z)$ the harmonic stretch ray that ends at $X$. 
\end{definition}

We start with the following lemma.
\begin{lemma}\label{lem:endpoints}
Let $X,Y \in \T(S)$ with $X\neq Y$. For any $Z\in\env(X,Y)\backslash\{X\}$, the intersection 
	 $\HSL(X,Z)\cap\env(X,Y)$ is an oriented segment with initial point $X$.  For any $Z\in\env(X,Y)\backslash\{Y\}$, the intersection  $\HSL(Z,Y)\cap\env(X,Y)$ is an oriented segment with terminal point $Y$.
\end{lemma}
\begin{proof}
Note that a harmonic stretch line is also a piecewise harmonic stretch line. The lemma then follows directly from Lemma \ref{lem:concatenation}. \end{proof}

  \begin{definition}[Extendability]\label{def:extendability}
 Let $X,Y \in\T(S)$ with $X\neq Y$ and $Z\in\env(X,Y)$.  
 \begin{itemize}
 	\item The point $Z$ is said to be right extendable if $\HSL(X,Z)\cap \env(X,Y)$ strictly contains $[X,Z]$. Otherwise, the point $Z$ is said to be a right boundary point.
 	\item The point $Z$ is said to be left extendable if there exists $\HSL(Z,Y)\cap \env(X,Y)$ strictly contains $[Z,Y]$. Otherwise, the point $Z$ is said to be a left boundary point.
 	\item The point $Z$ is said to be  bi-extendable if it is both left extendable and right extendable. 
 \end{itemize} 
  \end{definition}

 Equivalently, $Z\in\env(X,Y)$ is a right boundary point if $Z$ is the terminal point of the oriented segment $\HSL(X,Z)\cap\env(X,Y)$.  Similarly, $Z\in\env(X,Y)$ is a left boundary point if $Z$ is the initial point of the oriented segment $\HSL(Z,Y)\cap\env(X,Y)$. 
	Let $\romega(X,Y)$ and  $\lomega(X,Y)$ be respectively the set of right boundary points and the set of left boundary points of $\env(X,Y)$.  Let $\bnd(X,Y)=\lbd(X,Y)\cap \rbd(X,Y)$. In particular, we have $$X\in\lbd(X,Y)\backslash\bnd(X,Y) \text{ and } Y\in\rbd(X,Y)\backslash\bnd(X,Y).$$ We will show in Lemma \ref{lem:env:top} that $\env(X,Y)$ is a cone over $\rbd(X,Y)$ as well as a cone over $\lbd(X,Y)$, via harmonic stretch lines.

	For the punctured torus case, Dumas-Lenzhen-Rafi-Tao \cite{DLRT2020} proved the following (in terms of the current terminology):
	\begin{itemize}
		\item if $\Lambda(X,Y)$ is a simple closed curve, then $\env(X,Y)$ is geodesic quadrilateral region with 	 $\bnd(X,Y)\cup\{X,Y\}$ being the set of vertices and  $\rbd(X,Y)\cup\lbd(X,Y)$ being the sides of the quadrilateral;
		\item if $\Lambda(X,Y)$ is not a simple closed curve, then $\env(X,Y)$ is a geodesic with $\bnd(X,Y)=\emptyset$, $\lbd(X,Y)=\{X\}$, and $\rbd(X,Y)=\{Y\}$.
	\end{itemize}
Here we will generalize the above description. We begin with the following observation, characterizing extendability in terms of the maximally stretched lamination.

 	\begin{lemma}[Characterization of extendablity]\label{lem:extendability}
   Let  $X,Y\in\T(S)$ with $X\neq Y$ and $Z\in\env(X,Y)\backslash\{X,Y\}$.  Then 
\begin{enumerate}
	\item $Z$ is right extendable if and only if  $\Lambda(Z,Y)=\Lambda(X,Y)$,
	\item  $Z$ is left extendable if and only if   $ \Lambda(X,Z)=\Lambda(X,Y)$,
	\item  $Z$ is bi-extendable if and only if   $\Lambda(X,Z)=\Lambda(Z,Y)=\Lambda(X,Y)$.
\end{enumerate}
\end{lemma}
\begin{proof}
(i)   We start with the \enquote{only if} direction. Suppose   $Z$ is right extendable.  By the definition of right extendability, there exists $Z' \in \HSR(X,Z)\cap \env(X,Y)$ outside of $[X,Z]$, such that $[X,Z]\cup [Z,Z']\cup[Z',Y]$ is a geodesic from $X$ to $Y$.  In particular, $[Z,Z']\cup[Z',Y]$ is a geodesic from $Z$ to $Y$. Hence, by the definition of envelope and Corollary \ref{cor:env in terms of stretch lams}, we have $\Lambda(Z,Z')\supset\Lambda(Z,Y)$. Recall that, by Theorem \ref{thm:PanWolf2022:uniqueness}(ii), points on the same harmonic stretch ray share the same maximally stretched lamination. In particular, that $Z'\in\HSR(X,Z)$  implies that $\Lambda(Z,Z')=\Lambda(X,Z)$.  Hence, $\Lambda(X,Z)\supset\Lambda(Z,Y)$, which yields that $$\Lambda(Z,Y)=\Lambda(X,Z)\cap \Lambda(Z,Y) =\Lambda(X,Y),$$ 
	where the first equality follows from the inequality $\Lambda(X,Z)\supset\Lambda(Z,Y)$ while the second equality follows from Corollary~\ref{cor:env in terms of stretch lams} and the assumption that $Z\in\env(X,Y)$.

 For the \enquote{if} direction,  suppose that $\Lambda(Z,Y)=\Lambda(X,Y)$. Since, by assumption, $Z\in\env(X,Y)$, we see from Corollary~\ref{cor:env in terms of stretch lams}  that $\Lambda(X,Y)\subset \Lambda(X,Z)$. Therefore, $\Lambda(X,Z)\supset \Lambda(X,Y) = \Lambda(Z,Y)$.  Consider the harmonic stretch line $\HSL(X,Z)$, which maximally stretches exactly $\Lambda(X,Z)\supset\Lambda(Z,Y)$.  By Lemma \ref{lem:concatenation},  the intersection $\HSL(X,Z)\cap \env(Z,Y)$ contains a geodesic segment from $Z$ to some point $Z'$ distinct from $Z$. In particular,  this segment is contained in $\env(Z,Y)$, and hence also in $\env(X,Y)$ since  $\env(Z,Y)\subset \env(X,Y)$. This means that $Z$ is right extendable. 
	  	
	(ii) The proof is similar to the proof of (i).
	
	(iii) This follows from (i) and (ii). 
\end{proof}

The Lemma~\ref{lem:extendability} will be a principal tool in the proof of our results on the shape of the envelope. In particular, we may now prove one of the first such results, Theorem~\ref{thm:env:shape}(ii), which describes the envelope in terms of {\it the closure} of the set of surfaces which are both maximally stretched by $\Lambda(X,Y)$ from $X$ and also similarly stretched towards $Y$.

\begin{proposition}\label{prop:out-in-envelopes}
  	For any distinct  $X,Y\in\T(S)$, 
  	$$\env(X,Y)=\overline{\outenv(X,\Lambda(X,Y))\cap\inenv(Y,\Lambda(X,Y))}.$$
  \end{proposition}
 \begin{proof}
 Let $\lambda=\Lambda(X,Y)$. 
  It is clear that 
  $$\env(X,Y)\supset\overline{\outenv(X,\lambda)\cap\inenv(Y,\lambda)}:$$
  if $Z \in \outenv(X,\lambda) \cap \inenv(Y, \lambda)$, then $[X,Z] \cup [Z,Y]$ is a geodesic from $X$ to $Y$ (as each has stretch lamination $\lambda$), and $\env(X,Y)$ is a closed set.
  It remains to show the inverse direction.

 To begin, notice that by Corollary \ref{cor:env in terms of stretch lams} and Lemma \ref{lem:extendability}, we can decompose $\env(X,Y)$ as follows:
 \begin{eqnarray}
     &&\env(X,Y)\\
    \nonumber &=&(\outenv(X,\lambda)\cap \inenv(Y,\lambda))
     \cup (\lbd(X,Y)\backslash\bnd(X,Y))\\&&\nonumber \cup (\rbd(X,Y)\backslash\bnd(X,Y)) \cup \bnd(X,Y).
 \end{eqnarray}
Any  $Z\in \lbd(X,Y)\backslash\bnd(X,Y) $ is right extendable. Therefore, by Lemma \ref{lem:extendability}(i), we have $\Lambda(Z,Y)=\lambda$. Notice that any point $Z'\in [Z,Y]\backslash\{Z\}$ is left extendable. By Lemma \ref{lem:extendability}(ii), we have
$\Lambda(X,Z')=\lambda$.
Combining with the fact that $\Lambda(Z',Y)=\Lambda(Z,Y)=\lambda$, we see that $Z'\in\outenv(X,\lambda)\cap \inenv(Y,\lambda)$. The arbitrariness of $Z'$ in $[Z,Y]\backslash\{Z\}$ implies that $Z\in\overline{\outenv(X,\lambda)\cap\inenv(Y,\lambda)}$. This yields that
\begin{equation*}
   \lbd(X,Y)\backslash\bnd(X,Y)\subset \overline{\outenv(X,\lambda)\cap\inenv(Y,\lambda)}.
\end{equation*}
Similarly, we have\begin{equation*}
   \rbd(X,Y)\backslash\bnd(X,Y)\subset \overline{\outenv(X,\lambda)\cap\inenv(Y,\lambda)}.
\end{equation*}
  
  Now consider $Z \in \bnd(X,Y)$. By Lemma \ref{lem:extendability}(ii), we have that $\Lambda(X,Z)$ strictly contains $\lambda$. Then for any point $Z'$ in $[X,Z]\backslash\{Z\}$, since $\Lambda(X,Z')=\Lambda(X,Z)$ (cf. Theorem \ref{thm:PanWolf2022:uniqueness}(ii)),  we see that $\Lambda(X,Z')$ also strictly contains $\lambda$. So by Lemma \ref{lem:extendability}(ii), the point $Z'$ is not left extendable. On the other hand, $Z'$  is right extendable since it can be extended further to $Z$ along $[X,Z]$. Consequently, $Z'\in\lbd(X,Y)\backslash\bnd(X,Y)$. Letting $Z'$ approach $Z$ along $[X,Z]$,  we see that $Z$ can be approximated by $\lbd(X,Y)\backslash\bnd(X,Y)$. The arbitrariness of $Z\in\bnd(X,Y)$ then implies that
  $\bnd(X,Y)$ can be approximated by $\lbd(X,Y)\backslash\bnd(X,Y)$.  
   Therefore, $\env(X,Y) \subset \overline{\outenv(X,\lambda)\cap \inenv(Y,\lambda)} $. This completes the proof.
 \end{proof}

 We recall that a chain recurrent geodesic lamination is maximally CR if it is not contained in any other chain recurrent geodesic lamination.  
\begin{corollary}\label{cor:unique:geodesic}
	Let $X,Y\in\T(S)$ with $X\neq Y$. Then $\env(X,Y)$ is a geodesic if and only if $\Lambda(X,Y)$ is maximally CR.
\end{corollary}
\begin{proof}
	If $\Lambda(X,Y)$ is a maximally CR geodesic lamination, then by \cite[Corollary 2.3]{DLRT2020} (see also \cite[Proof of Theorem~8.5]{Thurston1986}), we have that $\env(X,Y)$ is a geodesic.
	
	We now assume that $\Lambda(X,Y)$ is not a maximally CR geodesic lamination. Let $\lambda'$ be a chain recurrent geodesic lamination which strictly contains $\Lambda$. Let $R$ be the harmonic stretch ray starting at $X$ which is determined by $\lambda'$ and $0\in\cone(\lambda')$ (Lemma \ref{lem:hs}).  By Lemma \ref{lem:concatenation}, we see that there exists some $Z\in R$ different from $X$ such that $[X,Z]\subset \env(X,Y)$. The uniqueness of harmonic stretch lines (Theorem \ref{thm:PanWolf2022:uniqueness}) and the assumption that $\lambda'$ strictly contains $\Lambda(X,Y)$ imply that $[X,Z]$ is not contained in the harmonic stretch segment $[X,Y]$. In particular, $\env(X,Y)$ strictly contains $[X,Y]$, which completes the proof. 
	\end{proof}

 	Define $\rpi: \env(X,Y)\backslash\{X\}\to \romega(X,Y)$ by sending $Z\in\env(X,Y)\backslash\{X\}$ to the endpoint of $\HSL(X,Z)\cap \env(X,Y)$ other than $X$.  Similarly, define $\lpi: \env(X,Y)\backslash\{Y\}\to \lomega(X,Y)$ by sending $Z\in\env(X,Y)$ to the endpoint of $\HSL (Z,Y)\cap \env(X,Y)$ other than $Y$.
 	 	\begin{lemma}\label{lem:pi:continuous}
 		Both $\rpi$ and $\lpi$ are continuous. 	\end{lemma}
   \begin{proof}
   The proofs for $\pi_l$ and $\pi_r$ are similar. In the following, we shall prove that $\rpi$ is continuous, and omit that of $\lpi$.

	  Let $Z_n, Z\in\env(X,Y)\backslash\{X\}$ with $Z_n\to Z$.   Notice that  $\env(X,Y)$ is compact.  By Theorem \ref{thm:maximally:stretched:lamination} and Lemma   \ref{lem:extendability}(i),   any limit of  $\{\rpi(Z_n)\}_{n\geq1}$ is a right boundary point.  On the other hand, by  Theorem \ref{prop:PanWolf2022}, we see that any such limit is contained in  $\HSL(X,Z)$. Therefore,   as $n\to\infty$, $\rpi(Z_n)$ converges to the (unique) right boundary point on $\HSL(X,Z)$, which is exactly $\rpi(Z)$.   \end{proof}

 		\begin{lemma}\label{lem:bd2}
Let  $X,Y\in\T(S)$ with $X\neq Y$. 
 \begin{enumerate}
  \item  $\lomega(X,Y)$, $\romega(X,Y)$ and $\bnd(X,Y)$ are compact.
 	\item $\bnd(X,Y)$ is empty if and only if  $\env(X,Y)$ is a geodesic.
 \end{enumerate}
		\end{lemma}
		\begin{proof} 
		 (i)  Recall that $\env(X,Y)$ is compact. To prove the statement, it suffices to show that $\lbd(X,Y)$, $\rbd(X,Y)$ and $\bnd(X,Y)$ are closed subsets of $\env(X,Y)$.  Since the proofs are similar, we shall prove the statement about $\rbd(X,Y)$, and omit the proof of the other two cases. 
		 
		  Let $Z_n\in\rbd(X,Y)$ be a sequence that converges to some $Z\in\env(X,Y)$. Then by Lemma \ref{lem:pi:continuous}, we see that $\pi_r(Z_n)$ converges to  $\pi_r(Z)$. Then since $\pi_r(Z_n)=Z_n$, we see that the convergence statements imply that $\pi_r(Z) = Z$, i.e. that  $Z\in\rbd(X,Y)$.

			(ii) 	It is clear that  $\bnd(X,Y)$ is empty if $\env(X,Y)$ is a  geodesic.  We now assume that $\env(X,Y)$ is not a geodesic. Note that  $\rbd(X,Y)$ contains $Y$ (Lemma \ref{lem:concatenation}). Moreover,  $\rbd(X,Y)$ strictly contains $Y$, since for any $Z \in \env(X,Y) \setminus [X,Y]$, the terminal point of the oriented segment $\HSL(X,Z)\cap\env(X,Y)$ is a right boundary point distinct from $Y$. Notice that $\rbd(X,Y)$ is compact and does not contain $X$. Let $Z\in \rbd(X,Y)$ be a point which realizes the largest distance to $Y$, i.e
   \begin{equation*}
       \dth(Z,Y)=\max \{\dth(Z',Y):Z'\in \rbd(X,Y)\}.
   \end{equation*}
     We shall prove that  $Z\in \lbd(X,Y)$. 
       By Lemma \ref{lem:concatenation},  $\HSL(Z,Y)\cap \env(X,Y)$ is an oriented segment with terminal point $Y$.  For any $Z'\in\HSL(Z,Y)\cap \env(X,Y)$, we have $\Lambda(Z',Y)=\Lambda(Z,Y)$ which strictly contains $\Lambda(X,Y)$, so  $Z'\in\rbd(X,Y)$ by Lemma \ref{lem:extendability}(i). The assumption that $Z$ realizes the maximal distance to $Y$ in $\rbd(X,Y)$ implies that $Z$ is the initial point of $\HSL(Z,Y)\cap\env(X,Y)$, hence $Z\in\lbd(X,Y)$ by definition.    In particular, $\bnd(X,Y)=\rbd(X,Y)\cap \lbd(X,Y)$ contains $Z$, and  is thus nonempty. 
		\end{proof}

  \subsection{Shapes of envelopes}
 	\label{sec:env:topology}
 	For $0\leq t\leq 1$, define 
 	$$M_t(X,Y)=\{Z\in\env(X,Y): \dth(X,Z)=t\dth(X,Y)\}.$$
 	There is a natural topology on $M_t(X,Y)$ induced from $\T(S)$.

The next lemma will be the final ingredient in the proofs of Theorem~\ref{thm:env:shape} and Corollary~\ref{cor:bnd:description} on the shapes of the envelopes. The lemma will also analyze these level sets $M_t$.
 	\begin{lemma}\label{lem:env:top}
 		Let $X,Y\in\T(S)$ with $X\neq Y$.
 		\begin{enumerate}
 		 \item Both $\lbd(X,Y)$ and $\rbd(X,Y)$ are cones over $\bnd(X,Y)$.
 		   \item $\env(X,Y)$ is a cone over $\rbd(X,Y)$. 
 		   \item   $\env(X,Y)$ is a cone over $\lbd(X,Y)$. 
 		   	\item For any $0\leq t\leq 1$, the space $M_t(X,Y)$ is contractible and compact, and $\env(X,Y)$ retracts onto $M_t(X,Y)$.
       \item $\env(X,Y)$ retracts onto the segment $[X,Y]$.
 			\end{enumerate}
 	\end{lemma}
 	\begin{proof}
 	(i) We present the proof for $\rbd(X,Y)$. The proof for $\lbd(X,Y)$ follows similarly.  
	Let $Z$ be an arbitrary point of $\rbd(X,Y)$.  By the first item in Lemma \ref{lem:extendability},  $\Lambda(Z,Y)$ strictly contains $\Lambda(X,Y)$. Notice that  for any $Z'\in\HSL(Z,Y)\cap\env(X,Y)$, we have $\Lambda(Z',Y)=\Lambda(Z,Y)$, which strictly contains $\Lambda(X,Y)$.  Again by the first item in Lemma \ref{lem:extendability}, this implies that $Z'\in \rbd(X,Y)$. It then follows from the arbitrariness of $Z'$ that $\HSL(Z,Y)\cap \env(X,Y)\subset\rbd(X,Y)$.  Let  $Z''\in \HSL(Z,Y)\cap \env(X,Y)$ be a left boundary point of $\env(X,Y)$.  Then, $Z''\in\bnd(X,Y)$ and   $[Z'',Y]=\HSL(Z,Y)\cap \env(X,Y)\subset \rbd(X,Y)$. We conclude that all elements $Z\in \rbd(X,Y)$ lie on segments between $\bnd(X,Y)$ and $Y$.

On the other hand, segments $[Z'',Y]$ with $Z'' \in \bnd(X,Y)$ are clearly contained in $\rbd(X,Y)$, as the maximally stretched lamination $\Lambda(Z,Y)$ with $Z\in [Z'',Y]$ strictly contains $\Lambda(X,Y)$ by Lemma~\ref{lem:extendability}. 

Consider the quotient $\bnd(X,Y)\times[0,1]/_\sim$ , where the equivalence relation is defined by $(U,1)\sim (V,1)$ for $U,V\in \rbd(X,Y)$. The discusion above defines a surjective map $\Phi: \bnd(X,Y)\times[0,1]/_\sim\to\rbd(X,Y)$ by sending  $(Z,t)\in\bnd(X,Y)\times [0,1]$ to the point on $[Z,Y]\hs$ whose distance to  $Y$ is $t\dth(Z,Y)$, i.e.  
 		\begin{equation}\label{eq:BdToRbdHomeo}
		{\dth(Z,\Phi(Z,t))}=t\cdot {\dth(Z,Y)}.
 		\end{equation}
  By Theorem \ref{thm:PanWolf2022:uniqueness}(i), this map is injective. Next, we will prove the map $\Phi$ and its inverse map are both continuous. Indeed, let $(Z_n,t_n), (Z,t)\in \bnd(X,Y)\times [0,1]$ be such that $(Z_n,t_n)\to (Z,t)$ as $n\to\infty$. By Proposition \ref{prop:PanWolf2022}, we see that $[Z_n,Y]\to [Z,Y]$ as $n\to\infty$. Since $\Phi(Z_n,t_n)\in[Z_n,Y],~ \Phi(Z,t)\in[Z,Y]$ and $\dth(Z_n,\Phi(Z_n,t_n))=t_n\dth(Z_n,Y)\to t \dth(Z,Y)$, we see that $\Phi(Z_n,t_n)\to\Phi(Z,t)$. This shows that $\Phi$ is continuous. For the inverse direction, let $Z_n, Z\in \rbd(Z,Y)$ such that $Z_n\to Z$ as $n\to\infty$. The discussion 
 in the first two paragraphs of this proof implies that $\HSL(Z_n,Y)\cap \env(X,Y)\subset \rbd(X,Y)$ and $\HSL(Z,Y)\cap \env(X,Y)\subset \rbd(X,Y)$. In particular, their initial points $\pi_l(Z_n)$ and $\pi_l(Z)$, which are contained in $\lbd(X,Y)$ by definition, are also contained in $\rbd(X,Y)$. Therefore, we have $\pi_l(Z_n),\pi_l(Z)\in\bnd(X,Y)$.  By Lemma \ref{lem:pi:continuous}, we see that $\pi_l(Z_n)\to\pi_l(Z)$. Combining this with the facts $\dth(\pi_l(Z_n),Z_n)\to \dth(\pi_l(Z),Z)$ and $\dth(\pi_l(Z_n),Y)\to\dth(\pi_l(Z),Y)$, we see that the inverse of $\Phi$ is also continuous.   Therefore, the map $\Phi$ is a homeomorphism, and affinely sends $\{Z\}\times [0,1]$ to a segment in $\rbd(X,Y)$ for any $Z\in\bnd(X,Y)$.

    (ii) Consider the quotient $\rbd(X,Y)\times [0,1]/_\sim$ where the equivalence relation is defined by $(Z,0)\sim (W,0)$ for any $Z,W\in\rbd(X,Y)$.
    Define a map 
 		$\Pi:\rbd(X,Y)\times[0,1]/_\sim\to \env(X,Y) $ by sending $(Z,t)\in\rbd(X,Y)\times [0,1]$ to the point on $[X,Z]\hs$ whose distance from $X$ is $t\dth(X,Z)$, i.e.  
 		\begin{equation*}
 			{\dth(X,\Pi(Z,t))}=t\cdot {\dth(X,Z)}.
 		\end{equation*}
 		It is clear that $\Pi$ is a bijection: injectivity follows from Theorem~\ref{thm:PanWolf2022:uniqueness} while surjectivity is a consequence of Lemma~\ref{lem:extendability} (as any segment $[X,Z]$ with $Z \in \env(X,Y)$ extends to a segment that meets $\rbd(X,Y)$).
 		That $\Pi$ and its inverse are continuous follows from Proposition \ref{prop:PanWolf2022}. 
 		
 		(iii) It follows similarly as (ii).

 		(iv) That $M_t(X,Y)$ is compact follows from the definition and the compactness of $\env(X,Y)$.  It remains to show the contractibility.

   For any fixed $t\in(0,1)$, define   
 		$$\phi_t:\env(X,Y)\times[0,1]\to \env(X,Y)$$
 		in the following way. 
 		\begin{itemize}
 			\item  If $\dth(X,Z)\geq t\cdot \dth(X,Y)$, define $\phi_t(Z,r)$ to be the point on $[Z_t,Z]\hs$ such that  
 		\begin{equation*}
 			{\dth(Z_t,\phi_t(Z,r))} =r\cdot {\dth(Z_t,Z)}
 		\end{equation*}
 		where $Z_t:=[X,Z]\hs\cap M_t(X,Y)$.
 		\item If $\dth(X,Z)< t\cdot \dth(X,Y)$, define $\phi_t(Z,r)$ to be the point on $[Z,Z_t]\hs$ such that  
 		\begin{equation*}
 			{\dth(\phi_t(Z,r), Z_t)} =r\cdot {\dth(Z,Z_t)}
 		\end{equation*}
 		where $Z_t:=[Z,Y]\hs\cap M_t(X,Y)$.
 		\end{itemize}
 		It is clear that $\phi_t$ is a deformation retract as $r\to 0$ from $\env(X,Y)$ to $M_t(X,Y)$. 
 		
   Finally, we assert that $M_t(X,Y)$ is contractible. Let $Z_0 = [X,Y] \cap M_t(X,Y)$ be the point where the segment $[X,Y]$ meets $M_t(X,Y)$: we will contract $M_t(X,Y)$ to $Z_0$. To that end, choose $Z \in M_t(X,Y)$; of course, the segment $[X,Z]$ may be extended to meet $\rbd(X,Y)$, say at $W \in \rbd(X,Y)$.  As in the earlier parts of this argument, we observe that $[W,Y] \subset \rbd(X,Y)$, and moreover, for $W_s \in [W,Y]$, we have that $d_{Th}(X,W_s)= d_{Th}(X,W) + d_{Th}(W,W_s)$, and so in particular $d_{Th}(X,W_s) >t$.  Thus there is some point $Z_s \in M_t(X,Y)$ on the segment $[X,W_s]$ between $X$ and $W_s$.  The family $Z_s$ defines a (continuous, cf. Proposition~\ref{prop:PanWolf2022}) retraction on $M_t(X,Y)$ that takes $Z_s$ to $Z_0$.

   (v) The retraction described in the previous paragraph that retracts $M_t(X,Y)$ to the singleton $M_t(X,Y)\cap[X,Y]$ for every $0\leq t\leq 1$  also retracts $\env(X,Y)$ to $[X,Y]$.
  \end{proof}

 An immediate corollary of Lemma~\ref{lem:env:top} is the following.

\begin{corollary}\label{cor:env:contractible}
    For any $X,Y \in \T(S)$, the envelope $\env(X,Y)$ is contractible.
\end{corollary}

Note that, as our methods do not address the question whether $\bnd(X,Y)$ is a manifold, we are not asserting that $\env(X,Y)$ is a topological ball.

		\subsection{Proof of Theorem \ref{thm:env:shape} and Corollary \ref{cor:bnd:description}}\label{sec:proof:continuity}
		\begin{proof}[Proof of Theorem \ref{thm:env:shape}]
			Item (i)  in Theorem \ref{thm:env:shape} follows from Corollary \ref{cor:unique:geodesic},   item (ii) follows from Proposition \ref{prop:out-in-envelopes}, and both items (iii) and  (iv) follow from Lemma \ref{lem:env:top}. Item (v) follows from Corollary \ref{cor:env:contractible}.
		\end{proof}

    \begin{proof}
        [Proof of Corollary \ref{cor:bnd:description}]
        By Lemma \ref{lem:extendability}, $\bnd(X,Y)$ is the set of points $Z$ such that both $\Lambda(X,Z)$ and $\Lambda(Z,Y)$ strictly contain $\Lambda(X,Y)$. Moreover, such a point $Z$ is the intersection point between the harmonic stretch line $\HSL(X,Z)$ and $\bnd(X,Y)$. 
         Hence, there is a bijection between  $\bnd(X,Y)$ and the set $\mathcal{H}$ of those harmonic stretch lines, $$\mathcal{H}:=\{\HSL(X,Z): \Lambda(X,Z) \text{ strictly contains } \Lambda(X,Y)\}.$$ By Proposition \ref{prop:PanWolf2022}, this bijection is a homeomorphism. By Theorem \ref{thm:out:in:env}, we see that there is a bijection between $\mathcal{H}$ and the union $\cup_{\lambda'}\cone(\lambda')$, where $\lambda'$ ranges over all chain recurrent geodesic laminations strictly containing $\Lambda(X,Y)$.   By Proposition \ref{prop:hsl:extended}, this bijection is a homeomorphism.
         This completes the proof.
    \end{proof}

  \section{Continuity of envelopes}\label{sec:continuity}
  The goal of this section is to show the following continuity.
  \begin{theorem}\label{thm:continuity}
  Let $X,Y,X_n,Y_n\in\T(S)$ be such that $X\neq Y$, $X_n\to X$ and $Y_n\to Y$. Then $\env(X_n,Y_n)\to \env(X,Y)$ in the  Hausdorff topology, as $n\to\infty$.
  	  \end{theorem}

  	We will prove  Theorem 
 \ref{thm:continuity} in two steps. In the first step,   we  prove that $\env(X_n,Y_n)$ converges to $\env(X,Y)$ pointwise (Lemma \ref{lem:approximate}). In the second step, we improve the pointwise convergence to uniform convergence, i.e. convergence in Hausdorff topology (See Section \ref{subsec:hausdorff} for definitions).

       Notations:
       \begin{itemize}
       	\item For $M\subset \T(S)$, the notation $\mathscr N_\epsilon(M)$ represents the $\epsilon$ neigbourhood of $M$ under the symmetrization $\dth^S$ of the Thurston metric (see Section \ref{subsec:hausdorff}). To simplify notations, we set $\mathscr{N}_\epsilon(X):=\mathscr{N}_\epsilon(\{X\})$. For a geodesic lamination $\lambda$, the notation $N_\epsilon(\lambda)$  represents the $\epsilon$ neighbourhood of $\lambda$ on $X$ (see the clarification in Section \ref{subsec:nghd}). 
        \item  The \emph{closed out-going ball}  centered at $X$ of radius $r$ is defined as:
        \begin{equation*}
          \mathscr{K}=  \mathscr{K}(X,r)=  \{Y\in\T(S):\dth(X,Y)\leq  r\}.
        \end{equation*} 
       \end{itemize}

 \subsection{Pointwise convergence}

We begin our proof of Theorem~\ref{thm:continuity} in this subsection by showing (Lemma~\ref{lem:approximate}) that the envelope $\env(X_n, Y_n)$, defined by approximates $X_n$ and $Y_n$ of the defining initial and terminal points $X$ and $Y$ of $\env(X,Y)$, lies near all of $\env(X,Y)$.  The proof of this is through using an algorithm to find -- for each point $Z \in \env(X,Y)$ -- an approximate $Z_n \in \env(X_n, Y_n)$. We will give an informal description of that algorithm at the start of the proof of Lemma~\ref{lem:approximate}; the estimates in that lemma rely on the initial technical Lemma~\ref{lem:sequence:epsilon}.
  \begin{lemma}\label{lem:sequence:epsilon}
       Let $X, Y, X_n,Y_n\in\T(S)$ be such that  $X\neq Y$, $X_n\to X$ and $Y_n\to Y$. For $\lambda= \Lambda(X,Y)$, let $\mathscr K$ be the closed out-going ball centered at $X$ with radius $2\dth(X,Y)$.    Then for any $Z\in \outenv(X,\lambda)\cap \inenv(Y,\lambda)$ and any $\epsilon>0$,    there exists a   sequence  of positive numbers 
  	     $$\epsilon>\epsilon_0>\epsilon_1>\epsilon_2>\cdots $$
       such that
        \begin{enumerate}
           \item for any $X',X'',Y',Y''\in \mathscr K$ with $X'\in \mathscr N_{\epsilon_{n+1}}(X'')$ and   $Y'\in \mathscr 
           N_{\epsilon_{n+1}}(Y'')$, we have $d_H([X',Y'],[X'',Y''])<\epsilon_n $, and
        	\item for any $Z'\in \mathscr N_{\epsilon_{n+1}}([X,Z])$ and any $Y'\in \mathscr N_{\epsilon_{n+1}}(Y)$, we have $\Lambda(Z',Y')\subset N_{\epsilon_{n}}(\lambda)$, and
        	\item for any $X'\in[X,Z]$,  any $Z'\in \mathscr N_{\epsilon_{n+1}}(X')$ and any $\lambda'\subset N_{\epsilon_{n+1}}(\lambda)$, there exists a harmonic stretch ray $\mathscr{R}(Z')\subset \overline{\outenv(Z',\lambda')}$ starting at $Z'$ such that 
  	   \begin{equation}\label{eq:approximate:Z'Z}
  	   	d_H(\mathscr{R}(Z')\cap \mathscr K,\HSR(X',Z)\cap \mathscr K)<\epsilon_n.
  	   \end{equation}
        \end{enumerate}
         where  $\HSR(X', Z)$ is the harmonic stretch ray that begins at $X'$ and passes through $Z$.
  \end{lemma}

We remind the reader of the discussion around Proposition~\ref{prop:PanWolf2022} which explains that Hausdorff convergence of segments $[X_m,Y_m]$ implies Hausdorff convergence of subsegments $[W_m,Z_m] \subset [X_m,Y_m]$. In particular, the paths in Lemma \ref{lem:sequence:epsilon}(i) (resp. Lemma \ref{lem:sequence:epsilon}(iii))  not only have small Hausdorff distance, but also travel in the same direction.

\begin{proof}
     We start with the construction of $\epsilon_0$. 
  	   By Proposition \ref{prop:PanWolf2022}, and the compactness of $\mathscr K$,  there exists $\delta_1=\delta_1(\mathscr K,\epsilon)<\epsilon$, such that for any $X'$, $X''$, $Y'$, $Y''\in \mathscr K$ with $X'\in \mathscr{N}_{\delta_1}(X'')$ and   $Y'\in \mathscr{N}_{\delta_1}(Y'')$, we have $$d_H([X',Y'],[X'',Y''])<\epsilon. $$
  	To see this, suppose to the contrary that there exists $\epsilon>0$ such that for any $\delta_1>0$ there exist $X'$, $X''$, $Y'$, $Y''\in \mathscr K$ with $X'\in \mathscr{N}_{\delta_1}(X'')$,   $Y'\in \mathscr{N}_{\delta_1}(Y'')$, and $d_H([X',Y'],[X'',Y''])\geq\epsilon$. In particular, taking $\delta_n=1/n$ yields a sequence $X'_n$, $X''_n$, $Y'_n$, $Y''_n\in \mathscr K$ with $X'_n\in \mathscr{N}_{1/n}(X''_n)$,   $Y'_n\in \mathscr{N}_{1/n}(Y''_n)$, and $d_H([X'_n,Y'_n],[X''_n,Y''_n])\geq\epsilon$. Up to taking subsequences if necessary, we may assume that $X''_n\to X''_\infty\in\mathscr{K}$ and  $Y''_n\to Y''_\infty\in\mathscr{K}$. Then assumption that $X'_n\in \mathscr{N}_{1/n}(X''_n)$ and $Y'_n\in \mathscr{N}_{1/n}(Y''_n)$ implies that $X'_n\to X''_\infty$ and $Y'_n\to Y''_\infty$ as well. By Proposition \ref{prop:PanWolf2022}, as $n\to\infty$, we have $[X'_n,Y'_n]\to [X''_\infty,Y''_\infty]$ and $[X''_n,Y''_n]\to [X''_\infty,Y''_\infty]$ uniformly. In particular, $d_H([X'_n,Y_n'],[X''_n,Y_n''])\to 0$ as $n\to \infty$, contradicting to the aforementioned inequality $d_H([X'_n,Y_n'],[X''_n,Y_n''])\geq\epsilon$.

  	    From Theorem \ref{thm:maximally:stretched:lamination} and the compactness of $[X,Z]$, it follows that  there exists $\delta_2<\epsilon$ such that for any $Z'\in \mathscr{N}_{\delta_2}([X,Z])$ and any $Y'\in \mathscr{N}_{\delta_2}(Y)$, we have $\Lambda(Z',Y')\subset N_{\epsilon}(\lambda)$.  
  	    
  	  Recall that from the assumption $Z\in\outenv(X,\lambda)\cap\inenv(Y,\lambda)$, we see that $\Lambda(X,Z)=\Lambda(X,Y)=\lambda$; moreover, for $X' \in [X,Z]$, we then also have $\Lambda(X',Z)=\Lambda(X,Z)=\lambda$.   Then from Corollary \ref{cor:out:in:approximation} and the compactness of $[X,Z]$, it follows that there exists $\delta_3<\epsilon$ such that for any $X'\in[X,Z]$, any  $Z'\in \mathscr{N}_{\delta_3}(X')$ and any $\lambda'\subset N_{\delta_3}(\lambda)$, there exists a harmonic stretch ray $\mathscr{R}(Z')\subset \overline{\outenv(Z',\lambda')}$ starting at $Z'$ such that 
  	   \begin{equation*}
  	   	d_H(\mathscr{R}(Z')\cap \mathscr K,\HSR(X',Z)\cap \mathscr K)<\epsilon.
  	   \end{equation*}
      Setting  $\epsilon_0=\min\{\delta_1,\delta_2,\delta_3\}$ completes the construction of $\epsilon_0$.  

      Inductively applying the preceding argument yields the desired sequence $\{\epsilon_n\}$, which completes the proof.
\end{proof}

  \begin{lemma}\label{lem:length:kappa}Every strictly increasing sequence of chain recurrent geodesic laminations $\eta_0\subsetneq\eta_1\subsetneq\eta_2\subsetneq\cdots \subsetneq\eta_m$ on $S$  has length at most  $\kappa:=2|\chi(S)|$ where $\chi(S)$ is the characteristic of $S$, i.e. $m\leq \kappa$. Moreover, if $m=\kappa$ then $\eta_m$ is maximally CR.
  \end{lemma}
  \begin{proof}
  Let $X\in\T(S)$. For  $i\geq 0$, each component $(X\backslash\lambda_i)\backslash\lambda_{i+1}$ has area at least $\pi$, so the sequence would terminate after at most $\kappa:=\frac{2\pi|\chi(S)|}{\pi}=2|\chi(S)|$ steps, where $\chi(S)$ is the characteristic of $S$.
  \end{proof}\begin{lemma}\label{lem:approximate}
  	      Let $X, Y, X_n,Y_n\in\T(S)$ be such that  $X\neq Y$, $X_n\to X$ and $Y_n\to Y$. 
  	    	Then  any point in $\env(X,Y)$ is a limit of some sequence $\{Z_n\}$ with  $Z_n\in\env(X_n,Y_n)$.
  	    \end{lemma}
  	    \begin{proof}
  	
  	    Let $\lambda=\Lambda(X,Y)$ and $\lambda_n=\Lambda(X_n,Y_n)$.
  	    By Proposition \ref{prop:out-in-envelopes},  it suffices to prove the lemma for any  point in $\outenv(X,\lambda)\cap \inenv(Y,\lambda).$ Let $Z$ be such a point distinct from $X$ and $Y$.
       
  	   Let $\mathscr K$ be the closed out-going ball centered at $X$ with radius $2\dth(X,Z)$.   
  	   
  	   Let $\epsilon$ be any positive real number with $5\epsilon<\dth(X,Y)-\dth(X,Z)$. Let $\epsilon>\epsilon_0>\epsilon_1>\cdots$ be the sequence obtained from Lemma \ref{lem:sequence:epsilon}, and let $\kappa$ be the fixed constant from Lemma~\ref{lem:length:kappa}.  	   Since $X_n\to X$ and $Y_n\to Y$, there exists $D>0$ such that 
  	\begin{equation}\label{eq:app:XnYn}
  	    X_n\in \mathscr{N}_{\epsilon_{2\kappa+1}}(X)\subset \mathscr{N}_{\epsilon_{2\kappa+1}}([X,Z]),~ Y_n\in \mathscr{N}_{\epsilon_{2\kappa+1}}(Y)
  	\end{equation}
      and $$\lambda_n=\Lambda(X_n,Y_n)\subset N_{\epsilon_{2\kappa+1}}(\lambda)$$ for all $n>D$. 
  	The lemma then follows from  Lemma \ref{lem:approximate:Z} below and the arbitrariness of $\epsilon$: 
     \end{proof}
    
  \begin{lemma}\label{lem:approximate:Z} With notations as above, for any $n>D$, the surface $Z$ is contained in the $\epsilon$ neighbourhood of $\env(X_n,Y_n)$.
     \end{lemma}
     \begin{proof}
     We will describe an algorithm for finding a point $Z_n$$\in \env(X_n,Y_n)$ near $Z$ where distance is measured with respect to the symmetrization of the Thurston metric.  The precise description of this process runs for a few pages, so we first give a rough outline.  We start at $X_n$ and shoot off a ray that shadows the ray from $X$ to $Z$.  It is possible that approximate ray ends up close to $Z$ and we may just use a point $Z_n$ on that ray that is near $Z$.   
     
     But it is also possible that this approximate ray is a bit \enquote{short}, in that it hits a boundary point, say $W_n^1$ in $\rbd(X_n, Y_n)$ before it gets sufficiently close to $Z$. In that case, we find, near $W_n^1$, a shadow point $X_n^1$ on the segment $[X,Z]$, and we start over: i.e. we look for a second approximating ray, this time beginning at this boundary point $W_n^1 \in \rbd(X_n,Y_n)$ and shadowing the ray from the new point $X_n^1$ to $Z$.  
     
    Now, in this iteration, two things have happened:  first, since $X_n^1$ is along the segment $[X,Z]$, we have made progress towards $Z$.  Second, because we are starting our new ray from a right boundary point $W_n^1$, the complexity of the maximally stretched lamination $\Lambda(W_n^1, Y_n)$ has increased over that of the original maximally stretched lamination $\Lambda(X_n,Y_n)$.

    We continue this process, always making progress towards $Z$ and always increasing the complexity of the maximally stretched lamination.  But by Lemma~\ref{lem:length:kappa}, we eventually have to stop as the maximally stretched lamination can no longer be enlarged.  That will mean that somewhere in this process, one of our rays was close to $Z$.  
    
  Here, all of the uniform estimates were from Lemma~\ref{lem:approximate}.

   We now describe this algorithm more carefully.
  	 
  	   \textbf{Step 1:} By Lemma \ref{lem:sequence:epsilon}, there exists  a harmonic stretch ray $\mathscr{R}(X_n)\subset\overline{\outenv(X_n,\lambda_n)}$ starting at $X_n$ and shadowing $\HSR(X,Z)$ in the sense that
  	   \begin{equation}\label{eq:app:Xn}
  	   	d_H(\mathscr{R}(X_n)\cap \mathscr K,\HSR(X,Z)\cap \mathscr K)<\epsilon_{2\kappa}.
  	   \end{equation}
  	   Let $W^1_n$ be the right boundary point of $\mathscr{R}(X_n)$ in $\env(X_n,Y_n)$.  Then $[X_n,W^1_n]\subset\env(X_n,Y_n)$ and  also $\Lambda(W^1_n,Y_n)$ strictly contains $\lambda_n$ (by Lemma \ref{lem:extendability}(i)).  If $Z$ is contained in the $\epsilon_{2\kappa}$ neighbourhood of $[X_n,W^1_n]$, then the algorithm terminates.  (In this case,  claim  1 follows since $[X_n,W^1_n]=\mathscr{R}(X_n)\cap\env(X_n,Y_n)\subset\env(X_n,Y_n)$ and $\epsilon_{2\kappa}<\epsilon$.) Otherwise, we continue the algorithm.
  	   \bigskip
  
  \textbf{Step 2:}   By \eqref{eq:app:Xn}, there exists   $X^1_n\in[X,Z]$ such that $W^1_n\in \mathscr{N}_{\epsilon_{2\kappa}}(X^1_n)$.  By the second item of  Lemma \ref{lem:sequence:epsilon}, this implies that $\Lambda(W^1_n,Y_n)\subset N_{\epsilon_{2\kappa-1}}(\lambda)$.   Applying  the third item of Lemma \ref{lem:sequence:epsilon} to $X^1_n$,  $W^1_n$ and $\Lambda(W^1_n,Y_n)$, we see that there exists  a harmonic stretch ray $\mathscr{R}(W^1_n)\subset\overline{\outenv(W^1_n,\Lambda(W^1_n,Y_n))}$ starting at $W^1_n$ such that
  	   $$ 	d_H(\mathscr{R}(W^1_n)\cap \mathscr K,\HSR(X^1_n,Z)\cap \mathscr K)<\epsilon_{2\kappa-2}. $$
  Let $W^2_n$ be the right boundary point of $\mathscr{R}(W^1_n)$ in the envelope $\env(W^1_n,Y_n)$ (here note the (new) left endpoint $W^1_n$ in the envelope under consideration).  Then, since $W^1_n\in  \env(X_n,Y_n)$,  we have $$[W^1_n,W^2_n]\subset\env(W^1_n,Y_n)\subset \env(X_n,Y_n)$$ and  that $\Lambda(W^2_n,Y_n)$ strictly contains $\Lambda(W^1_n,Y_n)$.  If $Z$ is contained in the $\epsilon_{2\kappa-2}$ neighbourhood of $[W^1_n,W^2_n]$, then the algorithm terminates.  (In this case,  claim 1 follows since $[W^1_n,W^2_n]=\mathscr{R}(W^2_n)\cap\env(W^1_n,Y_n)\subset\env(X_n,Y_n)$ and $\epsilon_{2\kappa-2}<\epsilon$) Otherwise, we continue the algorithm.
  
  	   \bigskip
  \textbf{Termination of the algorithm.} We claim:
  
  \begin{quotation}
      \textbf{Claim:} the algorithm terminates in at most $\kappa+1$ steps.
  \end{quotation}
  The claim  follows if the algorithm terminates at the $m$-th step with $m\leq \kappa$. Now suppose that the algorithm does not terminate in $\kappa$ steps. In this case, we have
   \begin{itemize}
   \item $\{W^i_n\}_{i\leq \kappa}\subset \env(X_n,Y_n)$ and harmonic stretch rays $\mathscr{R}(W^i_n)$ with  $W^0_n=X_n$  and $W^{i+1}_n$ being the right boundary point of $\mathscr{R}(W^i_n)$ in $\env(W^i_n,Y_n)$, 
   \item  $\{X^i_n\}_{i\leq \kappa}\subset [X,Z]$ with $X^0_n=X$ and $X^i_n\in [X^{i-1}_n,Z]\cap \mathscr{N}_{\epsilon_{2\kappa-2i+2}}(W^{i}_n)$,
  \item a strictly increasing sequence of  geodesic laminations
   \begin{equation}
   \label{eq:increasing:Lambda}\Lambda(X_n,Y_n)=	\Lambda(W^0_n,Y_n)\subsetneq \Lambda(W^1_n,Y_n)\subsetneq\cdots
  	\subsetneq \Lambda(W^{\kappa}_n,Y_n) \end{equation}
  with $ \Lambda(W^{i}_n,Y_n)$ contained in the $\epsilon_{2\kappa-2i+1}$ neighbourhood of $\lambda$.  	   \end{itemize}  	    
    In particular, $X^\kappa_n\in \mathscr{N}_{\epsilon_2}(W^\kappa_n)\subset \mathscr{N}_{\epsilon_1}(W^\kappa_n)$ and  $ \Lambda(W^{\kappa}_n,Y_n)\subset N_{ \epsilon_{1}}(\lambda)$.  Applying Lemma \ref{lem:sequence:epsilon} to $X^\kappa_n$, $W^{\kappa}_n$ and $ \Lambda(W^{\kappa}_n,Y_n)$, we see that there exists a harmonic stretch ray $\mathscr{R}(W^{\kappa}_n)\subset\overline{\outenv(W^{\kappa}_n,\Lambda(W^{\kappa}_n,Y_n))}$ starting at $W^{\kappa}_n$ such that 
  	\begin{equation}\label{eq:app:Wkn}
  	d_H(\mathscr{R}(W^{\kappa}_n)\cap\mathscr K, \HSR(X^\kappa_n,Z)\cap \mathscr K)<\epsilon_0.
       \end{equation}
Consider the strictly increasing sequence of chain recurrent geodesic laminations \eqref{eq:increasing:Lambda}. By Lemma \ref{lem:length:kappa}, the geodesic lamination $\Lambda(W_n^\kappa,Y_n)$ is maximally CR. Thus, by Corollary \ref{cor:unique:geodesic}, the envelope $\env(W_n^\kappa,Y_n)$ is exactly the geodesic segment $[W^\kappa_n,Y_n]$ and  $\overline{\outenv(W^\kappa_n,\Lambda(W^\kappa_n,Y_n))}$ is the harmonic stretch ray $\HSR(W_n^\kappa,Y_n)$ which begins at $W_n^\kappa$ and passes through $Y_n$. This implies that $\mathscr{R}(W^\kappa_n)=\HSR(W_n^\kappa,Y_n)$. It then follows from \eqref{eq:app:Wkn} that there exists $Z_n\in \HSR(W_n^\kappa,Y_n)\cap \mathscr K$ such that 
\begin{equation}
\label{eq:app:ZZ'}
    \dth(Z_n,Z)<\epsilon_0, \text{ and } \dth(Z,Z_n)<\epsilon_0.
\end{equation}
Moreover, not only is $Z_n\in \HSR(W_n^\kappa,Y_n)$, but we assert that the point $Z_n$ is contained in $[W_n^\kappa,Y_n]$. I.e. In the ray $\mathscr{R}(W_n^{\kappa})=\HSR(W_n^{\kappa}, Y_n)$, we have $Z_n$ in between $W_n^{\kappa}$ and $Y_n$ and not \enquote{after} $Y_n$. Otherwise, since the union $[X_n,W^1_n]\cup [W_n^1,W_n^2]\cup \cdots \cup \HSR(W_n^\kappa,Y_n)$ is a geodesic, we have $\dth(X_n,Z_n)\geq \dth(X_n,Y_n)$, so \begin{eqnarray*}
    \dth(X,Z)&\geq& \dth(X_n,Z_n)-\dth(Z,Z_n)-\dth(X_n,X)\\
    &\geq& \dth(X_n,Y_n)-\dth(Z,Z_n)-\dth(X_n,X)\\
     &\geq& \dth(X,Y)-\dth(X,X_n)-\dth(Y_n,Y)-\dth(Z,Z_n)-\dth(X_n,X)\\
    &\geq & \dth(X,Y)-4\epsilon,
\end{eqnarray*}
where the last inequality follows from \eqref{eq:app:XnYn} and \eqref{eq:app:ZZ'}. 
This contradicts the choice of $\epsilon$. 
Hence, the point $Z_n$ is contained in $[W_n^\kappa,Y_n]$. This proves the claim. The lemma then follows from the claim. 
  \end{proof}

  	\subsection{Finishing the proof of Theorem \ref{thm:continuity}}
   
  	   \begin{proof}[Proof of Theorem \ref{thm:continuity}]

     Note that, since  $X_n\to X$ and $Y_n\to Y$, so, 
     $$c(X,Y):=\sup\{\dth(X,Y),\dth(X,X_n),\dth(X_n,Y_n):n\geq 1\}<\infty.$$ In particular, for any $n\geq1$ and any $Z_n\in\env(X_n,Y_n)$, we see that $\dth(X_n,Z_n)\leq \dth(X_n,Y_n)$ and
     \begin{eqnarray*}
         \dth(X,Z_n)&\leq& \dth(X,X_n)+\dth(X_n,Z_n)\\
         &\leq &\dth(X,X_n)+\dth(X_n,Y_n)\\
         &\leq& 2c(X,Y).
     \end{eqnarray*}
      In particular, the envelope $\env(X_n,Y_n)$ is contained in the closed out-going ball centered at $X$ of radius $2c(X,Y)$. 
       
  	   Suppose to the contrary that  $\env(X_n,Y_n)$ does not converge to $\env(X,Y)$ with respect to the Hausdorff topology. Then there exists $\epsilon>0$ such that, up to taking a subsequence if necessary,
  	     	\begin{enumerate}
  		\item $\env(X_n,Y_n)$ is not contained in the $\epsilon$-neighbourhood of $\env(X,Y)$ for all $n\geq1$ , or
  		\item $\env(X,Y)$ is not contained in the $\epsilon$-neighbourhood of $\env(X_n,Y_n)$ for all $n\geq1$.
  	\end{enumerate} 	   	
  	We now proceed case by case, the first case being a general metric space argument and the second using the constructions in this section.

  	Case 1: $ \env(X_n,Y_n)$ is not contained in the  $\epsilon$ neighbourhood of $\env(X,Y)$ for all $n\geq1$. It follows that there exists $Z_n\in \env(X_n,Y_n)$ such that $\dth(\env(X,Y), Z_n)\geq \epsilon$.  Let $Z_\infty$ be an arbitrary accumulation point of $\{Z_n\}$ in $\T(S)$ (which exists by our discussion at the beginning of the proof).  Then $\dth(\env(X,Y),Z_\infty)\geq \epsilon$.  In particular, $Z_\infty\notin \env(X,Y)$. On the other hand, since $X_n\to X$ and $Y_n\to Y$ by assumption, 
  	\begin{eqnarray*}
  		\dth(X,Z_\infty)+\dth(Z_\infty,Y)&=&\lim_{n\to\infty} \left(\dth(X_n,Z_n)+\dth(Z_n,Y_n)\right)\\
  		&=&\lim_{n\to\infty}\dth(X_n,Y_n)\\&=&\dth(X,Y)
  	\end{eqnarray*}  	
  	which implies that  $Z_\infty\in \env(X,Y)$.  This is a contradiction.
  	\vskip 5pt
  	Case 2: $\env(X,Y)$ is not contained in the  $\epsilon$ neighbourhood of $\env(X_n,Y_n)$ for all $n\geq1$.  There then exists $W_n\in \env(X,Y)$ such that 
  	$$\dth(\env(X_n,Y_n),W_n) \geq\epsilon.$$ 
  	Let $W$ be an accumulation point of $\{W_n\}$. Then, up to a subsequence, 
  	\begin{equation}\label{eq:mid:W}
  		\dth(\env(X_n,Y_n),W) \geq\epsilon/2
  	\end{equation}
  	 for every sufficiently large $n$, which contradicts Lemma \ref{lem:approximate}.
  	This completes the proof.
    	   \end{proof}

\section{Extension to the Thurston boundary}\label{sec:boundary}

The goal of this section is to extend the previous result Theorem \ref{thm:continuity} to the Thurston boundary. 
 Along the way, we find (Corollary~\ref{cor:con:compactification} and Proposition~\ref{prop:cont:compactification2}) that the \enquote{visual} compactification of \tec space from a pole is independent of the pole. (This is weaker than claiming a base-point independent visual compactification in the sense of \enquote{visual boundary} in metric geometry that one might find in the standard literature, e.g \cite{BridsonHaefliger1999}.) 

The first two subsections provide some setting and then establish some basic continuity results for harmonic stretch paths $[X, Y_n]$ and $[X,\eta)$, where $\eta$ is a Thurston boundary point and $Y_n \in \T(S)$ tends to the boundary. The next subsection then studies the maximally stretched laminations, say $\Lambda(X, \eta)$ and the geometry of limits of such laminations. 

In Section~\ref{sec:env:ThurstonCompact}, we give a detailed description of the accumulation set of the infinite envelope $\env(X,\eta)$, providing (Proposition~\ref{prop:msl:bnd2})  necessary and sufficient conditions for a Thurston boundary point to be an accumulation point. It turns out that this accumulation set is a star about $\eta$, and we give some examples of possible closures of these infinite envelopes.

In Section~\ref{sec:ext:bnd}, we extend our descriptions (Theorem~\ref{thm:env:shape})  of the shape of the envelopes to this setting of infinite envelopes.  The description becomes more involved, as the set $\lbd(X,\eta)$ can be non-compact, and the set $\rbd(X,\eta)$ might not meet all geodesics in the envelope. Theorem~\ref{thm:env:shape:tree} summarizes the situation.

Finally, in section~\ref{sec:continuityThurstonEnvelopes}, we prove a broad continuity theorem for envelopes $\env(X,\eta)$, for $X \in \T(s)$ and $\eta$ a Thurston boundary point.

\subsection{The Thurston compactification}

Let $\MF(S)$  be  the space of equivalence classes of measured foliations, where two measured foliations are said to be \emph{equivalent} if they have the same intersection number with every (homotopy class of) simple closed curves on $S$.  By using length functions,  Thurston embedded the Teichm\"uller space $\T(S)$  into the space $\R^\mathcal{S}$ of functionals over the set $\mathcal S$ of homotopy classes of simple closed curves on $S$. Moreover, the Teichm\"uller space can also be embedded into the projective space $P\R^\mathcal{S}:=\R^\mathcal{S}/\R_+$. The closure $\overline{\T(S)}$ of $\T(S)$ in  $P\R^\mathcal{S}$ is homeomorphic to the closed unit ball in $\R^{6g-6}$. The boundary of $\overline{\T(S)}$ is exactly the space $\PMF(S)$ of projective classes of measured foliations. For more details about the Thurston compactification, we refer to \cite{FLP}.

 Walsh \cite{Walsh2014} proved that the Thurston compactification is the same as the \emph{horofunction compactification} with respect to the Thurston metric. 

\subsection{Harmonic stretch lines and the Thurston compactification}
Recall that \cite[Proposition 7.11]{PanWolf2022} every harmonic stretch ray converges (in the forward direction) to a unique point in the Thurston compactification. Moreover, we have

\begin{theorem}[\cite{PanWolf2022}, Theorem 1.10]\label{thm:ray:uniqueness:tree}
    For any $X\in\T(S)$ and $\eta\in\PMF(S)$, there exists a unique harmonic stretch ray $[X,\eta)$ which starts at $X$ and converges to $\eta$ in the Thurston compactification. 
\end{theorem}

\begin{proposition}[\cite{PanWolf2022}, Proposition 13.10]\label{prop:cont:compactification1}
Let $X,X_n\in\T(S)$ and $\eta,\eta_n\in \PMF(S)$ be such that $X_n\to X$ and $\eta_n\to \eta$ as $n\to\infty$. Then 
$[X_n,\eta_n)$ converges to $[X,\eta)$ locally uniformly.
\end{proposition}

Using the compactness of the Thurston compactification, we have 
\begin{corollary}\label{cor:con:compactification}
    Let $X,X_n\in\T(S)$ and $\eta,\eta_n\in \PMF(S)$. Then the following are equivalent:
    \begin{enumerate}
        \item $X_n\to X$ and $\eta_n\to \eta$,
     \item $[X_n,\eta_n)$ converges to $[X,\eta)$ locally uniformly.   
    \end{enumerate}
\end{corollary}
\begin{proof}
    That (i) implies 
 (ii) is the content of Proposition \ref{prop:cont:compactification1}. It remains to consider the inverse direction. We assume that $[X_n,\eta_n)$ converges to $[X,\eta)$ locally uniformly. In particular, this yields $X_n\to X$ as $n\to\infty$. For any convergent subsequence  $\eta_{n_j}$ with limit $\eta'$, by Proposition \ref{prop:cont:compactification1}, we see that $[X_{n_j},\eta_{n_j})$ converges to $[X,\eta')$. The assumption that $[X_n,\eta_n)$ converges to $[X,\eta)$ then implies that $\eta'=\eta$. This proves that $\eta_n\to\eta$.
\end{proof}

The goal of this subsection is to prove the following companion to Corollary~\ref{cor:con:compactification}:
\begin{proposition}\label{prop:cont:compactification2}
    Let $X,X_n, Y_n\in\T(S)$ and $\eta\in\PMF(S)$. Then the following are equivalent:
    \begin{enumerate}
        \item $X_n\to X$ and $Y_n\to \eta$,
     \item $[X_n,Y_n]$ converges to $[X,\eta)$ locally uniformly.   
    \end{enumerate}
\end{proposition}

\begin{remark}
    In some sense, by letting $X_n=X$, Proposition \ref{prop:cont:compactification2} asserts that not only does the visual compactification of the Thurston metric via harmonic stretch rays from a pole coincide with the Thurston compactification (\cite[Theorem 13.1]{PanWolf2022}), but is independent of the choice of the base point in the sense that Proposition \ref{prop:cont:compactification2}(i) describes the convergence in the Thurston compactification and Proposition \ref{prop:cont:compactification2}(ii) describes the convergence in the visual compactification from the pole/base point $X$. (As noted at the outset of the section, there is a nuance here in that this compactification that relies on a base point is different from the metric space visual compactification. In particular, the visual compactification adjoins geodesics, up to two such geodesics being equivalent if they remain \enquote{asymptotic}, i.e. with a finite distance of one another; In a complete CAT(0) space, this results in a compactification by proper geodesics originating at a fixed pole, as given any pair of poles and any geodesic originating at the first pole, there is a unique geodesic asymptotic to that given one originating at the second pole.  We are aware of no such result that might hold for the Thurston metric.)
\end{remark}
The proof will be divided into several lemmas. Let us start with the following observation.
\begin{lemma}\label{lem:Yn:eta0}
  Let $Y_n \to \eta$ as above.  Let $\eta_n$ be the endpoint in $\PMF(S)$ of the harmonic stretch ray $\HSR(X_n,Y_n)$ proceeding from $X_n$ through $Y_n$.  Then $[X_n,Y_n]$ converges to $[X,\eta)$ if and only if $[X_n,\eta_n)$ converges to $[X,\eta)$.  
\end{lemma}
\begin{proof}
    Without loss of generality, we may assume that $X_n\to X$ as $n\to\infty$. By assumption, we see that $[X_n,Y_n]\subset [X_n,\eta_n)$.   Moreover, since $Y_n$ diverges in $\T(S)$, for any closed ball $\mathscr K\subset\T(S)$, which is of positive radius and which is centered at $X$, we have $[X_n,Y_n]\cap \mathscr K=[X_n,\eta_n)\cap \mathscr K$  for $n$ large enough. Hence, the harmonic stretch segments $[X_n,Y_n]$ converge to $[X,\eta)$ locally uniformly   if and only if $[X_n,\eta_n)$ converges to $[X,\eta)$.
\end{proof}

\begin{lemma}\label{lem:Yn:eta1}
    Let $X,X_n, Y_n\in\T(S)$ and $\eta\in\PMF(S)$. Suppose that $[X_n,Y_n]$ converges to $[X,\eta)$ locally uniformly. Then $X_n\to X$ and $Y_n\to \eta$, as $n\to\infty$.
\end{lemma}

\begin{proof}
   Let $\eta_n$ be the endpoint in $\PMF(S)$ of the harmonic stretch ray $\HSR(X_n,Y_n)$ proceeding from $X_n$ through $Y_n$. Combining the assumption in the statement, Lemma \ref{lem:Yn:eta0},  and Corollary \ref{cor:con:compactification}, we see that $X_n\to X$ and $\eta_n\to \eta$. In the following, we shall prove that $Y_n$ converges to $\eta$.

Recall that the harmonic stretch ray $[X_n,\eta_n)$  is defined by a chain recurrent geodesic lamination $\lambda_n$ and a surjective harmonic diffeomorphism $f_{n,\infty}:V_{n,\infty}\to X_n\backslash \lambda_n$ (see Section \ref{subsec:PHSL}). Let $\eta_{X_n}$ be the pushforward of the vertical foliation of the Hopf differential $\Hopf(f_{n,\infty})$.   The lemma is now a consequence of the following lemma.

\begin{lemma}\label{lem:Yn:length} 
Let  $X_n,X\in\T(S)$ and $\eta_n,\eta\in\PMF(S)$ with $X_n\to X$ and $\eta_n\to \eta$. Then, for any simple closed curve $\alpha$ on $S$, there exist  constants $C=C(X,\eta,\alpha)$ and $N=N(X,\eta,\alpha)$ depending only on $X$, $\eta$, and $\alpha$, such that for any $n>N$ and for any $Y_n\in[X_n,\eta_n)$, we have
\begin{equation}\label{eq:length:Yn}
    |\ell_\alpha(Y_n)-2{r_n}\cdot i(\eta_{X_n},\alpha)|\leq C,
\end{equation}
where $r_n=\exp(\dth(X_n,Y_n))$.
\end{lemma}

We remark on the constant $C$. Since both $\ell_\alpha(Y_n)$ and $i(\eta_{X_n},\alpha)$ are linear in $\alpha$, if $\alpha$ is scaled by a constant, then the constant $C$ is also scaled by the same constant.

\begin{proof}
    Note that the ray $[X_n,\eta_n)$ can also be constructed from a harmonic diffeomorphism $g_{n,\infty}:V_{n,\infty}\to Y_n\setminus \lambda_n$ with Hopf differential $\Hopf(g_{n,\infty})=(r_n)^2 \Hopf(f_{n,\infty})$ where $r_n=\exp(\dth(X_n,Y_n))$. Let $\eta_{Y_n}$ be the pushforward of the vertical foliation of $\Hopf(g_{n,\infty})$. In particular, we have $\eta_{Y_n}={r_n}\eta_{X_n}$. By \cite[Equation (3.2)]{PanWolf2022}, we see that 
    \begin{equation}\label{eq:length:Yn:1}
        \ell_\alpha(Y_n)\geq 2i(\eta_{Y_n},\alpha). 
    \end{equation}

    It remains to control $\ell_\alpha(Y_n)$ from above:  \begin{equation}\label{eq:length:Yn:2}
        \ell_\alpha(Y_n)\leq 2i(\eta_{Y_n},\alpha)+C. 
    \end{equation}
    Note that the harmonic stretch ray $[X,\eta)$  is defined by a chain recurrent geodesic lamination $\lambda$ and a surjective harmonic diffeomorphism $f_{\infty}:V_{n,\infty}\to X\backslash \lambda$. Let $\eta_{X}$ be the pushforward of the vertical foliation of the Hopf differential  $\Hopf(f_{\infty})$. By Proposition \ref{prop:cont:compactification1} and \cite[Proposition 8.1(i)]{PanWolf2022} and using a diagonal argument via harmonic map ray approximates, we see that the $f_{n,\infty}$ converges to $f_\infty$ in the sense of \cite[Definition 4.5]{PanWolf2022}. Therefore, the Hopf differentials $\Hopf(f_{n,\infty})$ converge to $\Hopf(f_\infty)$. 

    Now let us fix a sufficiently large positive constant $R$ (cf. \cite[Theorem 5.1 and Theorem 7.1]{Minsky1992}). Consider the Minsky region $\mathscr{P}_{n,R}$ of $\Hopf(g_{n,\infty})$ (\cite[Theorem 5.1]{Minsky1992}, see also the comments between Lemma 3.7 and Lemma 3.8 in \cite{PanWolf2022} for a summary of the construction of $\mathscr{P}_R$). Each component of $\mathscr P_{n,R}$ contains the $R$ neighbourhood of at least one zero of $\Hopf(g_{n,\infty})$, and  is a locally convex subset of $V_{n,\infty}$ with boundary either a horizontal closed geodesic or a polygonal curve comprising alternatively horizontal and vertical segments of $\Hopf(g_{n,\infty})$. Furthermore, both the diameter and the boundary length of each component of 
    $\mathscr{P}_{n,R}$, with respect to the singular flat metric induced from $\Hopf(g_{n,\infty})$, are at most $cR$ for some constant depending only on the topology of $V_{n,\infty}$, hence is determined by the topology of $S$. Since  $\Hopf(g_{n,\infty})=(r_n)^2 \Hopf(f_{n,\infty})$, we see that both the diameter and the boundary length of each component of $\mathscr{P}_{n,R}$, with respect to the singular flat metric induced from $\Hopf(f_{n,\infty})$, are at most $cR/r_n$.  Since $\Hopf(f_{n,\infty})$ converges to $\Hopf(f_\infty)$, we see that there exists $N=N(X,\eta)>0$ and $\mathfrak{r}=\mathfrak{r}(X,\eta)$ such that for $n>N$ and $r_n>\mathfrak{r}$, each component of $\mathscr{P}_{n,R}$ is a (simply connected) polygon. (Otherwise, suppose there were one component of $\mathscr P_{n,R}$ that is not simply connected, then it contains at least one saddle connection whose length,  measured with respect to the singular flat metric induced from $\Hopf(f_{n,\infty})$, is at least $\delta/2$ for $n$ large enough, where $\delta$ is the length of the shortest saddle connection on $\Hopf(f_{\infty})$. Hence, the diameter of this component, again measured with respect to the singular flat metric induced from $\Hopf(f_{n,\infty})$, is at least $\delta/2$ for $n$ large enough. This contradicts the fact above that both the diameter and the boundary length of each component of $\mathscr{P}_{n,R}$, with respect to the singular flat metric induced from $\Hopf(f_{n,\infty})$, are at most $cR/r_n$.)  Combining this with \cite[Theorem 5.1 (iv) and (vi) and Theorem 7.1]{Minsky1992} (see also \cite[Lemma 3.1(iii)]{Wolf1991a}), we infer that each component of the image $g_n(\mathscr{P}_{n,R})$ on $Y_n$ is nearly a hyperbolic $k$ polygon with each edge of length comparable to $R$, where $k\leq 4g-2$ with $g$ being the genus of $S$. In particular, the diameter of each component of the image $g_n(\mathscr{P}_{n,R})$ is at most $c'R$ for some constant $c'$ depending on the topology of $S$.  \color{black}

    Next, we consider a representative on $Y_n$ in the homotopy class $[\alpha]$ of the simple closed curve $\alpha$ that is quasi-transversal\footnote{Here, by \enquote{quasi-transversal} we mean the representative realizes the minimal intersection number with $\eta_{Y_n}$ but possibly contains some leaf segment of $\eta_{Y_n}$.} to $\eta_{Y_n}$. That representative will pass through singular points of the foliation a finite number of times, with a bound on those number of passages depending only on $[\alpha]$ and the topology of $S$.  Inside each component of the image $g_n(\mathscr{P}_{n,R})$  of the Minsky region $\mathscr{P}_{n,R}$, we replace each subsegment of the representative by a hyperbolic geodesic segment with endpoints on the boundary of $g_n(\mathscr{P}_{n,R})$ . Each such geodesic segment has length at most the diameter of the underlying component of $g_n(\mathscr{P}_{n,R})$, which is at most $c'R$ by the discussion in the preceding paragraphs.

  Away from the image $g_n(\mathscr{P}_{n,R})$ of the Minsky region $\mathscr{P}_{n,R}$, we may replace the representative by a representative that is alternately horizontal and vertical: this is straightforward as the remaining segments of the curve, being outside of neighborhoods of singularities, may be taken to live in a single euclidean chart, foliated by horizontal and vertical segments.

The relevant estimates now follow from \cite[Lemma 3.1(iii)]{Wolf1991a}, see also \cite[Lemma 2.3]{Wolf1995} for a convenient summary:  the image of a horizontal arc is nearly a geodesic with length approximately twice the vertical measure of the arc, and the length of the image of a vertical arc outside the neighborhood of the singularity is negligible, with both estimates having exponentially small error in the Hopf differential $\Hopf(g_{n,\infty})$-distance $R$ to the singularity. 
Moreover, the length of each component of the representative in $g_n(\mathscr{P}_{n,R})$ is at most the diameter of $g_n(\mathscr{P}_{n,R})$, which is at most $c'R$ by the discussion in the preceding paragraphs. 
The displayed estimate \eqref{eq:length:Yn:2} then follows from adding up the lengths of the images of the horizontal portions,  the vertical portions, and the portion in $g_n(\mathscr{P}_{n,R})$, as long as $n>N(X,\eta)$ and $r_n>\mathfrak{r}(X,\eta)$.  Finally, we may choose a constant $N(X,\eta,\alpha)$ larger than $N(X,\eta)$ so that for each $n>N(X,\eta,\alpha)$ we have $\ell_\alpha(X_n)\leq 2\ell_\alpha(X)$. This implies $\ell_\alpha(Y_n)\leq e^{r_n}\ell_\alpha(X_n) \leq 2e^{\mathfrak{r}}\ell_\alpha(X)$. In summary, the displayed equation \eqref{eq:length:Yn:2} holds for any $n>N(X,\eta,\alpha)$ and any $Y_n\in[X_n,\eta)$ with a new constant $C'=\max\{C,2e^\mathfrak{r}\ell_\alpha(X)\}$.

 The displayed equation \eqref{eq:length:Yn} now follows from \eqref{eq:length:Yn:1}, \eqref{eq:length:Yn:2} and the fact $\eta_{Y_n}=r_n\eta_{X_n}$.
\end{proof}

We now return to the proof of Lemma \ref{lem:Yn:eta1}. By Lemma \ref{lem:Yn:length}, we see that for any simple closed curve $\alpha$ on $Y_n$, we have \begin{equation*}
    \left|\frac{1}{r_n}\ell_{\alpha}(Y_n)-2i(\eta_{X_n},\alpha)  \right|\leq \frac{C}{r_n}.
\end{equation*}
Letting $n\to\infty$, since $\eta_{X_n}\to \eta_X$ and $r_n\to\infty$, we see that for any simple closed curve $\alpha$,
\begin{eqnarray*}
    &&\left|\frac{1}{r_n}\ell_{\alpha}(Y_n)-2i(\eta_{X},\alpha)  \right|\\
    &\leq& \left|\frac{1}{r_n}\ell_{\alpha}(Y_n)-2i(\eta_{X_n},\alpha)  \right|+|2i(\eta_{X_n},\alpha)-2i(\eta_X,\alpha)|\\
    &\leq & \frac{C}{r_n}+|2i(\eta_{X_n},\alpha)-2i(\eta_X,\alpha)|\\
     &\to& 0.
\end{eqnarray*}
Since $\eta_X$ is a representative of $\eta$, and the above estimate holds for all simple closed curves $\alpha$, this implies that $Y_n$ converges to $\eta$.
\end{proof}

\begin{lemma}\label{lem:Yn:eta2}
    Let $X,X_n, Y_n\in\T(S)$ and $\eta\in\PMF(S)$. Suppose that $X_n\to X$ and $Y_n\to \eta$, as $n\to\infty$. Then  $[X_n,Y_n]$ converges to $[X,\eta)$ locally uniformly. 
\end{lemma}
 \begin{proof}
     Let $\eta_n$ be the endpoint of the harmonic stretch ray $\HSR(X_n,Y_n)$. To prove the lemma, by Lemma \ref{lem:Yn:eta0} and Corollary \ref{cor:con:compactification}, it suffices to prove that $\eta_n\to \eta$.  

     Let $\eta_{n_j}$ be a convergent subsequence of $\eta_n$. Let $\eta'$ be the corresponding limit. 
     Then by Proposition \ref{prop:cont:compactification1}, we see that $[X_{n_j},\eta_{n_j})$  converges to $[X,\eta')$. It then follows from  Lemma \ref{lem:Yn:eta0} and Lemma \ref{lem:Yn:eta1} that $Y_{n_j}\to \eta'$. Combined with the assumption that $Y_n\to\eta$, this yields that $\eta'=\eta$. The arbitrariness of $\eta'$ then implies that $\eta_n\to\eta$, which completes the proof. 
 \end{proof}

\begin{proof}[Proof of Proposition \ref{prop:cont:compactification2}] This follows from Lemma \ref{lem:Yn:eta1} and Lemma \ref{lem:Yn:eta2}.
\end{proof}

\subsection{Dual trees and maximally stretched laminations}\label{subsec:msl:trees}   Let $\eta\in\MF(S)$ be a measured foliation. Consider the lift $\widetilde{\eta}$ of $\eta$ to the universal cover $\widetilde{S}$ of $S$. The leaf space $T_\eta$ of $\widetilde{\eta}$,  consisting of equivalence classes of points on the universal cover $\widetilde S$ of $S$ where two
points are equivalent if they are contained in a connected leaf of  $\widetilde \eta$ (including
leaves which branch at singularities of $\widetilde\eta$), is an $\R$-tree, whose metric is defined as \emph{twice} the transverse measure of $\eta$. The projection map along leaves defines a $\pi_1(S)$-equivariant map from $\widetilde{S}$ to $T_\eta$. 

Recall that, by Theorem \ref{thm:ray:uniqueness:tree}, for any $X\in\T(S)$ and any $\eta\in\MF(S)$,  there exists a unique harmonic stretch ray $[X,\eta)$ which starts at $X$ and converges to (the projective class of) $\eta$ in the Thurston compactification.  

\begin{definition}
    Let $\Lambda(X,\eta)$ be the maximally stretched lamination along the ray $[X,\eta)$. 
   
\end{definition}
Here we define $\Lambda(X,\eta)$ in terms of the full ray $[X,\eta)$, which is not quite only the data of $X$ and $\eta$. However, we will soon show in Lemma~\ref{lem:optimal:lip:tree:} that not only is $\Lambda(X, \eta)$ the maximal stretch lamination for the surfaces along the ray $[X,\eta)$, but it is also the intersection of (the quotient of) the maximally stretched laminations for all optimal equivariant maps from the lift $\tilde{X}$ to the tree $T_\eta$; this will later allow us to rephrase the definition of $\Lambda(X,\eta)$ in more suggestive terms that only refer to the endpoints $X$ and $\eta$ (Definition \ref{def:msl:tree}).

Define 
\begin{equation}
\label{eq:L(X,eta)}
L(X,\eta):=\sup_{\alpha}\frac{2i(\eta,\alpha)}{\ell_\alpha(X)},
\end{equation}
where $\alpha$ ranges over all (homotopy classes of) simple closed curves. 
\begin{proposition}[\cite{PanWolf2022}, Proposition 13.2]\label{prop:optimallip:trees}
 Let $X\in\T(S)$ and $\eta\in\MF(S)$.  Then there exists an equivariant  Lipschitz projection map $f:\widetilde{X}\to T_\eta$ whose restriction to the lift of $\Lambda(X,\eta)$ is an affine map of factor $L(X,\eta)$ and which has (pointwise) Lipschitz constant strictly less than $L(X,\eta)$ outside the lift of $\Lambda(X,\eta)$. 
 \end{proposition} 
  
 \begin{remark}\label{rmk:lip:def}
 Here the pointwise Lipschitz constant of $f:\widetilde{X}\to T_\eta$ at $p\in \widetilde{X}$ is defined as:
 \begin{equation*}
     L_f(p):=\limsup_{q\to p,~q\neq p}\frac{d_{T_\eta}(f(q),f(p))}{d_{\widetilde{X}}(q,p)}.
 \end{equation*}
 \end{remark}
 \begin{remark}\label{rmk:lip:constuction}
 \begin{enumerate}
     \item The Lipschitz map in the statement of Proposition \ref{prop:optimallip:trees} is obtained as follows.  Recall that the harmonic stretch ray $[X,\eta)$ is defined by a harmonic diffeomorphism $f_\infty:X_\infty\to X\setminus \Lambda(X,\eta)$. The pushforward of the vertical measured foliation of $f_\infty$ extends to a unique measured foliation $\eta_X$ on $X$ (\cite[Lemma 7.7]{PanWolf2022}), which is measure equivalent to a scaling of $\eta$ (\cite[Proposition 7.11]{PanWolf2022} and Theorem \ref{thm:ray:uniqueness:tree}). The Lipschitz map $f:\widetilde{X}\to T_\eta$ is then the projection map along leaves of the lift of the extended foliation on $X$.
    \item The (equivalence class of the) measured foliation $\eta$ has many realizations (via isotopy and/or Whitehead moves). Each realization gives rise to a projection map from $\widetilde X$ to $T_\eta$.  \end{enumerate}
\end{remark} 

\color{black}
Considering the maximally stretched loci of  projection maps from $\widetilde X$ to $T_\eta$, we have 
\begin{lemma}\label{lem:optimal:lip:tree:}
 Let $X\in\T(S)$ and $\eta\in \MF(S)$. The maximally stretched locus of any $L(X,\eta)$-Lipschitz equivariant map $f:\widetilde X\to T_\eta$ contains the lift of $\Lambda(X,\eta)$.
\end{lemma}
\begin{proof}
Let $L:=L(X,\eta)$. 
 Suppose to the contrary that there exists a $L$-Lipschitz equivariant map $\phi:\widetilde X\to T_\eta$ whose maximally stretched locus does not contain the lift $\widetilde {\Lambda(X,\eta)}$ of $\Lambda(X,\eta)$. In particular, there exists a point $\tilde p$ in the lift $\widetilde {\Lambda(X,\eta)}$ such that the Lipschitz constant $L_\phi(\tilde p)$ of $\phi$ strictly less than $L$. Let $p$ be the projection of $\tilde p$ under the projection map $\widetilde X\to X$. 

 Notice that the Lipschitz map $f:\widetilde X\to T_\eta$ obtained in Proposition \ref{prop:optimallip:trees} induces a realization $\eta_f$ of $\eta$ on $X$ (cf. Remark \ref{rmk:lip:constuction}(i), $\eta_f$ is a scaling of $\eta_X$ in Remark \ref{rmk:lip:constuction}(i)). Moreover, the lamination $\Lambda(X,\eta)$ is transverse to this realization (cf. \cite[Section 13.1]{PanWolf2022}). 
 
 Recall that $\Lambda(X,\eta)$ is chain recurrent, meaning that $\Lambda(X,\eta)$ can be approximated by geodesic multicurves with respect to the Hausdorff topology. 
 Without loss of generality, we may assume that $\Lambda(X,\eta)$ is connected. Otherwise, we treat each component independently. In particular, we may approximate $\Lambda(X,\eta)$ by a sequence $\alpha_n$ of simple closed geodesics.
 Suppose that $\Lambda(X,\eta)$ has $k\geq 1$ isolated leaves. Then by \cite[Corollary 4.7.1 and Corollary 4.7.2]{CassonBleiler1988}, the lamination $\Lambda(X,\eta)$ is the union of several minimal sublaminations and those isolated leaves. Hence, by fellow traveling $\alpha_n$ on $\Lambda(X,\eta)$, we see that for any positive constant $R$ and any sufficiently large $n$, there exist  $k$ geodesic arcs $\{I^n_1,I^n_2,\cdots,I^n_k\} \subseteq \Lambda(X,
     \eta)$, and $k$ short arcs $\{\delta^n_1,\cdots,\delta_k^n\}$ such that 
 \begin{itemize}
     \item the point $p$ is the center of the (first) segment $I_1^n$,
     \item  each segment $I_i^n$ is of length $\ell(I_i^n)>R$ and is contained in a leaf of $\Lambda(X,
     \eta)$,
     \item each segment $\delta_i^n$ is of length $\ell(\delta_i^n)<1/n$ and is contained in a leaf of the realization $\eta_f$ of $\eta$,
         \item the concatenation $I^n_1*\delta^n_1* I_2^n*\cdots *I_k^n*\delta_k^n$ is a closed curve 
    which is homotopic to $\alpha_n$ and whose length on $X$ satisfies
    \begin{equation*}   \left|\sum_i\left(\ell(I_i^n)+\ell(\delta_i^n)\right)-\ell(\alpha_n)\right|<\frac{1}{n}.
    \end{equation*}
 \end{itemize}
 Here the third property follows from the fact that $\Lambda(X,\eta)$ is transverse to the realization $\eta_f$ of $\eta$. 
(When there are no isolated leaves, i.e. $k=0$, we may instead use a single long geodesic arc $I^n$ and a single short segment $\delta^n$.) The fourth property follows from the Anosov Closing Lemma \cite[Lemma 13.1]{Anosov1969}.

 Since the map $f$ is an affine map of factor $L$ on the lift of $\Lambda(X,\eta)$ (cf. Proposition \ref{prop:optimallip:trees}) and that $\delta_i^n$ is a leaf segment of the realization $\eta_f$ of $\eta$, while each $I_i^n \subset \Lambda(X, \eta)$ is transverse to $\eta$, we see that 
 \begin{equation}\label{eq:length:alpha:tree}
     2\cdot i(\alpha_n,\eta)=L\cdot \sum_i \ell(I^n_i).
 \end{equation}

  Next, we consider the map $\phi:\widetilde X\to T_\eta$ introduced from the beginning of the proof. 
  The assumption that the Lipschitz constant $L_\phi(p)$ of $\phi$ at $p$ is strictly less than $L$ implies that there exists a fixed interval $I\subset I^n_1$ centered at $p$ (independent of $n$) such that the length $\mathrm{Length}(\phi(I))$ of $\phi(I)$ with respect to the metric on $T_\eta$  satisfies: 
 \begin{equation*}
   \mathrm{Length}(\phi(I))<\frac{L_\phi(p)+L}{2} \cdot \ell(I),
 \end{equation*}
 where $\ell(I)$ denotes the hyperbolic length of $I$ on $X$. Since $\phi$ is $L$-Lipschitz, it follows that the length of the image of the closed curve $I^n_1*\delta^n_1* I_2^n*\cdots *I_k^n*\delta_k^n$ under the map $\phi$ satisfies  
\begin{eqnarray*}
&&\mathrm{Length}(\phi(I^n_1*\delta^n_1* I_2^n*\cdots*I_k^n*\delta_k^n*I_1^n))\\
&=&\sum_{i=1}^k \mathrm{Length}(\phi(I^n_i))+ {\sum_{i=1}^k} \mathrm{Length}(\phi(\delta^n_i))\\
   &\leq & L\cdot \ell(I^n_1\backslash I)+ \mathrm{Length}(\phi(I)) +L\cdot \sum_{i=2}^k \ell(I^n_i) +L{\sum_{i=1}^k} \ell(\delta^n_i)\\
   &\leq & L\cdot \ell(I^n_1\backslash I)+ \frac{L_\phi(p)+L}{2} \cdot \ell(I) +L\cdot \sum_{i=2}^k \ell(I^n_i) +\frac{kL}{n} \\&=&
   L \sum_{i=1}^k \ell(I^n_i) +\frac{L_\phi(p)-L}{2} \cdot \ell(I) +\frac{kL}{n} 
\end{eqnarray*}
This implies that, for $n$ sufficiently large, because, $L_\phi(p)<L$, we have the strict inequality, 
\begin{equation*}
    \mathrm{Length}\left(\phi(I^n_1*\delta^n_1* I_2^n*\cdots*I_k^n* \delta_k^n)\right)< L \sum_{i=1}^k \ell(I^n_i).
\end{equation*}
Since $I^n_1*\delta^n_1* I_2^n*\cdots *I_k^n*\delta_k^n$ is homotopic to $\alpha_n$, it follows that
\begin{equation*}
    2\cdot i(\alpha_n,\eta)\leq  \mathrm{Length}\left(\phi(I^n_1*\delta^n_1* I_2^n*\cdots*I_k^n*\delta_k^n)\right)< L \sum_{i=1}^k \ell(I^n_i). 
\end{equation*}
This contradicts \eqref{eq:length:alpha:tree}.
\end{proof}
\begin{remark}
    The proof here extends the argument in \cite[Section 9]{GueritaudKassel2017} where they treated Lipschitz maps between hyperbolic manifolds. The corresponding estimate of $\mathrm{Length}(\phi(I^n_1*\delta^n_1* I_2^n*\cdots*I_k^n*\delta_k^n))$ in \cite{GueritaudKassel2017} was done by applying the Anosov Closing Lemma.
\end{remark}

As an immediate consequence of Proposition~\ref{prop:optimallip:trees} and Lemma~\ref{lem:optimal:lip:tree:}, we find the following.

\begin{corollary}\label{cor:optimal:lip:tree} 
    The lift of the lamination $\Lambda(X, \eta)$ is the intersection of the maximally stretched loci of all $L(X,\eta)$-Lipschitz equivariant maps from $\tilde{X}$ to the tree $T_\eta$.
\end{corollary}

In view of Corollary~\ref{cor:optimal:lip:tree}, we make the following definition.
\begin{definition}\label{def:msl:tree}
    For $X\in\T(S)$ and $\eta\in\MF(S)$.
    \begin{enumerate}
        \item
     The lift $\widetilde{\Lambda(X,\eta)}$ of the chain recurrent geodesic lamination $\Lambda(X,\eta)\subsetneq X$ is called the \emph{maximally stretched lamination} from $\widetilde X$ to the tree $T_\eta$.
     \item The chain recurrent geodesic lamination  $\Lambda(X,\eta)$ is called   the \emph{maximally stretched lamination} from $X$ to (the projective class of) $\eta$.
     \end{enumerate}
\end{definition}   

 Regarding the continuity of maximally stretched laminations, we have
 \begin{proposition}[\cite{PanWolf2022}, Proposition 13.11]
\label{prop:maximally:stretched:lamination:tree} 
 Let $X,X_n\in\T(S)$ and $\eta,\eta_n\in\PMF(S)$ be  such that $\lim\limits_{n\to\infty}X_n=X$ and $\lim\limits_{n\to\infty}\eta_n=\eta$.  Then $\Lambda(X,\eta)$ contains the limit of any convergent subsequence of $\Lambda(X_n,\eta_n)$ with respect to the Hausdorff topology.
 \end{proposition}

Combined with the previous discussion about harmonic stretch lines, this gives
\begin{corollary}
    \label{cor:maximally:stretched:lamination:bnd} 
 Let $X,X_n,Y_n\in\T(S)$ and $\eta\in\PMF(S)$ be such that $\lim\limits_{n\to\infty}X_n=X$ and $\lim\limits_{n\to\infty}Y_n=\eta$.  Then $\Lambda(X,\eta)$ contains the limit of any convergent subsequence of $\Lambda(X_n,Y_n)$ with respect to the Hausdorff topology.
\end{corollary}
\begin{proof}
    By Proposition \ref{prop:cont:compactification2}, we see that $[X_n,Y_n]$ converges to $[X,\eta)$ locally uniformly. Consider the limit $\eta_n$ of $[X_n,Y_n]$ in the Thurston boundary. It then follows from Lemma \ref{lem:Yn:eta0} that $[X_n,\eta_n)$ converges to $[X,\eta)$, yielding that $\eta_n\to\eta$ by Corollary \ref{cor:con:compactification}. Applying Proposition \ref{prop:maximally:stretched:lamination:tree}, we see that $\Lambda(X,\eta)$ contains any geodesic lamination in the limit set of $\Lambda(X_n,Y_n)$ (=$\Lambda(X_n,\eta_n$)) with respect to the Hausdorff topology.  
\end{proof}

Recall from Remark \ref{rem:EnvelopeByStretchLamination} that for any triple $X,Y,Z\in\T(S)$ the union $[X,Z]\cup[Z,Y]$ is a geodesic if and only if the maximally stretched laminations satisfy
\begin{equation*}
 \Lambda(X,Y)=\Lambda(X,Z)\cap \Lambda(Z,Y).
\end{equation*}
We can extend this relation to the setting where $Y$ is replaced by a point in the Thurston boundary. That is
\begin{lemma}\label{lem:msl:bnd}
    Let $X,Z\in\T(S)$ and $\eta\in\PMF(S)$. Then the union $[X,Z]\cup[Z,\eta)$ is a geodesic if and only if the maximally stretched laminations satisfy
    \begin{equation*}
    \Lambda(X,\eta)=\Lambda(X,Z)\cap \Lambda(Z,\eta).
    \end{equation*}
\end{lemma}
\begin{proof}

  (1) For one direction, suppose that $\Lambda(X,\eta)=\Lambda(X,Z)\cap \Lambda(Z,\eta)$, that is, both $[X,Z]$ and $[Z,\eta)$ maximally stretches $\Lambda(X,\eta)$. Hence the union $[X,Z]\cup[Z,\eta)$ is a geodesic.

   (2) For the reverse direction, suppose that the union $[X,Z]\cup[Z,\eta)$ is a geodesic. Let $Z_n\in[Z,\eta)$ be a divergent sequence. Then by Proposition \ref{prop:cont:compactification2}, since $[Z,Z_n]\subsetneq [Z,\eta)$ converges to $[Z,\eta)$, we have $Z_n\to \eta$ in the Thurston compactification. Combined with Corollary \ref{cor:maximally:stretched:lamination:bnd}, this implies that 
    $\Lambda(X,
    \eta)$ contains the Hausdorff limit of any convergent subsequence of $\Lambda(X,Z_n)$. Notice that $[X,Z]\cup[Z,Z_n]$ is a geodesic, yielding that $\Lambda(X,Z_n)=\Lambda(X,Z)\cap\Lambda(Z,Z_n)=\Lambda(X,Z)\cap \Lambda(Z,\eta)$, as $[Z,Z_n] \subset[Z, \eta)$. 
    Therefore, $\Lambda(X,
    \eta)$ contains $\Lambda(X,Z)\cap\Lambda(Z,\eta)$.

    On the other hand, by \cite[Theorem 1.5 and Theorem 1.6]{PanWolf2022}, the harmonic stretch map from $\widetilde X$ to $\widetilde Z$ maximally stretches exactly along the lift of $\Lambda(X,Z)$. By  Proposition \ref{prop:optimallip:trees}, there exists an equivariant Lipschitz map from the universal cover of $Z$ to the tree dual to $\eta$, which maximally stretches exactly the lift of $ \Lambda(Z,\eta)$ to the universal cover. Therefore,  the composition of these two maps gives rise to an equivariant Lipschitz map from the universal cover of $X$ to the tree dual to $\eta$, which maximally stretches exactly the lift of $\Lambda(X,Z)\cap \Lambda(Z,\eta)$ to the universal cover. 
    Hence, by the definition (Definition \ref{def:msl:tree}) of $\Lambda(X,\eta)$, we see that $\Lambda(X,\eta)$ is contained in $\Lambda(X,Z)\cap \Lambda(Z,\eta)$. This completes the proof.
\end{proof}

\subsection{Envelopes and the Thurston compactification} \label{sec:env:ThurstonCompact}
In this subsection, we extend the definition of the envelope to include terminal points on the Thurston boundary. This then invites the question of the accumulation points of this envelope: we describe the locus of these accumulation points, which is often larger than the initially-defined terminus of the envelope.

   \begin{definition}
         For any $X\in\T(S)$ and $\eta\in\PMF(S)$, define the envelope $\env(X,\eta)\subset \T(S)$ to be the union of all geodesic rays which start at $X$ and which converge to $\eta$ in the Thurston compactification.
     \end{definition}
 
  Here by a \enquote{geodesic ray}, we mean (the image of) a continuous map $\iota:[0,\infty)\to\T(S)$ such that for any $0\leq s<t<r<\infty$ we have $\dth(\iota(s),\iota(r))=\dth(\iota(s),\iota(t))+\dth(\iota(t),\iota(r))$. 

  \color{black}\begin{lemma}\label{lem:env:tree}
      Let  $X\in\T(S)$ and $\eta\in\PMF(S)$. 
      \begin{enumerate}
          \item The envelope $\env(X,\eta)$ satisfies
      \begin{eqnarray*}
          \env(X,\eta)&=&\{Z\in\T(S): [X,Z]\cup[Z,\eta) \text{ is a geodesic}\}\\
          \env(X,\eta)&=&\{Z\in\T(S):\Lambda(X,\eta)=\Lambda(X,Z)\cap\Lambda(Z,\eta)\}.
      \end{eqnarray*}
      \item The envelope $\env(X,\eta)$ is a closed subset of $\T(S)$.
      \item For any  $Z\in\env(X,\eta)$, we have
      \begin{equation*}       \env(X,Z)\cup \env(Z,\eta)\subset \env(X,\eta).
      \end{equation*}
      \end{enumerate}
  \end{lemma}

 \begin{proof}

     (i). By Lemma \ref{lem:msl:bnd}, we see that for any $Z\in\T(S)$, the union $[X,Z]\cup[Z,\eta)$ is a geodesic if and only if $\Lambda(X,\eta)=\Lambda(X,Z)\cap\Lambda(Z,\eta)$. In view of this, it suffices to prove the first equation. For any $Z\in\T(S)$ with the property that $[X,Z]\cup [Z,\eta)$ is a geodesic ray, since $[Z,\eta)$ converges to $\eta$ in the Thurston compactification, we see that $Z$ is contained in the envelope $\env(X,\eta)$.  For the converse, let $Z\in\env(X,\eta)$. Then, by definition, there is a  geodesic ray $\mathscr R$ 
     through $Z$ which starts at $X$ and converges to $\eta$ in the Thurston compactification.  Let $Y_n\in \mathscr R$ be a divergent sequence such that $\mathscr R$ proceeds from $Z$ to $Y_n$. Then by Lemma \ref{lem:Yn:eta2}, we see that $[X,Y_n]$ converges to $[X,\eta)$ locally uniformly.  On the other hand, since the geodesic ray $\mathscr R$ proceeds from $Z$ to $Y_n$, it follows that $Z\in\env(X,Y_n)$. Hence, the union  
     $[X,Z]\cup [Z,Y_n]$ is a geodesic from $X$ to $Y_n$. Letting $n\to\infty$, we see that $[X,Z]\cup [Z,\eta)$ is also a geodesic. This completes the proof. 

     (ii). Let $Z_n\in\env(X,\eta)$. Then by the item (i), we see that the union $[X,Z_n]\cup[Z_n,\eta)$ is a geodesic. Suppose that $Z_n$ converges to a point $Z$ in $\T(S)$. Then by Proposition \ref{prop:PanWolf2022} and Corollary \ref{cor:con:compactification}, we see that $[X,Z_n]\cup [Z_n,\eta)$ converges to $[X,Z]\cup[Z,\eta)$ locally uniformly, yielding that $[X,Z]\cup[Z,\eta)$ --as the limit of a sequence of geodesics-- is a also geodesic. Hence by item (i), the point $Z$ is also contained in the envelope $\env(X,\eta)$. 

     (iii). For any $Z'\in \env(X,Z)$ and any $Z''\in\env(Z,\eta)$, we have  
      \begin{eqnarray*}
         \Lambda(X,\eta)&=&\Lambda(X,Z)\cap\Lambda(Z,\eta)  \qquad (\text{by Lemma \ref{lem:msl:bnd}})\\
         &=& \Lambda(Z,Z')\cap\Lambda(Z',Z) \cap \Lambda(Z,Z'')\cap\Lambda(Z'',\eta), 
     \end{eqnarray*}
     where the second equation follows from Corollary \ref{cor:env in terms of stretch lams} and Lemma \ref{lem:msl:bnd}.
     Hence, the union
     \begin{equation*}
      [X,Z']\cup[Z',Z]\cup [Z,Z'']\cup [Z'',\eta)
     \end{equation*}
     is a geodesic converging to $\eta$. Therefore, both $Z'$ and $Z''$ are contained in $\env(X,\eta)$. 
 \end{proof}

\color{black}
We aim to describe the accumulation set of $\env(X, \eta)$. This will require some preparation. We begin with an extension of Lemma \ref{lem:concatenation} to the current setting.  (We will also use this lemma when we discuss the notions of extendability in the setting of these infinite envelopes.)

\begin{lemma}\label{lem:concatenation:tree}
	Let $S$ be an orientable closed surface of genus at least two. Let $X\in\T(S)$ and $\eta\in\PMF(S)$. Let  $\lambda$ be an arbitrary geodesic lamination which contains $\Lambda(X,\eta)$. Then 
	\begin{enumerate}
		\item for any piecewise harmonic stretch line $\mathrm{PHSL}(X,\lambda)$  through $X$ which maximally stretches $\lambda$,  the intersection 
	\begin{equation*}
		\mathrm{PHSL}(X,\lambda)\cap \env(X,\eta)
	\end{equation*}
 is a non-trivial oriented geodesic segment/ray with initial point $X$;
 \item for any piecewise harmonic stretch line $\mathrm{R}$  ending at $\eta$ which maximally stretches $\lambda$,  the intersection 
	\begin{equation*}
		\mathrm{R}\cap \env(X,\eta)
	\end{equation*}
 is a non-trivial oriented geodesic ray with terminal point $\eta$.
	\end{enumerate}
	 	\end{lemma}
\begin{proof}
    (i) Let $Z\in [X,\eta)$ be a point distinct from $X$. It follows from Lemma \ref{lem:env:tree}(iii) that $\env(X,Z)\subsetneq\env(X,\eta)$. By Lemma \ref{lem:concatenation}, the intersection $\PHSL(X,\lambda)\cap\env(X,Z)$ is a non-trivial oriented geodesic segment with initial point $X$. Hence the intersection $\PHSL(X,\lambda) \cap \env(X,\eta)$ is nontrivial. It remains to show that it is an oriented geodesic segment/ray. Let $Z'\in \PHSL(X,\lambda)\cap\env(X,\eta)$ be an arbitrary point distinct from $X$. Let $\PHSL(X,\lambda)|_{[X,Z']}$ be the oriented geodesic subsegment of $\PHSL(X,\lambda)$ with initial point $X$ and terminal point $Z'$. Then  $$\PHSL(X,\lambda)|_{[X,Z']}\subset \PHSL(X,\lambda)\cap\env(X,Z')\subsetneq \PHSL(X,\lambda)\cap\env(X,\eta),$$ where the second inclusion follows from Lemma \ref{lem:env:tree}(iii).
    Therefore, the intersection $\PHSL(X,\lambda)\cap\env(X,\eta)$ is an oriented geodesic segment/ray of $\PHSL(X,\lambda)$.

  (ii). The proof is essentially that of Lemma \ref{lem:concatenation} (modified for example to avoid the use of formulas that only hold for compact sets in $\T(S)$, for instance Equation \eqref{eq:ratio:ZX} with $X$ or $Z'$ being at infinity). Consider the horocyclic (partial) measured foliation near $\Lambda(X,\eta)\subset\lambda$. From the construction of harmonic stretch lines, we see that the horocyclic (partial) measured foliation near $\Lambda(X,\eta)$ on $X$ is a scaling of the restriction of $\eta$.  Similarly, since by assumption $\mathrm{R}$ maximally stretches $\lambda$ and converges to $\eta$ in the Thurston compactification, for any $Z\in\mathrm{R}\cap\T(S)$, the horocyclic (partial) measured foliation near $\Lambda(X,\eta)\subset\lambda$ on $Z$ is also a scaling of the restriction of $\eta$.  Therefore, there exists a Lipschitz homeomorphism from a neighbourhood of $\Lambda(X,\eta)$ on $X$ to a neighbourhood of $\Lambda(X,\eta)$ on $Z$, which maps leaves of the horocyclic partial foliation on $X$ to leaves of  the horocyclic partial foliation on $Z$  and takes leaves of $\Lambda(X, \eta)$ on $X$ linearly to corresponding leaves of $\Lambda(X,\eta)$ on $Z$. In particular, this gives a Lipschitz homeomorphism from a neighbourhood of $\Lambda(X,\eta)$ on $X$ to a neighbourhood of $\Lambda(X,\eta)$ on $Z$ which maximally stretches exactly along $\Lambda(X,\eta)$. On the other hand, by Corollary \ref{cor:maximally:stretched:lamination:bnd}, for any $Z$ sufficiently close to $\eta$ in the Thurston compactification, the maximally stretched lamination $\Lambda(X,Z)$ is contained in a small neighbourhood of $\Lambda(X,\eta)$ on $X$, for example $\Lambda(X,Z)$ could be contained in that neighborhood described just above on which we have already decided that the maximal stretch lamination is $\Lambda(X,\eta)$.
  Thus, together with the previous discussion, this implies that for any $Z\in\mathrm{R}$ sufficiently close to $\eta$, the maximally stretched lamination $\Lambda(X,Z)$ is exactly $\Lambda(X,\eta)$, proving that $[X,Z]\cup[Z,\eta)$ is a geodesic since $\Lambda(Z,\eta)=\lambda$ contains $\Lambda(X,\eta)$. Hence, by Lemma \ref{lem:env:tree}, we see that $[Z,\eta)\subset\env(X,\eta)$.
\end{proof}

Next, we characterize the accumulation points of $\env(X, \eta)$ on the Thurston boundary  by proving the following companion to Lemma \ref{lem:env:tree}. Let $\overline{\env(X,\eta)}^{Th}$ be the closure of $\env(X,\eta)$ in the Thurston compactification. Before stating the proposition, let us fix a notation. For any $\eta\in\PMF(S)$ and any $X\in\T(S)$, we denote by $\eta_X$ the particular representation of $\eta$ which is the extension of the pushforward of the vertical foliation of Hopf differential of the harmonic map defining the harmonic stretch ray $[X,\eta)$ (cf. Remark \ref{rmk:lip:constuction}(i)).
In particular,  by \cite[Lemma 7.8]{PanWolf2022}, we have 
\begin{equation}\label{eq:length:intersection}
  \ell_\alpha(X)=2i(\eta_X,\alpha) 
\end{equation}
for any measured lamination $\alpha$ supported on $\Lambda(X,\eta)$. On the other hand, by \cite[Equation 3.2]{PanWolf2022} we have $\ell_\alpha(X)\geq 2i(\eta_X,\alpha)$ for any measured lamination $\alpha$ on $X$. Combined with \eqref{eq:length:intersection}, we have
\begin{equation}\label{eq:lip=1}
  L(X, \eta_X) = 1 
\end{equation}
 where we refer to \eqref{eq:L(X,eta)} for the definition of $L$.  Furthermore, for $Y\in\env(X,\eta)$ at distance $t$ from $X$, we have 
\begin{equation}\label{eq:eta:YX}
    \eta_{Y} = e^t\eta_X.
\end{equation}
To see this, note that $\eta_{Y_t}$ is a scaling of $\eta_X$. By  $Y\in\env(X,\eta)$ and  Lemma \ref{lem:env:tree}(i), we see that $\Lambda(X,\eta)\subset \Lambda(X,Y)$. Hence, for any measured lamination $\alpha_0$ supported on $\Lambda(X,\eta)$, we have $2i(\eta_{Y_t},\alpha_0)=\ell_{\alpha_0}(Y_t)=e^t \ell_{\alpha_0}(X)=2e^ti(\eta_X,\alpha_0)$. This proves \eqref{eq:eta:YX}.

\begin{lemma}\label{lem:env:mu}
     Let $X\in\T(S)$ and $\eta,\mu\in\PMF(S)$. If $\mu\in\overline{\env(X,\eta)}^{Th}$, then $[X,\mu)\subset\env(X,\eta)$ and $\Lambda(X,\eta)\subset\Lambda(X,\mu)$.
\end{lemma}
\begin{proof}
   Since $\mu$ is an accumulation point of $\env(X,\eta)$ in the Thurston boundary, there exists a sequence $Y_n\in\env(X,\eta)$ which converges to $\mu$ in the Thurston compactification. By Corollary \ref{cor:con:compactification}, the harmonic stretch segment $[X,Y_n]$ converges to $[X,\mu)$.  Since $[X,Y_n]\subsetneq \env(X,\eta)$, by Lemma \ref{lem:env:tree}(ii), we see that the whole harmonic stretch ray $[X,\mu)$ is contained in $\env(X,\eta)$. Hence, by Lemma \ref{lem:msl:bnd},  for any $Y\in[X,\mu)$, we have
   $\Lambda(X,\eta)=\Lambda(X,Y)\cap\Lambda(Y,\eta)$. In particular, $\Lambda(X,\eta)\subset \Lambda(X,Y)=\Lambda(X,\mu)$.  This proves the lemma. 
\end{proof}

Let $Y_t\in[X,\mu)$ be the point such that $\dth(X,Y_t)=t$. Define 
\begin{eqnarray}
\nonumber	\mathbf{L}:[0,\infty)&\longrightarrow&\R 
	\\ t&\longmapsto &\max\limits_{\alpha\in\ML(S)}\frac{2i(\eta_X,\alpha)}{e^{-t}\ell_{\alpha}(Y_t)}. \label{eq:L:def}
\end{eqnarray}

\begin{proposition}\label{prop:msl:bnd2}
    Let $X\in\T(S)$ and $\eta,\mu\in\PMF(S)$. Then the following are equivalent. 
  \begin{enumerate}
      \item $\mu\in\overline{\env(X,\eta)}^{Th}$.
      \item  $\Lambda(X,\eta)\subset\Lambda(X,\mu)$, and $\max\limits_\alpha\frac{i(\eta_X,\alpha)}{i(\mu_X,\alpha)}$ is achieved by any measured lamination supported on $\Lambda(X,\eta)$.
      \item $\max\limits_\alpha\frac{i(\eta_X,\alpha)}{i(\mu_X,\alpha)}=1$.
      \item $\mathbf{L}(t)$ is constant.
      \end{enumerate} 
  \end{proposition}
The proof is a consequence of the following lemma.
 
\begin{lemma}\label{lem:monotone}
	Let $X\in\T(S)$ and $\eta,\mu\in\PMF(S)$. As above, we set $Y_t \in [X,\mu)$ to be the point on $[X,\mu)$ with $\dth(X,Y_t)=t$. Then we have the following.
    \begin{enumerate}
        \item The function $\mathbf L(t)$ is increasing in $t$,
        \item  $\mathbf{L}(0)=1$ and is realized by any measured lamination supported on $\Lambda(X,\eta)$.
        \item  $\lim\limits_{t\to\infty}\mathbf{L}(t)= \max\limits_{\alpha\in\ML(S)}\frac{i(\eta_X,\alpha)}{i(\mu_X,\alpha)}$.  
        \item If, moreover, $Y_{t_0}\in\env(X,\eta)$ for some $t_0>0$, then $\mathbf{L}(t_0)=1$  and  is attained by any measured lamination  supported on $\Lambda(X,\eta)$.  
      \end{enumerate}
\end{lemma}

\begin{remark}
  Note that if $\alpha$ is supported on $\mu$, then for $Y_t \in [X, \mu)$,  the quantity $e^{-t}\ell_{\alpha}(Y_t)$ is constant in $t$, as the geodesic representative of $\alpha$ remains in $\Lambda(X,\mu)$ and so is  stretched exactly by the stretching factor $e^t$ as $Y_t$ travels along $[X,\mu)$; moreover, that quantity limits on $2i(\mu,\alpha)$.
\end{remark}

\begin{proof}
	(i) Let $0\leq s<t<\infty$. Let $\alpha_s$ be a measured lamination realizing $\mathbf{L}(s)$, that is,
	\begin{equation*}
		\mathbf{L}(s)=\frac{2i(\eta_X,\alpha_s)}{e^{-s}\ell_{\alpha_s}(Y_s)}. 
	\end{equation*}
	Since $\dth(Y_s,Y_t)=t-s$, it follows that for any measured lamination $\alpha$ we have, 
    \begin{equation}\label{eq:length:st}
        \ell_{\alpha}(Y_t)\leq e^{t-s}\ell_{\alpha}(Y_s) \quad \text{ and so} \quad e^{-t}\ell_{\alpha}(Y_t)\leq e^{-s}\ell_{\alpha}(Y_s).
    \end{equation}
     Hence,
	\begin{equation*}
		\mathbf{L}(s)=\frac{2i(\eta_X,\alpha_s)}{e^{-s}\ell_{\alpha_s}(Y_s)}\leq \frac{2i(\eta_X,\alpha_s)}{e^{-t}\ell_{\alpha_s}(Y_t)}\leq \max_{\alpha\in\ML(S)}\frac{2i(\eta_X,\alpha)}{e^{-t}\ell_{\alpha}(Y_t)}=\mathbf{L}(t).
	\end{equation*}
	This proves that $\mathbf L(t)$ is increasing in $t$. 
    
    (ii) This  follows from the fact $Y_0=X$, \eqref{eq:length:intersection}, and \eqref{eq:lip=1}.
	
	(iii) By the monotonicity in item (i) of the current lemma, we see that the limit $\lim\limits_{t\to\infty}\mathbf{L}(t)$ exists in $[0,\infty]$. It is straightforward to estimate that limit.  Note that, by \cite[Proposition 7.11]{PanWolf2022}, especially the last displayed equation in the proof of \cite[Proposition 7.11]{PanWolf2022}, we have $e^{-t}Y_t\to 2\mu_X$,  Then, 
  $$\lim_{t\to\infty}\mathbf{L}(t)\geq \lim_{t\to\infty} \frac{2i(\eta_X,\alpha)}{e^{-t}\ell_\alpha(Y_t)}=\frac{i(\eta_X,\alpha)}{i(\mu_X,\alpha)},$$
 yielding that 
 \begin{equation}\label{eq:limit:L:low}
     \lim\limits_{t\to\infty}\mathbf{L}(t)\geq \max\limits_{\alpha\in\ML(S)}\frac{i(\eta_X,\alpha)}{i(\mu_X,\alpha)}
 \end{equation}
On the other hand, combining the monotonicity estimate \eqref{eq:length:st} and the fact that $e^{-t}Y_t\to 2\mu_X$ mentioned above, we see that $e^{-t}\ell_\alpha(Y_t)\geq 2i(\mu_X,\alpha)$ for any measured lamination $\alpha$ on $S$. Accordingly, 
\begin{eqnarray}\label{eq:L:up}
 \nonumber   \mathbf{L}(t)&=&\max_{\alpha\in\ML(S)}\frac{2i(\eta_X,\alpha)}{e^{-t}\ell_\alpha(Y_t)}\\
 \nonumber&\leq &\max_{\alpha\in\ML(S)}\frac{2i(\eta_X,\alpha)}{2i(\mu_X,\alpha)}\\
 &=&\max_{\alpha\in\ML(S)}\frac{i(\eta_X,\alpha)}{i(\mu_X,\alpha)}. \label{eq:L:up}
\end{eqnarray}
Combined with \eqref{eq:limit:L:low}, this proves (iii).

 (iv) Suppose that $Y_{t_0}\in \env(X,\eta)$.  Then by \eqref{eq:eta:YX}, we see that   $\eta_{Y_{t_0}} = e^{t_0} \eta_X$. By Lemma \ref{lem:env:tree}(i), we have $\Lambda(X,\eta)\subset \Lambda(Y_{t_0},\eta)$. Let $\alpha_0$ be an arbitrary measured lamination  supported on $\Lambda(X,\eta)\subset\Lambda(Y_{t_0},\eta)$. Then  Proposition~\ref{prop:optimallip:trees} implies that 
\begin{equation*}
      \frac{2i(\eta_{Y_{t_0}},\alpha)}{\ell_{\alpha}(Y_{t_0})}= \frac{2i(\eta_{X},\alpha)}{e^{-t_0}\ell_{\alpha}(Y_{t_0})}
   \end{equation*}
   is maximized by $\alpha_0$. Thus, since $Y_0=X$, we have
   \begin{eqnarray}
   \nonumber	\mathbf{L}(t_0)&=&\frac{2i(\eta_X,\alpha_0)}{e^{-t_0}\ell_{\alpha_0}(Y_{t_0})}\\  \nonumber&=&\frac{2i(\eta_X,\alpha_0)}{\ell_{\alpha_0}(Y_0)}\\ \nonumber&=&\frac{2i(\eta_X,\alpha_0)}{\ell_{\alpha_0}(X)}\\
    &=&\frac{i(\eta_X,\alpha_0)}{i(\mu_X,\alpha_0)} \qquad\qquad\qquad(\text{by }\eqref{eq:lip=1})\label{eq:L:alpha0}\\
    &=& \frac{\ell_{\alpha_0}(X)}{\ell_{\alpha_0}(X)} \qquad\qquad\qquad(\text{by }\eqref{eq:lip=1})\nonumber\\
    &=& 1  \label{eq:L:alpha0:1}.
   \end{eqnarray}
   This completes the proof.
 \end{proof}

   The above passage leads to an elementary proof of the following theorem: 
   \begin{theorem}[\cite{Papadopoulos1991, Walsh2014}] \label{thm:GeodesicLimit} 
       Every Thurston geodesic ray converges to a unique point in the Thurston boundary.
   \end{theorem}
   \begin{proof}
    Let $X_t$ be an arbitrary geodesic ray in the Thurston metric. Then the normalized length function $e^{-t}\ell_\alpha(X_t)$ is decreasing in $t\geq0$, hence is convergent as $t\to\infty$. To complete the proof, it remains to show that the limit $\lim\limits_{t\to\infty}e^{-t}\ell_{\alpha_0}(X_t)$ is positive for some measured lamination $\alpha_0$. To see this, notice that $\Lambda(X_0,X_s)\supset \Lambda(X_0,X_t)$ for any $0\leq s\leq t$.  Combined with the fact that  $\Lambda(X_0,X_t)$ is a closed subset of the underlying surface $S$ for any $t\geq0$, this implies that $\bigcap_{t\geq0}\Lambda(X_0,X_t)$ is a (nonempty) geodesic lamination. Hence, for any measured lamination $\alpha_0$ supported on $\bigcap_{t\geq0}\Lambda(X_0,X_t)$, the geodesic length of $\alpha_0$ is stretched exactly by the factor $e^t$ from $X_0$ to $X_t$. Equivalently, the scaled function $e^{-t}\ell_{\alpha_0}(X_t)\equiv \ell_{\alpha_0}(X_0)$ is a positive constant. This completes the proof.
     \end{proof}
   \begin{remark}\label{rem:GeodesicLimit}
        (Miyachi \cite[Lemma 1]{Miyachi} applied a similar idea to prove that every Teichm\"uller geodesic ray converges in the Gardiner-Masur compactification of the \tec space.
\end{remark}

\begin{proof}[Proof of Proposition \ref{prop:msl:bnd2}]
 \underline{(i) $\Rightarrow$ (ii)}. Suppose that $\mu$ is an accumulation point of $\env(X,\eta)$. By Lemma \ref{lem:env:mu}, we see that $[X,\mu)\subset\env(X,\eta)$. Combined with Lemma \ref{lem:env:tree}, this implies that $\Lambda(X,\eta)\subset\Lambda(X,\mu)$.  Furthermore, by Lemma \ref{lem:monotone}(iv), we infer that $\mathbf{L}(t)\equiv1$ for all $t>0$. Combining with Lemma \ref{lem:monotone}(iii), we see that $\max\limits_{\alpha\in\ML(S)}\frac{i(\eta_X,\alpha)}{i(\mu_X,\alpha)}=1$. It then follows from \eqref{eq:L:alpha0} and \eqref{eq:L:alpha0:1} that  $\max\limits_\alpha\frac{i(\eta_X,\alpha)}{i(\mu_X,\alpha)}$ is realized by any measured lamination supported on $\Lambda(X,\eta)$.

\underline{(ii) $\Rightarrow$ (iii)}. 
     By \eqref{eq:length:intersection}, for any measured lamination $\alpha_0$ supported on $\Lambda(X,\eta)\subset\Lambda(X,\mu)$, we have $2i(\eta_X,\alpha)=\ell_\alpha(X)=2i(\mu_X,\alpha)$. The assumption that $\max\limits_\alpha\frac{i(\eta_X,\alpha)}{i(\mu_X,\alpha)}$ is achieved by any measured lamination supported on $\Lambda(X,\eta)$ then implies that $\max\limits_\alpha\frac{i(\eta_X,\alpha)}{i(\mu_X,\alpha)}=\frac{i(\eta_X,\alpha_0)}{i(\mu_X,\alpha_0)}=\frac{\ell_{\alpha_0}(X)}{\ell_{\alpha_0}(X)}=1$.

     \underline{(iii) $\Rightarrow$ (iv)}. 
     Suppose $\max\limits_{\alpha\in\ML(S)}\frac{i(\eta_X,\alpha)}{i(\mu_X,\alpha)}=1$.  By 
     Lemma \ref{lem:monotone}(i-iii),  we see that 
     \begin{equation}
         \mathbf{L}(t)\equiv\max\limits_{\alpha\in\ML(S)}\frac{i(\eta_X,\alpha)}{i(\mu_X,\alpha)}=1
     \end{equation}  for all $t\geq0$. 

       \medskip \underline{(iv) $\Rightarrow$ (i)}. By Lemma \ref{lem:monotone}(ii), we see that $\mathbf{L}(t)\equiv 1$ for all $t\geq0$. 
     Now for $t>0$, by Lemma \ref{lem:monotone}(ii),  \eqref{eq:length:st} and the fact $Y_0=X$, we see that for any measured lamination $\alpha_0$ supported on $\Lambda(X,\eta)$, 
     \begin{equation*}
         1=\mathbf{L}(0)= \frac{2i(\eta_X,\alpha_0)}{\ell_{\alpha_0}(Y_{0})}\leq \frac{2i(\eta_X,\alpha_0)}{e^{t}\ell_{\alpha_0}(Y_t)} \leq\mathbf{L}(t)=1.
     \end{equation*}
     From this we infer that $\mathbf{L}(t)$ is attained by any measured lamination $\alpha_0$ supported on $\Lambda(X,\eta)$, meaning that $\alpha_0$ is maximally stretched along the ray $[Y_t,\eta)$ for any $t\geq0$.  Since by assumption $\alpha_0$ is contained in $\Lambda(X,\eta)\subset\Lambda(X,\mu)$, it is also maximally stretched along $[X,Y_t]$. Hence, the union $[X,Y_t]\cup [Y_t,\eta)$ is a geodesic  with respect to the Thurston metric. According to the item (i) in Lemma \ref{lem:env:tree}, this proves that $Y_t\in\env(X,\eta)$ for any $t\geq0$.  Since $Y_t\to\mu$ in the Thurston compactification, this implies that $\mu$ is an accumulation point of $\env(X,\eta)$. 
\end{proof}

Here are some direct consequences of Proposition \ref{prop:msl:bnd2}. 
\begin{lemma}\label{rmk:env:accumulation}

  Let $X\in\T(S)$ and $\eta,\mu\in\PMF(S)$. Then we have the following.
  \begin{enumerate}
      \item  If $\mu\in\overline{\env(X,\eta)}^{Th}\cap\PMF(S)$, then the maximum     {$\max\limits_\alpha\frac{i(\eta_X,\alpha)}{i(\mu_X,\alpha)}$}  is finite, meaning that $\eta$ is absolutely continuous with respect to $\mu$. (For the relation between $\eta,\mu$ and $\eta_X$, $\mu_X$, we refer to Remark \ref{rmk:lip:constuction}(i).)         
       \item  The intersection $\overline{\env(X,\eta)}^{Th}\cap\PMF(S)$ is a star over $\eta$, that is, if $\mu\in\overline{\env(X,\eta)}^{Th}$ then $s\eta+(1-s)\mu\in\overline{\env(X,\eta)}^{Th}$ for any $0\leq s\leq 1$.   
        
        \item  For any $X\in\T(S)$ and any pair of distinct equivalence classes of measured laminations $\mu$ and $\eta$ in $\PMF(S)$,  we always have $\mu\notin\overline{\env(X,\eta)}^{Th}$ or $\eta\notin\overline{\env(X,\mu)}^{Th}$.
        \end{enumerate}
    \end{lemma}
    \begin{proof}
    (i) This follows directly from Proposition \ref{prop:msl:bnd2}(ii). 
    
    (ii) By Proposition \ref{prop:msl:bnd2}, we see that $\max\limits_{\alpha\in\ML(S)}\frac{i(\eta_X,\alpha)}{i(\mu_X,\alpha)}=1$ and is realized by any measured lamination supported on $\Lambda(X,\eta)$. Note that $\eta_X$ is a representative of $\eta$ and that $\mu_X$ is a representative of $\mu$, so the representative $(s\eta+(1-s)\mu)_X$ (cf. Remark \ref{rmk:lip:constuction}(i)) of $s\eta+(1-s)\mu$ is a linear sum of $s'\eta_X+s''\mu_X$ for some non-negative numbers $s'$ and $s''$. Then 
    \begin{eqnarray*}
       && \max_\alpha \frac{i(\eta_X,\alpha)}{i((s\eta+(1-s)\mu)_X,\alpha)}\\
       &=&\max_\alpha \frac{i(\eta_X,\alpha)}{s'i(\eta_X,\alpha)+s''i(\mu_X,\alpha)}\\
       &=&\max_\alpha \frac{\frac{i(\eta_X,\alpha)}{i(\mu_X,\alpha)}}{s'\frac{i(\eta_X,\alpha)}{i(\mu_X,\alpha)}+s''}\\
       &=&\frac{\max\limits_\alpha\frac{i(\eta_X,\alpha)}{i(\mu_X,\alpha)}}{s'\max\limits_\alpha\frac{i(\eta_X,\alpha)}{i(\mu_X,\alpha)}+s''}
    \end{eqnarray*}
    is also realized by any measured lamination $\alpha_0$ supported on $\Lambda(X,\eta)$. Combined with the fact $\Lambda(X,\eta)\subset\Lambda(X,\mu)$ (Lemma \ref{lem:env:mu}) and Proposition \ref{prop:msl:bnd2}(ii), this implies that $s\eta+(1-s)\mu\in\overline{\env(X,\eta)}^{Th}$. 
    
        (iii)  Suppose to the contrary that $\mu\in\overline{\env(X,\eta)}^{Th}$ and $\eta\in\overline{\env(X,\mu)}^{Th}$. Then by Proposition \ref{prop:msl:bnd2}, we see that $\max\limits_\alpha\frac{i(\eta_X,\alpha)}{i(\mu_X,\alpha)}=\max\limits_\alpha\frac{i(\mu_X,\alpha)}{i(\eta_X,\alpha)}=1$, yielding that $\eta_X=\mu_X$, hence $\eta=\mu$. 
    \end{proof}

    \begin{remark}
    Regarding accumulation points in $\PMF(S)$ of the envelope $\env(X,\eta)$, there are two cases where $\env(X,\eta)$ has exactly one such accumulation point (in $\PMF(S)$): either (i) $\Lambda(X,\eta)$ is maximally CR, or (ii) $\eta$ is maximal and uniquely ergodic. For case (i) the envelope $\env(X,\eta)$ is the harmonic stretch ray $[X,\eta)$ (See Corollary \ref{cor:unique:geodesic:tree}), hence $\eta$ is the only accumulation point in $\PMF(S)$ of $\env(X,\eta)$. For case (ii), any measured lamination, with respect to which $\eta$ is absolutely continuous, is a scaling of $\eta$. Hence, by  Remark~\ref{rmk:env:accumulation}(1), $\eta$ is the only accumulation point in $\PMF(S)$ of $\env(X,\eta)$. On the other hand, the following proposition indicates that there exists some envelope which contains more than one accumulation point in $\PMF(S)$.
    \end{remark}
    \begin{proposition}\label{prop:accumulation}
    	For any simple closed curve $\alpha$, there exist a hyperbolic surface $X\in\T(S)$ and a measured lamination $\mu$ whose support strictly contains $\alpha$   such that $\mu\in\overline{\env(X,\alpha)}^{Th}$.
    \end{proposition}

We defer the proof to the end of this subsection. Before that, we need two lemmas.

\begin{lemma}\label{lem:alpha:ray} Let $\alpha$ be a simple closed curve on $S$ and let $\lambda$ be a simple closed curve which intersects $\alpha$ once if $\alpha$ is non-separating or twice if $\alpha$ is separating. Then there exists a harmonic stretch line $\mathscr R$ in $\T(S)$ which converges to $\alpha$ and whose maximally stretched lamination is $\lambda$.  
\end{lemma}
    \begin{proof}  
 We extend $\lambda$ to a measured lamination $\lambda+\beta$ so that 
    \begin{itemize}
    \item  $\beta$ is disjoint from $\lambda$, and 
    \item $\lambda+\beta$ and $\alpha$ fill up the underlying surface.
    \end{itemize}
    (We may take $\beta$ to be a measured lamination that fills up $X\setminus\lambda$.) Consider the sequence of quadratic differentials $\Phi_n$ whose horizontal and vertical foliations have transverse measures $n\lambda+\beta$  and $\alpha$,  respectively, i.e. $\Phi_n$ is a quadratic differential on the unique Riemann surface $X_n$ which admits a quadratic differential with the specified horizontal and vertical foliations. Let $\HR_n$ be the sequence of harmonic map rays defined by $\Phi_n$. 
        
     For $X_n$ that Riemann surface underlying $\Phi_n$, let $Y_n\in\HR_n$ be the hyperbolic surface such that the Hopf differential of the harmonic map from $X_n$ to $Y_n$ is exactly $\Phi_n$.    To prove the lemma,  by the definition of harmonic stretch lines, it suffices to prove that $\{Y_n\}$ contains a convergent subsequence.  Recall that finite  length spectrum rigidity (see \cite{FLP} for instance) implies that there exists finitely many simple closed curves $\alpha_1,\alpha_2,\cdots,\alpha_k$ such that if two hyperbolic surfaces $V$ and $W$ in $\T(S)$ satisfy $\ell_{\alpha_i}(V)=\ell_{\alpha_i}(W)$ for every $i$ then $V=W$. Therefore, to prove that $\{Y_n\}$ contains a bounded subsequence, it suffices to prove that there exists a subsequence $n_j$ such that $\{\ell_{\alpha_i}(Y_{n_j})\}_{j}$  is bounded for each $1\leq i\leq k$. Notice that  $\|\Phi_n\|=i(\alpha,n\lambda+\beta)\to\infty$.  Consider the singular flat metrics induced by $\Phi_n$. For any homotopically non-trivial simple closed curve $\xi$ on $X\setminus\lambda$, the $|\Phi_n|$-length  of $\xi$ is at least $\max\{i(\alpha,\xi),i(n\lambda+\beta,\xi)\}\geq \frac{i(\alpha,\xi)+i(\beta,\xi)}{2}$. Since $\alpha$ and $\beta+\lambda$ fill up the whole surface, it follows that $\alpha$ and $\beta$ fill up the subsurface $X\setminus\lambda$. Then there exists a positive constant depending on $\alpha$ and $\beta$ such that $\frac{i(\alpha,\xi)+i(\beta,\xi)}{2}>\delta$ for any homotopically non-trivial and non-peripheral simple closed curve $\xi$ on $X\setminus\lambda$. Accordingly, the $|\Phi_n|$-length of $\xi$ is at least $\delta$. Note that $\Phi_n$ contains a maximal horizontal cylinder whose core curve is homotopic to $\lambda$ and whose circumference and height are respectively $i(\alpha,\lambda)$ and $n$. From these discussions, we see that there exists a subsequence $\Phi_{n_j}$ converging to a meromorphic differential $\Phi_\infty$ with two second-order poles corresponding to pinching $\lambda$.    Similarly as in the proof of Lemma \ref{lem:Yn:length}, we see that for each $\alpha_i$  there exists a positive constant  $C_{\alpha_i}$ such that 
     \begin{equation*}
     	\left|\ell_{\alpha_i}(Y_{n_j})-2i(\alpha,\alpha_i)\right|\leq C_{\alpha_i}
     \end{equation*}
     for $j$ sufficiently large.
    Hence  $\{\ell_{\alpha_i}(Y_{n_j})\}_{j}$  is bounded for each $1\leq i\leq k$. Let $Y$ be an accumulation point of $\{Y_{n_j}\}$. Let $\Phi$ be the limit of $\Phi_{n_j}$. Then the horizontal foliation of $\Phi$ comprises two half-infinite cylinders corresponding to $\infty\cdot\lambda$ and a finite number of compactly supported components corresponding to $\beta$.  Moreover, the sequence $\{\HR_{n_j}\}$ of harmonic map rays contains a subsequence which converges to a harmonic stretch line $\mathscr R$ determined by  $Y$, $\alpha$, and $\Phi$.  That the horizontal foliation of $\Phi$ comprises two half-infinite cylinders corresponding to $\infty\cdot\lambda$ and a finite number of compactly supported components corresponding to $\beta$ implies that the maximally stretched lamination of $\mathscr R$ is exactly $\lambda$.  Since the vertical foliation of $\Phi_{n_j}$ is identically $\alpha$, it follows that $\mathscr R$ converges to $\alpha$ in the Thurston compactification. 
    This completes the proof.
    \end{proof}

    \begin{lemma}\label{lem:open:closure}  
        Let $X\in\T(S)$ and let $\lambda$ be a simple closed curve on $X$. Then the set $\mathscr{E}(X,\lambda)$ of endpoints of harmonic stretch lines $\HSR(X,Y)$ for $Y\in\outenv(X,\lambda)$ is an open subset of $\PMF(S)$.
    \end{lemma}
   \begin{proof}
   Recall that because $\outenv(X, \lambda)$ is a collection of rays $[X, \mu)$, we see that a projective measured foliation $\nu$ is contained in $\mathscr E(X,\lambda)$ if and only if any intersection $\outenv(X,\lambda)\cap [X,\nu)\neq\emptyset$. Now, it suffices to prove that the complement of $\mathscr E(X,\lambda)$ in $\PMF(S)$ is closed. Let $\mu_n$ be a sequence in the complement such that $\mu_n\to\mu_\infty$. Then by Proposition \ref{prop:cont:compactification1}, we see that $[X,\mu_n)$ converges to $[X,\mu_\infty)$ locally uniformly. Since $\outenv(X,\lambda)$ is an open subset of $\T(S)$ and is disjoint from the entire ray $[X,\mu_n)$ for every $n$, it follows that $\outenv(X,\lambda)$ is also disjoint from the entire ray $[X,\mu_\infty)$, proving that $\mu_\infty$ is not contained in $\mathscr E(X,\lambda)$.  Therefore, the complement of $\mathscr E(X,\lambda)$ in $\PMF(S)$ is a closed subset. Accordingly, the set $\mathscr E(X,\lambda)$ itself is an open subset of $\PMF(S)$.
    \end{proof}
    
    \begin{proof}[Proof of Proposition \ref{prop:accumulation}]
        Let $\lambda$ and $\mathscr R$ be as in Lemma \ref{lem:alpha:ray}. Let $X\in\mathscr R$. Then $\alpha\in \overline{\outenv(X,\lambda)}^{Th}$.  By Lemma \ref{lem:open:closure}, we see that the set $\mathscr{E}(X,\lambda)$ of endpoints of harmonic stretch lines $\HSR(X,Y)$ for $Y\in\outenv(X,\lambda)$ is an open subset of $\PMF(S)$.  Hence, for any measured lamination $\gamma$ disjoint from both $\alpha$ and $\lambda$, and for any sufficiently small positive constant $s>0$, the measured lamination $\alpha+s\gamma$ belongs to $\mathscr E(X,\lambda)$.  Furthermore, that $\gamma$ is disjoint from both $\alpha$ and $\lambda$ implies that $\lambda$ realizes the maximum 
    $$\max_{\xi}\frac{i(\alpha,\xi)}{i(\alpha+s\gamma,\xi)}:$$
     here the fraction is clearly at most one, but for $\lambda$ disjoint from $\gamma$, we have that the numerator and denominator of the fraction are both $i(\alpha,\lambda)$, so that the fraction is one.
    By Proposition \ref{prop:msl:bnd2}(i) and (iii), this implies that $\alpha+s\gamma\in\overline{\env(X,\alpha)}^{Th}$.
    \end{proof}

\begin{remark} 
    Using the same set of ideas, we can prove that for any measured geodesic lamination $\lambda$ intersecting $\alpha$ transversely, there exists $W$ such that $\Lambda(W,\alpha)=\lambda$. In particular, we may let $\lambda$  be a maximal measured lamination. Using Corollary \ref{cor:unique:geodesic:tree}, we will see that the corresponding envelope $\env(W,\alpha)$ is the harmonic stretch ray $[X,\alpha)$, proving that $\alpha$ is the only accumulation point of $\env(W,\alpha)$ in the Thurston boundary.  In particular, the accumulation set of an envelope $\env(X, \eta)$ depends on both $X$ and $\eta$ and so does not admit a topological description from, say, $\eta$ alone.
\end{remark}

\subsection{Extendability}\label{sec:ext:bnd}
For any $X\in\T(S)$ and $\eta\in\PMF(S)$, define left (resp. right) extendable and left (resp. right) boundary points similarly as in Definition \ref{def:extendability}.   

We then prove an analogue of Lemma~\ref{lem:extendability}.

\begin{lemma}\label{lem:extendability:tree}
Let  $X\in\T(S)$, $\eta\in\PMF(S)$, and $Z\in\env(X,\eta)$.  Then 
\begin{enumerate}
	\item $Z$ is right extendable if and only if  $\Lambda(Z,\eta)=\Lambda(X,\eta)$,
	\item  $Z$ is left extendable if and only if   $ \Lambda(X,Z)=\Lambda(X,\eta)$,
	\item  $Z$ is bi-extendable if and only if   $\Lambda(X,Z)=\Lambda(Z,\eta)=\Lambda(X,\eta)$.
\end{enumerate} 
\end{lemma}
\begin{proof}
    (i) The proof is similar to the proof of item (i) in Lemma \ref{lem:extendability} with Corollary \ref{cor:env in terms of stretch lams} replaced by Lemma \ref{lem:msl:bnd} and  Lemma \ref{lem:concatenation}(i) replaced by Lemma \ref{lem:concatenation:tree}(i). 
    
    (ii) For this second statement, we cannot immediately extend the proof of the second item of Lemma~\ref{lem:extendability}. Nevertheless, we may leverage the compact case in this argument, using that  Lemma~\ref{lem:extendability}(ii) holds near interior points of \tec space. 
    
    In particular, for any $Z\in\env(X,\eta)$, let $Z'\in[Z,\eta)$ be distinct from $Z$. Then $\Lambda(Z,\eta)=\Lambda(Z,Z')$. Moreover, since $\env(X,Z')\subsetneq \env(X,\eta)$ and $[Z,Z']\subsetneq [Z,\eta)$, the point $Z$ is left extendable with respect to $\env(X,\eta)$ if and only if it is left extendable with respect to $\env(X,Z')$.  By Lemma \ref{lem:extendability}(ii), the latter happens if and only if $\Lambda(X,Z)=\Lambda(X,Z')$. Now 
     \begin{eqnarray*}
       && \Lambda(X,Z') \\
&=&\Lambda(X,Z)\cap\Lambda(Z,Z') \qquad (\text{ since } [X,Z]\cup [Z,Z'] \text{ is a geodesic})\\
        &=& \Lambda(X,Z)\cap\Lambda(Z,\eta) \qquad (\text{ Since }\Lambda(Z,\eta)=\Lambda(Z,Z') )\\&=& \Lambda(X,\eta)  \qquad(\text{ by Lemma \ref{lem:msl:bnd}}).
     \end{eqnarray*}
     In summary, the point $Z$ is left extendable with respect to $\env(X,\eta)$ if and only if $\Lambda(X,Z)=\Lambda(X,\eta)$.  This completes the proof.

    (iii) This follows from (i) and (ii).
\end{proof}
\begin{remark}
	One may wonder whether the characterization in Lemma \ref{lem:extendability:tree} extends to the points at infinity or not. However,  there is an obvious obstruction in this regard, that is, the notion of maximally stretched laminations between trees is not well-defined.
\end{remark}

 With this characterization in place, we 
 may extend Proposition \ref{prop:out-in-envelopes} and Corollary \ref{cor:unique:geodesic} to the current setting. For any $\eta\in\PMF(S)$ and any chain recurrent geodesic lamination $\lambda$, define $\inenv(\eta,\lambda)$ to be the union of harmonic stretch lines which converge to $\eta$ in the forward direction and which maximally stretch exactly $\lambda$. 
 
 \begin{proposition}\label{prop:out-in-envelopes:tree}
  For any $X\in\T(S)$ and $\eta\in\PMF(S)$, we have
  \begin{equation*}
      \env(X,\eta)=\overline{\outenv(X,\Lambda(X,\eta))\cap \inenv(\eta,\Lambda(X,\eta))}.
  \end{equation*}
\end{proposition}
\begin{proof}
Let $\lambda:=\Lambda(X,\eta)$.  For any $Z\in\outenv(X,\lambda)\cap\inenv(\eta,\lambda)$, the union $[X,Z]\cup[Z,\eta)$ is a geodesic which converges to $\eta$ in the Thurston compactification, proving   that $Z\in\env(X,\eta)$. That the envelope $\env(X,\eta)$ is a closed subset of $\T(S)$ (see Lemma \ref{lem:env:tree} (ii)) then implies  \begin{equation*}
   \overline{\outenv(X,\Lambda(X,\eta))\cap \inenv(\eta,\Lambda(X,\eta))}  \quad \subset \quad \env(X,\eta).
  \end{equation*}
  The converse direction follows exactly the same argument as that of  Proposition \ref{prop:out-in-envelopes} with Lemma \ref{lem:extendability} replaced by Lemma \ref{lem:extendability:tree}.	
\end{proof}
   \begin{corollary}\label{cor:unique:geodesic:tree}
	Let $X\in\T(S)$ and $\eta\in\PMF(S)$. Then $\env(X,\eta)$ is a geodesic if and only if $\Lambda(X,\eta)$ is maximally CR.
\end{corollary}
\begin{proof}
Let $\lambda:=\Lambda(X,\eta)$.
If $\lambda$ is maximally CR, then the set $\cone(\lambda)$, of  straightened geodesic laminations (see Definition \ref{def:ad:lamination}) complementary to $\lambda$, is the singleton $\{0\}$. By Theorem \ref{thm:structure:out:in},  the union $\{X\}\cup\outenv(X,\lambda)$ is a harmonic stretch ray starting at $X$. Combined with Proposition \ref{prop:out-in-envelopes:tree} and the fact $[X,\eta)\subset \env(X,\eta)$, this implies that   $\env(X,\eta)=[X,\eta)$.  

For the  converse, suppose that $\env(X,\eta)$ is a geodesic. Then for any $Z\in\env(X,\eta)$, the envelope $\env(X,Z)$, as a subset of $\env(X,\eta)$, is also a geodesic. By Corollary \ref{cor:unique:geodesic}, we see that $\Lambda(X,Z)$ is maximally CR. Since $\Lambda(X,\eta)=\Lambda(X,Z)$, it follows that $\Lambda(X,\eta)$ is maximally CR.
\end{proof}

Regarding the shape of $\env(X,\eta)$, we have the following. 

 \begin{theorem}\label{thm:env:shape:tree} 
   For any $X\in\T(S)$ and any $\eta\in\PMF(S)$, the envelope $\env(X,\eta)$ is a cone over $\lbd(X,\eta)$ with  the cone point $\eta$ removed.  The right boundary set $\rbd(X,\eta)$ is a cone over $\bnd(X,\eta)$ with  the cone point $\eta$ removed.
 \end{theorem}
\begin{proof}
  The statement of the theorem asserts an identification of a pair of precompact sets, the first ($\env(X,\eta)$) in the Thurston compactification, and the second a deleted cone over a (not necessarily compact) set ($\lbd(X,\eta)$) in $\T(S)$. 

With this in mind, we observe that the first part almost follows from the identical argument for Theorem~\ref{thm:env:shape}, where we first substitute Proposition~\ref{prop:out-in-envelopes:tree} for Proposition~\ref{prop:out-in-envelopes} and Corollary~\ref{cor:unique:geodesic:tree} for Corollary~\ref{cor:unique:geodesic}: the only issue is the use, in the proof of Lemma~\ref{lem:env:top}, of a parametrization of $\rbd(X,Y)$ by Thurston distances in equation \eqref{eq:BdToRbdHomeo}. (In the present setting since the distances to the point $\eta$ on the Thurston boundary are now infinite, the right hand side of equation \eqref{eq:BdToRbdHomeo} becomes identically infinite upon substituting $\eta$ for the terminal point $Y$.) To adapt that argument to the present situation, we replace the term $t \cdot d_{Th}(Z,Y)$ on the right-hand side of equation \eqref{eq:BdToRbdHomeo} by a function $\psi(t)$, where $\psi: [0,1] \to [0, \infty]$ is some increasing homeomorphism between the closed intervals $[0,1]$ and $[0, \infty]$. The rest of the arguments extend verbatim.
 \end{proof}
 
\begin{corollary}\label{cor:bnd:description:pmf}
     For any $X\in\T(S)$ and any $\eta\in\PMF(S)$,  the space $\bnd(X,\eta)$ is homeomorphic to the union $\cup_{\lambda}\cone(\lambda)$, where $\lambda$ ranges over all chain recurrent geodesic laminations strictly containing $\Lambda(X,\eta)$.
\end{corollary}
\begin{proof}
   This follows from the identical argument for the proof of Corollary~\ref{cor:bnd:description} (see Section \ref{sec:proof:continuity}), once we observe that that argument only relies on geodesic segments $[X,Z]$ which originate at $X$, i.e. in the portion of the argument that characterizes $\bnd(X,Y)$, there is no dependence on the terminal point $Y$.
\end{proof}
\begin{remark}
   (i) Unlike envelopes with endpoints in $\T(S)$, the envelope $\env(X,\eta)$ might not (necessarily) be a cone over $\rbd(X,\eta)$; also, $\lbd(X,\eta)$ might not (necessarily) be a cone over $\bnd(X,\eta)$.  Here is a typical counterexample. For any $X$ and $\eta$ with $\Lambda(X,\eta)$ being maximally CR, then by Corollary \ref{cor:unique:geodesic:tree}, we have $\env(X,\eta)=[X,\eta)$. Hence $\rbd(X,\eta)=\bnd(X,\eta)=\emptyset$ but $\lbd(X,\eta)=\{X\}$.  In particular, $\env(X,\eta)$ is not a cone over $\rbd(X,\eta)$ and $\lbd(X,\eta)$ is not be a cone over $\bnd(X,\eta)$. 
   As another example, by Proposition \ref{prop:accumulation}, there exists $X\in\T(S)$ and $\eta,\mu\in\PMF(S)$ with $\eta\neq\mu$ such that $[X,\mu)\subset\env(X,\eta)$. Any point in $[X,\mu)$ is right extendable (since the whole ray $[X, \mu)$ is included in the envelope), so the ray does not meet $\rbd(X, \eta)$ and hence the envelope $\env(X,\eta)$ is not a cone over $\rbd(X,\eta)$. 

   (ii) On the other hand, using Lemma \ref{lem:env:mu} and Lemma \ref{rmk:env:accumulation}, we see that $\overline{\env(X,\eta)}^{Th}$ is a cone over $\rbd(X,\eta)\cup \left(\overline{\env(X,\eta)}^{Th}\cap\PMF(S)\right)$. Here we may view the union as the \enquote{right boundary} of $\overline{\env(X,\eta)}^{Th}$. Similarly, the \enquote{left boundary} of $\overline{\env(X,\eta)}^{Th}$ should be the union of $\lbd(X,\eta)$ with the set of endpoints of  $\overline{\env(X,\eta)}^{Th}\cap\PMF(S)$ (as a star over $\eta$). The $\bnd$ analogy of  $\overline{\env(X,\eta)}^{Th}$ should be the union of $\bnd(X,\eta)$ with the set of endpoints of  $\overline{\env(X,\eta)}^{Th}\cap\PMF(S)$ (as a star over $\eta$). With these modifications, we may extend Theorem \ref{thm:env:shape}(iii) and (iv) to $\overline{\env(X,\eta)}^{Th}$.
\end{remark}

 \subsection{Continuity of envelopes} \label{sec:continuityThurstonEnvelopes} 
 Recall (Section \ref{sec:continuity:out:in:cpt}) that for $X\in\T(S)$ and  $X\in \Omega\subset\T(S) $,  the \emph{star} $\str(\Omega,X)$ of $\Omega$ centered at $X$ is the subset of points $Z\in\Omega$ such that the harmonic stretch segment $[X,Z]$ is also contained in $\Omega$:
\begin{equation*}
    \str(\Omega,X):=\{Z\in\Omega:[X,Z]\subset \Omega\}.
\end{equation*} 
 
 Using the discussion in the previous subsection and applying the same argument as in the proof of Theorem \ref{thm:continuity}, we have \begin{theorem}\label{thm:continuity2} 
  Let $X,X_n\in\T(S)$ and $\eta_n,\eta\in\PMF(S)$ be such that $X_n\to X$ and $\eta_n\to \eta$.  Then for any closed ball $\mathscr K\subset \T(S)$ centered at $X$ of positive radius, $$\str\left(\env(X_n,\eta_n)\cap \mathscr K,X_n\right) \to \str\left(\env(X,\eta)\cap \mathscr K,X\right)$$ in the  Hausdorff topology, as $n\to\infty$. Furthermore, if $Y_n\in\T(S)$ converges to $ \eta$ in the Thurston compactification, then we also have $\str(\env(X_n,Y_n)\cap \mathscr K,X_n) \to \str(\env(X,\eta)\cap \mathscr K,X)$ in the  Hausdorff topology.
 \end{theorem}
 \begin{proof}
     The proof is a modification of that of Theorem \ref{thm:continuity}. The key to the proof of Theorem \ref{thm:continuity} is Lemma \ref{lem:approximate}, which uses Lemma \ref{lem:sequence:epsilon}, Lemma \ref{lem:length:kappa}, Proposition \ref{prop:out-in-envelopes} and Corollary \ref{cor:unique:geodesic}.  Here, we use Lemma \ref{lem:sequence:epsilon} (as it is concerned with  out-envelopes, it continues to apply in the present case), Lemma \ref{lem:length:kappa},  Proposition \ref{prop:out-in-envelopes:tree} and Corollary \ref{cor:unique:geodesic:tree},  and then apply the same algorithm as in the proof of Lemma \ref{lem:approximate}. This enables us to prove an analog of Lemma \ref{lem:approximate}.  
 \end{proof}

\appendix \section{Proof of Lemma \ref{lem:image:cut2}}
In this appendix, we provide a proof of Lemma \ref{lem:image:cut2} using train tracks (see Section \ref{sec:traintrack}). 

 \begin{proof}[Proof of Lemma \ref{lem:image:cut2}]
 Let $\lambda_0\subset\lambda$ be the stump of $\lambda$, that is, the maximal sublamination which admits a transverse measure of full support. Let $\tau_0 \subset N_\epsilon(\lambda_0)$ be a train track approximation of $\lambda_0$.  Let $\gamma\in\far_\lambda$.  Consider the set  $\{\eta_j\}$ of isolated leaves of $\gamma$. For each $\eta_j$, we add a branch $b_j$ to $\tau_0$ in such a way that  $b_j$ fellow travels $\eta_j$.    Let $\tau$ be the resulting train track.     Without loss of generality, we may assume that every switch of $\tau$ has exactly  three incident   half-branches.  Let $\zeta\subset X\setminus N_\epsilon(\lambda)$ be a train track approximation of the maximal compact sublamination of $\gamma$ (recall that $\gamma\subset X\backslash\lambda$).   Since $\lambda\cup\gamma$ is a geodesic lamination on $X$, we may assume that $\tau$ and $\zeta$ are disjoint. 
 Notice that the transverse measure on $\gamma$ induces a weight $\mathbf{w}(\eta_i)$ on each $b_i$, as well as a weight system $\mathbf{W}_{\mathrm{cpt}}(\gamma)$ on  $\zeta$.

 (1) Suppose $\gamma\in\cone(\lambda)$. Let $\hat\tau_0$ be an oriented connected component of $\tau_0$. Then there exists a measured lamination $\beta\in U_\delta(\lambda)$  carried by  $\tau\cup \zeta$  such that   $\cut_\lambda(\beta)=\gamma$ (see Equation \eqref{eq:far} for  $U_\delta(\lambda)$ and $\cut_\lambda$). The transverse measure on $\beta$ induces a weight on each branch of $\tau\cup \xi$. Consider the set of  half-branches of $\tau$ which are  incident to some switch of $\hat\tau_0$. The orientation of $\hat\tau_0$ divides these half-branches into two types with respect to their incident switches: incoming or outgoing. The switch condition at each switch ensures that the total sum of  weights of incoming  half-branches is the same as   the  total sum of  weights of outgoing  half-branches.  
 Every branch of $\hat\tau_0$  is divided into an incoming half-branch and an outgoing half-branch with the same weight. As a consequence, the sum of weights of incoming half-branches not contained in $\hat\tau_0$, is the same as the sum of weights of outgoing half-branches not contained in $\hat\tau_0$. 
 
 (2) We now turn to the inverse direction. Suppose that for any orientable component of the stump $\lambda_0$ of $\lambda$, the total measures of leaves of $\gamma$ tending to that component in the two directions are equal.  Recall that the train track $\tau$ is obtained from $\tau_0$, an approximation of  $\lambda_0$, by adding a branch $b_i$ for every isolated leaf $\eta_i$ of $\gamma$, and that $\zeta$ is a train track approximation of the (maximal) compact sublamination of $\gamma$.
  Notice that the transverse measure on $\gamma$ induces a weight system $\mathbf{W}_{\mathrm{cpt}}(\gamma)$ on $\zeta$ coming from the maximal compact sublamination of $\gamma$, and a weight $\mathbf{w}(b_i)$ on each added branch $b_i$ which equals the intersection $i(\gamma,\kappa)$ for any arc $\kappa$ which is transverse to $\gamma$ in $X\backslash\lambda$ and intersects $\eta_i$ exactly once and is disjoint from $\gamma\backslash\eta_i$.     To prove the lemma, it suffices to construct
  a measured lamination $\hat\gamma$ on $X$ such that $\cut_\lambda(\hat{\gamma})=\gamma$. To this end, since $\cut_\lambda$ keeps any compact measured lamination in $X\backslash\lambda$ invariant, we need to construct a measured lamination $\hat\gamma'$ on $X$ such that $\cut_\lambda(\hat{\gamma}')$ equals the measured sublamination of $\gamma$ obtained by removing its compact part. Recall that $\cut_\lambda$ \enquote{removes} everything near $\lambda$. Accordingly, to construct $\hat\gamma'$, it is equivalent to construct
  a weight system $w$ on $\tau$ satisfying the switch condition such that $w(b_i)=\mathbf{w}(b_i)$ for any isolated leaf $\eta_i$ of $\gamma$. In the following, we shall construct such a weight system for each component of $\tau_0$.  The construction is divided into three steps.

 \textbf{Step 1: non-orientable components of $\tau_0$.} We start with non-orientable components of $\tau_0$. Let $\hat\tau_0$ be such a component. We divide each branch   $b_i$  into half-branches. Consider those half-branches incident to some switch of $\hat\tau_0$. Let $\bar b_i\subset b_i$ be such a half-branch with the induced weight $\mathbf{w}(b_i)$. Recall that we assume that every switch of $\tau$ has exactly three incident half-branches. Let $v$  be the switch on which $\bar b_i$ is incident.  Choose a consistent orientation for the three half-branches incident on $v$ so that $\bar b_i$ is an incoming half-branch with respect to the switch $v$. 
 Since $\hat\tau_0$ is non-orientable, it follows that $\hat\tau_0\cup \bar b_i$ is alo non-orientable. By the second fact in \cite[Proposition 1.3.7]{PennerHarer1992}, there exists a train path $L$ on $\hat \tau_0\cup \bar b_i$ which starts at $\bar b_i$ and ends at the same half-branch but with the opposite orientation.
 Now we construct a modified weight system $\mathbf{W}(\bar b_i)$ on $\tau$  which associates to the half-branch $\bar b_i$ the weight $\mathbf{w}(b_i)$ and which   associates to a branch $b$ of $\hat\tau_0\subset\tau$ the weight
  \begin{equation*}\label{eq:wgt1}
  	\frac{\mathbf{w}(b_i)}{2}\times \text{ the number of times that $L$ passes through $b$}. 
  \end{equation*}
  It is clear that $\mathbf{W}(\bar b_i)$ satisfies the switch condition for every switch (here the  endpoint of $\bar b_i$ other than $v$ is not viewed as a switch of $\tau$).

   \textbf{Step 2: orientable components of $\tau_0$.}   The idea is   similar to that of the non-orientable case, but the detail is much more involved. Let $\tilde\tau_0$ be an orientable component of $\tau_0$. We cut each branch $b_i$ (in $\tau\backslash\tau_0$) into half-branches. Choose an orientation of $\tilde\tau_0$. This orientation induces an orientation for every  half-branch incident to $\tilde\tau_0$. According to this orientation, we divide these half-branches into two types: incoming or outgoing. Let $\{ b_{k}^{\mathrm{in}}\}$ and $\{ b_{j}^{\mathrm{out}}\}$ be respectively the set of incoming half-branches and outgoing half-branches that are incident to $\tilde\tau_0$. The transverse measure on $\gamma$ induces a weight on each half-branch, say $\mathbf{w}( b_k^{\mathrm{in}})$ and $\mathbf{w}( b_j^{\mathrm{out}})$. Moreover, by assumption, we have
  \begin{equation}\label{eq:weight}
  	\sum_k  \mathbf{w}( b_k^{\mathrm{in}})=\sum_j \mathbf{w}( b_j^{\mathrm{out}}).
  \end{equation}

Consider the set of pairs $P_{kj}:=(b_k^{\mathrm{in}},b_j^{\mathrm{out}})$. By \eqref{eq:weight}, we may associate to each pair $P_{ij}$ a non-negative weight $\mathbf{w}(P_{ij})$ so that 
\begin{equation*}
   \mathbf{w}( b_k^{\mathrm{in}})= \sum_{j}  \mathbf{w}( P_{kj}), \qquad \mathbf{w}( b_j^{\mathrm{out}})=\sum_{k}\mathbf{w}(P_{kj}).
\end{equation*}

It is straightforward to construct the solution $\{\mathbf{w}(P_{kj})\}$. For example, the collection $\{\mathbf{w}( b_k^{\mathrm{in}})\}$ of incoming weights define a partition of an interval of length $\sum_k  \mathbf{w}( b_k^{\mathrm{in}})$, while the collection $\{\mathbf{w}( b_j^{\mathrm{out}})\}$ of outgoing weights defines a partition of an interval of length $\sum_j \mathbf{w}( b_j^{\mathrm{out}})$.  The flux condition \eqref{eq:weight} asserts that these two intervals have the same length, so identifying the two intervals provides that weight $\mathbf{w}(P_{kj})$ as the length of any intersection in this identification between $b_k^{\mathrm{in}}$ and $b_j^{\mathrm{out}}$.

  Recall that (by assumption) the component of $\tilde\tau_0$ is orientable.  By the first fact in \cite[Proposition 1.3.7]{PennerHarer1992}, we see that, for each pair $P_{kj}$, there exists a train path $L_{kj}$ on $\tilde\tau_0\cup b_k^{\mathrm{in}}\cup b_j^{\mathrm{out}}$ which starts at $b^\mathrm{in}_k$ and ends at $b_j^\mathrm{out}$ with the prescribed orientations. This path defines a modified weight system $\mathbf{W}(P_{kj})$ on $\tau$, which associates to both $b_{k}^{\mathrm{in}}$ and $b_{j}^{\mathrm{out}}$ the weights $\mathbf{w}(P_{kj})$ and which associates to a branch $b$ of $\tilde\tau_0\subset \tau$ the weight	\begin{equation*}\label{eq:wgt3}
   		\mathbf{w}(P_{kj})\times \text{ the number of times that $L_{kj}$ passes through $b$}.
   	\end{equation*}
It is clear that $\mathbf{W}(P_{kj})$ satisfies the switch condition for every switch of $\tau$ (here the  endpoints of the half-branches $ b_{k}^{\mathrm{in}}$ and $ b_{j}^{\mathrm{out}}$ outside $\tilde\tau_0$ are not viewed as  switches.

   	\textbf{Step 3: final construction.} 
   	  Let $\mathbf{W}_{\mathrm{iso}}(\gamma)$ be the weight system on $\tau$ as follows.
   	  \begin{itemize}
   	  	\item For each branch $b_j$ corresponding to an isolated leaf $\eta_j$ of $\gamma$,  the weight is defined to be $\mathbf{w}(\eta_i)$.
   	  	\item For each branch $b$ in non-orientable component $\hat \tau_0$ of $\tau_0$, the weight is defined to be the total sum of the weight on $b$ of 
   	  	\begin{equation*}
   	  		\sum_i\mathbf{W}(\bar b_i)
   	  	\end{equation*}
   	  	where $\bar b_i$ ranges over the subset of  half-branches of $\tau\setminus\tau_0$ which are incident to $\hat \tau_0$. 
       \item For each branch $b$ in orientable component $\tilde \tau_0$ of $\tau_0$, the weight is defined to be the total sum of the weight on $b$ of 
       \begin{equation*}
           \sum_{k,j}\mathbf{W}(P_{kj}),
       \end{equation*}
      where $P_{kj}$ ranges over all pairs   $(b_k^\mathrm{in},b_j^\mathrm{out})$ of half-branches of $\tau\backslash \tau_0$ which are incident to $\tilde\tau_0$.
      \end{itemize}
      
   	 From the constructions in step 1 and step 2, we see that $\mathbf{W}_{\mathrm{iso}}(\gamma)$ satisfies the switch conditions of $\tau$ and that $\cut_\lambda(\mathbf{W} _{\mathrm{iso}}(\gamma)+\mathbf{W} _{\mathrm{cpt}}(\gamma))=\gamma$, where  $\mathbf{W} _{\mathrm{cpt}}(\gamma)$ is the weight system on $\zeta$ corresponding to the transverse measure on the maximal compact sublamination of $\gamma$.  This completes the proof.
  \end{proof}


\bibliographystyle{alpha}
\bibliography{main}

\end{document}